\documentclass[a4,10pt]{article}
\usepackage{amsmath,amscd,amssymb}

\newcommand{\nbigf}{\mathcal{F}}

\newcommand{\nbigh}{\mathcal{H}}

\newcommand{\nbigl}{\mathcal{L}}

\newcommand{\nbigo}{\mathcal{O}}
\newcommand{\nbigp}{\mathcal{P}}

\newcommand{\nbigv}{\mathcal{V}}

\newcommand{\proj}{\mathbb{P}}
\newcommand{\seisuu}{{\mathbb Z}}

\newcommand{\cnum}{{\mathbb C}}
\newcommand{\real}{{\mathbb R}}

\newcommand{\gbige}{\mathfrak E}

\newcommand{\gbigs}{\mathfrak S}

\newcommand{\veca}{{\boldsymbol a}}

\newcommand{\vecx}{{\boldsymbol x}}

\newcommand{\lrarr}{\longrightarrow}

\newcommand{\pf}{{\bf Proof}\hspace{.1in}}
\newcommand{\qed}{\mbox{\rule{1.2mm}{3mm}}}

\def\End{\mathop{\rm End}\nolimits}

\def\Image{\mathop{\rm Im}\nolimits}

\def\Re{\mathop{\rm Re}\nolimits}

\def\rank{\mathop{\rm rank}\nolimits}

\def\Tr{\mathop{\rm Tr}\nolimits}
\def\vol{\mathop{\rm vol}\nolimits}
\def\dvol{\mathop{\rm dvol}\nolimits}

\def\id{\mathop{\rm id}\nolimits}

\newcommand{\del}{\partial}
\newcommand{\delbar}{\overline{\del}}

\newcommand{\barz}{\overline{z}}
\newcommand{\zbar}{\barz}

\newcommand{\barlambda}{\overline{\lambda}}
\newcommand{\lambdabar}{\barlambda}

\newcommand{\etabar}{\overline{\eta}}
\newcommand{\xibar}{\overline{\xi}}

\newcommand{\tildepsi}{\widetilde{\psi}}
\newcommand{\psitilde}{\tildepsi}

\newcommand{\closedopen}[2]{[#1,#2[}

\newcommand{\openopen}[2]{]#1,#2[}

\newcommand{\Vtilde}{\widetilde{V}}

\newcommand{\wbar}{\overline{w}}

\newcommand{\kappatilde}{\widetilde{\kappa}}

\newcommand{\htilde}{\widetilde{h}}

\newcommand{\pitilde}{\widetilde{\pi}}

\newcommand{\ftilde}{\widetilde{f}}

\newcommand{\Utilde}{\widetilde{U}}

\def\sgn{\mathop{\rm sign}\nolimits}

\newcommand{\tautilde}{\widetilde{\tau}}

\newcommand{\Mtilde}{\widetilde{M}}

\newcommand{\Ybar}{\overline{Y}}

\newcommand{\Ktilde}{\widetilde{K}}

\newcommand{\juebar}{\bar{j}}
\newcommand{\nbigpbar}{\overline{\nbigp}}

\newcommand{\ttC}{{\tt C}}
\newtheorem{thm}{Theorem}[section]
\newtheorem{cor}[thm]{Corollary}

\newtheorem{rem}[thm]{Remark}
\newtheorem{lem}[thm]{Lemma}
\newtheorem{prop}[thm]{Proposition}
\newtheorem{df}[thm]{Definition}

\newtheorem{condition}[thm]{Condition}
\newtheorem{assumption}[thm]{Assumption}

\begin{document}

\title{Kobayashi-Hitchin correspondence for
 analytically stable bundles}

\author{Takuro Mochizuki}
\date{}
\maketitle

\begin{abstract}
We prove the existence of a Hermitian-Einstein metric
on holomorphic vector bundles
with a Hermitian metric
satisfying the analytic stability condition,
under some assumption for 
the underlying K\"ahler manifolds.
We also study the curvature decay of
the Hermitian-Einstein metrics.
It is useful for the study 
of the classification of instantons and monopoles
on the quotients of $4$-dimensional Euclidean space
by some types of closed subgroups.
We also explain examples of doubly periodic monopoles
corresponding to some algebraic data.

\vspace{.1in}


\end{abstract}

\section{Introduction}

Let $Y$ be a K\"ahler manifold
equipped with a K\"ahler form $\omega_Y$.
Let $(E,\delbar_E)$ be a holomorphic vector bundle on $Y$.
Let $A^{p,q}(E)$ denote the space of 
$C^{\infty}$-sections of $E\otimes\Omega^{p,q}$.
Let $\Lambda:A^{p,q}(E)\lrarr A^{p-1,q-1}(E)$
denote the operator obtained as the adjoint of
the multiplication of $\omega_Y$.
(See \cite{Kobayashi-vector-bundle}.)

Let $h$ be a Hermitian metric of $E$.
Let $F(h)\in A^{1,1}(\End(E))$ denote the curvature of 
the Chern connection of $E$
which is a unique unitary connection
determined by $\delbar_E$ and $h$.
The metric $h$ is called a Hermitian-Einstein metric
if the following condition is satisfied
(see \cite{Kobayashi-Nagoya, Kobayashi-vector-bundle}):
\[
 \Lambda F(h)=\frac{\Tr\Lambda F(h)}{\rank E}\id_E.
\]
In other words,
the trace free part of $\Lambda F(h)$ is $0$.
According to Kobayashi \cite{Kobayashi-Academy}
and L\"{u}bke \cite{Lubke2},
if $(Y,\omega_Y)$ is a smooth projective manifold
or more generally a compact K\"ahler manifold,
a holomorphic vector bundle
with a Hermitian-Einstein metric
is poly-stable with respect to the K\"ahler form $\omega_Y$
in the sense of Takemoto \cite{Takemoto1, Takemoto2}.
Conversely, a stable bundle on a compact K\"ahler manifold
has a Hermitian-Einstein metric,
according to the celebrated theorem of 
Donaldson \cite{Donaldson-Narasimhan-Seshadri, 
 Donaldson-surface, Donaldson-infinite}
and Uhlenbeck-Yau \cite{Uhlenbeck-Yau}.
The correspondence is called
{\em Kobayashi-Hitchin correspondence},
{\em Hitchin-Kobayashi correspondence},
or {\em Donaldson-Uhlenbeck-Yau correspondence}.
In this paper, we adopt {\em Kobayashi-Hitchin correspondence}.
Note that the $1$-dimensional case is due to
Narasimhan-Seshadri \cite{narasimhan-seshadri}.
See \cite{Hitchin-note-vanishing, Lubke-Teleman}
for more details on the history of this theorem.

It is natural to study such correspondences
on the non-compact K\"ahler manifolds.
We mention two pioneering works.
One is due to Mehta-Seshadri \cite{mehta-seshadri}
on the generalization of the theorem of Narasimhan-Seshadri
to the context of punctured Riemann surfaces.
Namely, they established an equivalence between
irreducible unitary flat bundles
and stable parabolic bundles on punctured Riemann surface.
It is one of the origins of the study of parabolic bundles,
and also an important prototype of various correspondences
for vector bundles with a good metric on quasi-projective manifolds.
The other is due to Donaldson \cite{Donaldson-GIT}
on the correspondence between
instantons on $\real^4$ and 
holomorphic bundles on $\proj^2$ 
with trivialization at infinity.
It allows us to study instantons on $\real^4$
by using the method of algebraic geometry,
which has been quite fruitful.

However, 
it is Simpson who developed 
the differential geometric analysis
for Hermitian-Einstein metrics 
on non-compact K\"ahler manifolds.
His fundamental work \cite{Simpson88} has rich contents.
Among them, he generalized
the Kobayashi-Hitchin correspondence
to several directions.
He introduced the concept of Hermitian-Einstein metrics
for Higgs bundles.
Under some assumption for the base K\"ahler manifold
which are not necessarily compact \cite[\S2]{Simpson88},
he introduced the analytic stability condition
for Higgs bundles $(E,\delbar_E,\theta)$
equipped with a Hermitian metric $h_0$.
Then, he established the following theorem.
\begin{thm}[Simpson
 \mbox{\cite[Theorem 1]{Simpson88}}]
\label{thm;19.1.2.2}
Suppose that the K\"ahler manifold $(Y,\omega_Y)$
satisfies the assumptions in
{\rm\cite[\S2]{Simpson88}}.
If $(E,\delbar_E,\theta,h_0)$ is analytically stable,
then $(E,\delbar_E,\theta)$ has a Hermitian-Einstein metric
$h$ such that
(i) $h$ and $h_0$ are mutually bounded,
(ii) $\det(h)=\det(h_0)$,
(iii) $(\delbar+\theta)(h h_0^{-1})$ is $L^2$.
\hfill\qed
\end{thm}
In particular, 
for Higgs bundles on compact K\"ahler manifolds,
the polystability condition is equivalent
to the existence of a Hermitian-Einstein metric.
He also studied the $L^2$-property of 
the curvature of the Hermitian-Einstein metrics.

He applied Theorem \ref{thm;19.1.2.2}
to the study of Higgs bundles on punctured Riemann surfaces
\cite{Simpson90}.
Let $C$ be a compact Riemann surface,
and let $D\subset C$ be a finite subset.
He established the equivalence between
tame harmonic bundles and 
polystable regular filtered Higgs bundles of degree $0$ on $(C,D)$,
as a kind of generalization of the theorem of Mehta-Seshadri.
(See \cite{Simpson90} for the terminologies. 
 See also the introduction of \cite{Mochizuki-KH-periodic}
 for a summary of Simpson's correspondence for
 tame harmonic bundles.)
It follows from the following two equivalences.
\begin{description}
\item[(a)]
 Equivalence between polystable regular filtered Higgs bundles 
 of degree $0$ on $(C,D)$,
 and {\em analytically stable} Higgs bundles of degree $0$
 on $C\setminus D$.
 Here, we use a K\"ahler metric of $C\setminus D$,
 which has orbifold-like singularity at each point of $D$.
\item[(b)]
Equivalence between tame harmonic bundles on $C\setminus D$,
and  {\em analytically stable} Higgs bundles of degree $0$
on $C\setminus D$.
\end{description}
The equivalence (a) follows 
from the careful analysis on the local property of
tame harmonic bundles around the singularity,
which was done in \cite{Simpson88, Simpson90}.
The equivalence (b) is a direct consequence 
of Theorem \ref{thm;19.1.2.2}.
It is a general guideline to obtain an equivalence between
nice objects in algebraic geometry
(polystable regular filtered Higgs bundles)
and
nice objects in differential geometry
(tame harmonic bundles)
from the studies of singularity of differential geometric objects
and the Kobayashi-Hitchin correspondence for 
{\em analytically stable} objects.

The method and the result of Simpson have been quite useful
in the study of the Kobayashi-Hitchin correspondences
for tame and wild harmonic bundles on projective manifolds.
Wild harmonic bundles on punctured Riemann surfaces
were studied by Biquard-Boalch \cite{biquard-boalch}
and Sabbah \cite{sabbah3}.
The higher dimensional generalization has been studied by
Biquard \cite{b},
Li-Narasimhan \cite{Li-Narasimhan}, Li \cite{li2},
Steer-Wren \cite{steer-wren}
and the author \cite{Mochizuki-KHI, Mochizuki-KHII, Mochizuki-wild}.

As mentioned,
in Theorem \ref{thm;19.1.2.2},
the base K\"ahler manifold should satisfy 
the conditions in \cite[\S2]{Simpson88}.
In particular, the volume should be finite.
In the study of tame and wild harmonic bundles on 
projective manifolds,
it is not so restrictive.
Indeed, because the condition for pluri-harmonic metrics 
depends only on the complex structure
of the base space, the role of K\"ahler metrics
is rather auxiliary,
and we may choose an appropriate K\"ahler metric
satisfying the conditions.

However, the condition for Hermitian-Einstein metrics
depend on the K\"ahler metrics.
There are many natural non-compact K\"ahler manifolds
such that we cannot directly apply the theorem of Simpson
to the construction of Hermitian-Einstein metrics
for vector bundles on the spaces.
For example,
we may mention the quotient space $\cnum^2/\Gamma$,
where $\Gamma\simeq \seisuu^a$  $(a<4)$
because the volume of such a space is infinite.

There are many interesting studies 
of the Hermitian-Einstein metrics 
on vector bundles over K\"ahler manifolds
with infinite volume.
For example, see
\cite{Biquard-Jardim,
Mochizuki-doubly-periodic,
Munteanu-Sesum, Ni1, Ni-Ren}
and a more recent work
\cite{Charbonneau-Hurtubise2}.
We have already mentioned 
the pioneering work of Donaldson
\cite{Donaldson-GIT}.
See also 
\cite{GuoI, GuoII, Owens}\footnote{These references
were informed by one of the reviewers.}.
We may also mention the construction of monopoles
on $\real^3$ in \cite{Jarvis-construction} as a related work.
However, the author does not find any systematic study 
on the relation between the existence of a Hermitian-Einstein metric
and {\em the analytic stability condition} 
on K\"ahler manifolds with infinite volume,
that is our issue in this paper.

We introduce a condition for K\"ahler manifolds 
(Assumption \ref{assumption;17.11.29.40})
which is weaker than the conditions
in \cite[\S2]{Simpson88}.
In particular, the volume can be infinite.
The analytic stability condition is defined in a natural way
for holomorphic vector bundles with a Hermitian metric
on such spaces.
(See Definition \ref{df;17.12.23.1}.)
Then, we shall prove the following theorem
for the existence of a Hermitian-Einstein metric
on analytically stable holomorphic vector bundles.
\begin{thm}[Theorem 
\ref{thm;17.8.8.11}]
\label{thm;19.1.2.3}
Suppose that $(Y,g_Y)$ is a K\"ahler manifold
satisfying Assumption {\rm\ref{assumption;17.11.29.40}}.
Let $(E,\delbar_E)$ be a holomorphic vector bundle on $Y$
with a Hermitian metric $h_0$
satisfying the analytic stability condition.
Then,
$(E,\delbar_E)$ has a Hermitian-Einstein metric $h$
such that
(i) $\det(h)=\det(h_0)$,
(ii) $h$ and $h_0$ are mutually bounded,
(iii) $\delbar_E(h h_0^{-1})$ is $L^2$.
\end{thm}
We shall also study the curvature decay 
(Propositions \ref{prop;17.8.10.40}--\ref{prop;17.8.12.1}),
and the uniqueness  of such metrics
(Proposition \ref{prop;17.11.30.3}).
The Higgs case will be studied elsewhere.

As in the case of \cite{mehta-seshadri, Simpson90},
when we are given a natural compactification
$\Ybar$ of a non-compact K\"ahler manifold $Y$,
it is interesting to obtain 
an equivalence between holomorphic vector bundles
with a Hermitian-Einstein metric on $Y$
and parabolic vector bundles on $(\Ybar,H_{\infty})$
satisfying a kind of stability condition,
where $H_{\infty}:=\Ybar\setminus Y$.
Ideally, it should be established
as a consequence of two equivalences like (a) and (b) above.
We expect that Theorem \ref{thm;19.1.2.3}
will be useful to study an analogue of (b).

Because
$\real\times T^3$ and $\real^2\times T^2$ satisfy 
Assumption \ref{assumption;17.11.29.40}
for the underlying K\"ahler manifolds
as explained in \S\ref{section;17.11.14.1},
where $T^j$ denotes a $j$-dimensional real torus,
Theorem \ref{thm;19.1.2.3}
is useful in the study of instantons
and monopoles on $\cnum^2/\Gamma$
for some types of $\Gamma$.
Indeed, 
the author 
has already applied it to 
the study of monopoles with Dirac type singularity
on $S^1\times\real^2$ in
\cite{Mochizuki-KH-periodic}.
We shall also explain a way to construct
examples of doubly periodic monopoles
from some algebraic data in \S\ref{subsection;18.1.20.12}.

For the proof of Theorem \ref{thm;19.1.2.3},
instead of the method of the heat equation,
we apply the deep result of Donaldson
on the Dirichlet problem for Hermitian-Einstein metrics
\cite{Donaldson-boundary-value}.
Suppose that $(E,\delbar_E,h_0)$ is analytically stable.
For simplicity, let us consider the case $\Tr F(h_0)=0$.
We take a sequence 
of closed submanifolds with boundary $Y_i\subset Y$
$(i=1,2,\ldots)$
such that $\bigcup Y_i=Y$.
According to \cite{Donaldson-boundary-value},
there exists a unique Hermitian-Einstein metric $h_i$ of $E_{|Y_i}$
satisfying
$\Lambda F(h_i)=0$ and
$h_{i|\del Y_i}=h_{0|\del Y_i}$.
Note that we have $\det(h_i)=\det(h_{0|Y_i})$.
It is natural to ask the convergence 
of a subsequence of $\{h_i\}$ on any compact subset,
which should have the desired property.
In \cite{Ni1},
Ni gave a useful argument to obtain the convergence
from a boundedness of the sequence $\{h_i\}$
by using the idea of an effective $C^1$-estimate
due to Donaldson \cite{Donaldson-boundary-value}.
(See also Remark {\rm\ref{rem;19.1.2.1}}.)
Hence, our issue is to obtain a $C^0$-bound
from the analytic stability condition.
We consider the Donaldson functional 
$M(h^{(1)},h^{(2)})$ on each $Y_i$.
It turns out that 
the inequality $M(h_{0|Y_i},h_i)\leq 0$ holds.
It allows us to deduce the desired $C^0$-bound
(Proposition \ref{prop;17.8.8.10}),
for which we essentially apply the argument of Simpson
in the proof of \cite[Proposition 5.3]{Simpson88}
by adjusting to Assumption \ref{assumption;17.11.29.40}.

\vspace{.1in}

We should remark that 
the analysis in the infinite volume case is not always
more difficult than 
the analysis in the finite volume case.
For example,
$L^2$-instantons on $\cnum^2$ with the standard Euclidean metric
is equivalent to instantons on the compact space $S^4$,
as was used in \cite{Donaldson-GIT}.

We should also remark that
the analytic stability condition does not seem
essential in some cases.
It has been known that 
the solvability of the Poisson equation on $Y$
is useful for the construction of 
Hermitian-Einstein metrics
of holomorphic vector bundles $(E,\delbar_E)$ on $Y$.
The idea has been efficiently used in
\cite{Bando, Jarvis-construction, Munteanu-Sesum, Ni1, Ni-Ren}.
In those cases,
a Hermitian-Einstein metric $h$ on $(E,\delbar_E)$
is obtained 
if $(E,\delbar_E)$ has a Hermitian metric $h_0$
such that $\Lambda F(h_0)$ satisfies
a decay condition,
even if we do not assume $(E,\delbar_E,h_0)$
is analytically stable.
Hence, the efficiency of the analytic stability condition
depends on the geometric property of $(Y,\omega_Y)$
around infinity.

Let us mention some examples.
In the case of  $\cnum^2/\seisuu$,
the volume of balls radius $r$ rapidly grows as $r\to\infty$,
and the Green function decays rapidly enough around infinity.
Therefore, it is naturally expected
to construct a Hermitian-Einstein metrics
for holomorphic bundles on $\cnum^2/\seisuu$
(or more generally ALF-spaces)
by assuming only the existence of a metric $h_0$
satisfying only the curvature decay condition,
not the analytic stability condition.
But, 
for the construction of monopoles on $\real^3/\seisuu$ and $\real^3/\seisuu^2$,
instantons on $\cnum^2/\seisuu^2$ and $\cnum^2/\seisuu^3$,
or more generally
instantons on ALG-spaces and ALH-spaces,
it seems significant to study the analytic stability condition.
Note that instantons on
$\cnum^2/\seisuu^2$ 
and $\cnum^2/\seisuu^3$
have already been studied
\cite{Biquard-Jardim, Charbonneau-Hurtubise2, Mochizuki-doubly-periodic}.
However, even if a detailed proof is explained,
they seem to depend on the specific global property of
 $Y$ and the compactification,
and some additional assumptions on the bundles.
For example, the author \cite{Mochizuki-doubly-periodic}
studied such an issue
in the case $\cnum^2/\seisuu^2$,
but with an indirect method
using the Nahm transform and
the Kobayashi-Hitchin correspondence
for wild harmonic bundles on an elliptic curve.
It looks significant to the author to give detailed arguments
on the basis of the study of analytic stable bundles.

\begin{rem}
In the proof of
{\rm\cite[Theorem 0.12, Proposition 5.12]{Biquard-Jardim}},
a generalization of the theorem of Simpson 
{\rm\cite[Theorem 1]{Simpson88}}
is mentioned on the basis of {\rm\cite{b}}.
(See Remark {\rm\ref{rem;17.12.1.1}} below.)
It is not clear how it is related with our result in this paper.
Anyway, the author thinks it useful to give a general statement
with a proof for further studies on similar issues.
(See {\rm\cite{Mochizuki-KH-periodic}}, for example.)
\hfill\qed
\end{rem}

\paragraph{Acknowledgement}

A part of this study was done during my stay in 
the University of Melbourne,
and I am grateful to Kari Vilonen and Ting Xue
for their excellent hospitality and their support.
I am grateful to Carlos Simpson
whose works provide the most important foundation
with this study.
I thank Olivier Biquard for kind replies to my questions
on his works.
I am grateful to Claude Sabbah
for his kindness and discussions on many occasions.
I thank Indranil Biswas, Hisashi Kasuya,
Jacques Hurtubise
and Masaki Yoshino for some discussions.
I thank Yoshifumi Tsuchimoto and Akira Ishii
for their constant encouragement.
I thank the referees for their helpful comments
to improve this paper.

My interest to ``Kobayashi-Hitchin correspondence''
was renewed when I made a preparation for 
a talk in the 16th Oka Symposium (2017),
which drove me to this study.
I thank the organizers,
particularly 
Junichi Matsuzawa and Ken-ichi Yoshikawa.

I am partially supported by
the Grant-in-Aid for Scientific Research (S) (No. 17H06127),
the Grant-in-Aid for Scientific Research (S) (No. 16H06335),
and the Grant-in-Aid for Scientific Research (C) (No. 15K04843),
Japan Society for the Promotion of Science.

\section{Kobayashi-Hitchin correspondence for 
analytically stable bundles}

\subsection{Assumption on the base space}
\label{subsection;17.11.29.1}

Let $G$ be a compact Lie group.
Let $(X,g_X)$ be an $n$-dimensional connected K\"ahler manifold.
Suppose that $(X,g_X)$ is equipped with
a left $G$-action $\kappa$ in the following sense.
\begin{itemize}
\item
 $\kappa:G\times X\lrarr X$ is a $C^{\infty}$-map
 satisfying
 $\kappa(a,\kappa(b,x))=\kappa(ab,x)$
 and $\kappa(1,x)=x$.
 Set $\kappa(a,x):=\kappa_a(x)$.
\item
 $\kappa_a:X\lrarr X$ is holomorphic
  for each $a\in G$.
\item
 $\kappa_a^{\ast}g_X=g_X$.
\end{itemize}

Let $\omega_X$ denote the K\"ahler form,
and let $\dvol_X$ denote the volume form 
associated to $g_X$.
We shall often denote $\int_Xf\dvol_X$
by $\int_Xf$.
We set
$\Delta_X:=-\sqrt{-1}\Lambda\del_X\delbar_X$.
Note that it is $1/2$-multiple of
the Laplacian associated to the Riemannian metric.

\begin{assumption}\mbox{{}}
\label{assumption;17.11.29.40}
There exist a $G$-invariant function
$\varphi_X:X\lrarr \openopen{0}{\infty}$ 
with $\int_X\varphi_X<\infty$,
and a function $\ttC_X:\real_{>0}\lrarr\real_{>0}$ $(i=1,2)$
satisfying the following condition.
\begin{itemize}
\item
Let $f:X\lrarr \closedopen{0}{\infty}$ be 
a bounded function such that
$\Delta_Xf\leq B\varphi_X$ for a positive number $B$
as a distribution.
Then, the following holds:
\[
 \sup_{P\in X}f(P)
\leq
 \ttC_X(B)
 \Bigl(
1+ \int_X f\varphi_X
 \Bigr)
\]
\end{itemize}
Moreover, if a bounded function $f$ satisfies 
the stronger condition $\Delta_X f\leq 0$ on $X$,
$\Delta_Xf=0$ holds.
\hfill\qed
\end{assumption}

See \S\ref{section;17.11.14.1}
for examples satisfying the assumption.

\subsection{Analytic stability condition for $G$-equivariant bundles}

Let $(E,\delbar_E)$ be a holomorphic vector bundle
on $X$ with a Hermitian metric $h_0$.
It is assumed to be $G$-equivariant in the following sense.
\begin{itemize}
\item
We are given a $C^{\infty}$-isomorphism
$\Theta:\kappa^{-1}(E)\simeq p_2^{-1}(E)$,
here $p_2$ is the projection $G\times X\lrarr X$.
Let $\Theta_a$ denote the restriction of
$\Theta$ to $\{a\}\times X$.
\item
 $\Theta_{a_1\cdot a_2}=
 \Theta_{a_2}\circ
 \kappa_{a_2}^{-1}(\Theta_{a_1})$
holds for any $a_1,a_2\in G$.
\item
$\Theta_{a}$ induces a holomorphic isomorphism
 $\kappa_a^{-1}(E,\delbar_E)\simeq (E,\delbar_E)$
and an isometry
$\kappa_a^{-1}(E,h_0)\simeq (E,h_0)$
for any $a\in G$.
\end{itemize}
Let $F(h_0)$ denote the curvature of 
the Chern connection of $(E,\delbar_E,h_0)$.
We assume the following.
\begin{itemize}
\item
We have $B>0$ such that 
$\bigl|\Lambda F(h_0)\bigr|_{h_0}
 \leq B\varphi_X$.
\end{itemize}

As in \cite{Simpson88},
we set
\[
 \deg(E,h_0):=
 \sqrt{-1}\int_X
 \Tr\bigl(
 \Lambda F(h_0)
 \bigr).
\]
Let $V$ be any $\nbigo_X$-submodule of $E$
which is saturated,
i.e., $E/V$ is torsion-free.
There exists a closed complex analytic subset $Z(V)\subset X$
such that
(i) $\dim Z(V)\leq \dim X-2$,
(ii) $V_{|X\setminus Z(V)}$ is a locally free $\nbigo$-module.
We obtain the induced metric $h_{0,V}$ of
$V_{|X\setminus Z(V)}$.
We set
\[
 \deg(V,h_{0}):=
 \sqrt{-1}\int_{X\setminus Z(V)}
 \Tr\Bigl(
 \Lambda F(h_{0,V})
 \Bigr).
\]
The integral is well defined
in $\real\cup\{-\infty\}$
by the Chern-Weil formula:
\[
 \deg(V,h_0)=
 \sqrt{-1}
 \int_{X\setminus Z(V)}
 \Tr\bigl(
 \pi_V\Lambda F(h_0)
 \bigr)
-\int_{X\setminus Z(V)}
 \bigl|
 \delbar_E\pi_V
 \bigr|^2_{h_0,g_X}.
\]
Here, $\pi_V$ denote the orthogonal projection
of $E_{|X\setminus Z(V)}$ to $V_{|X\setminus Z(V)}$,
and $|\cdot|_{h_0,g_X}$ denote the norm 
induced by $h_0$ and $g_X$.

\vspace{.1in}
An $\nbigo_X$-submodule $V$ of $E$
is called $G$-equivariant
if $\Theta_a(\kappa_a^{\ast}V)=V$
for any $a\in G$.

\begin{df}
\label{df;17.12.23.1}
$(E,\delbar_E,h_0)$ is called analytically stable
with respect to the $G$-action
if 
\[
 \frac{\deg(V,h_{0})}{\rank V}
<
 \frac{\deg(E,h_0)}{\rank E}
\]
holds for any saturated $G$-equivariant $\nbigo_X$-submodule $V$ of $E$
such that $0<\rank V<\rank E$.
We say that $(E,\delbar_E,h_0)$ is analytically semistable
with respect to the $G$-action
if 
\[
 \frac{\deg(V,h_{0})}{\rank V}
\leq
 \frac{\deg(E,h_0)}{\rank E}
\]
holds for any saturated $G$-equivariant $\nbigo_X$-submodule $V$ of $E$.
We say that $(E,\delbar_E,h_0)$ is analytically polystable
if it is analytically semistable with respect to 
the $G$-action,
and $G$-equivariantly isomorphic to 
the direct sum
$\bigoplus (E_i,\delbar_{E_i},h_{0,i})$,
where each $(E_i,\delbar_{E_i},h_{0,i})$ is 
analytically stable with respect to the $G$-action.
\hfill\qed
\end{df}

\subsection{Hermitian-Einstein metrics}

Let $(E,\delbar_E)$ be a $G$-equivariant holomorphic vector bundle on $X$
with a $G$-invariant Hermitian metric $h$.
Let $F(h)^{\bot}$ denote the trace free part of $F(h)$.
We say that $h$ is a Hermitian-Einstein metric
if  $\Lambda F(h)^{\bot}=0$,
i.e.,
$\Lambda F(h)=
 \Tr\Lambda\bigl(F(h)\bigr)\id_E/\rank E$.

\begin{lem}
Suppose that $h$ is 
a Hermitian-Einstein metric of $(E,\delbar_E)$,
and that
$|\Tr \Lambda F(h)|<B\varphi_X$ for some $B>0$.
Then, $(E,\delbar_E,h)$ is analytically polystable
with respect to the $G$-action.
\end{lem}
\pf
Let $V$ be any saturated $G$-equivariant $\nbigo_X$-submodule of $E$
such that $0<\rank V<\rank E$.
By the Chern-Weil formula,
the following holds:
\[
 \deg(V,h)
=\frac{\rank V}{\rank E}\deg(E,h)
-\int_{X\setminus Z(V)}
 \bigl|
 \delbar_E\pi_V
 \bigr|^2_{h_0,g_X}
\leq
 \frac{\rank V}{\rank E}\deg(E,h).
\]
If the equality holds,
we obtain  $\delbar \pi_V=0$
on $X\setminus Z(V)$,
i.e.,
$\pi_V$ is holomorphic.
By the Hartogs theorem,
there exists the holomorphic endomorphism
$\pitilde_V$ of $E$ such that
$\pitilde_{V|X\setminus Z(V)}=\pi_V$.
Because $\pi_V$ is holomorphic and an orthogonal projection,
it satisfies $\del_{E,h}\pi_V=0$.
By the continuity, 
$\pitilde_V$ is also an orthogonal projection
and satisfies 
$\delbar_E\pitilde_V=0$ and $\del_{E,h}\pitilde_V=0$.
Let $V^{\bot}$ denote the orthogonal complement of
$V$ in $E$ with respect to $h$.
It is also a $G$-equivariant holomorphic subbundle.
Let $h_V$ and $h_{V^{\bot}}$ denote
the induced Hermitian metrics of $V$
and $V^{\bot}$, respectively.
We can easily observe that 
$h_V$ and $h_{V^{\bot}}$
are $G$-invariant Hermitian-Einstein metrics.
Then, we obtain the claim of the lemma
by an easy induction on the rank of the bundle.
\hfill\qed

\begin{prop}
\label{prop;17.11.30.5}
Let $h_1$ and $h_2$ be $G$-equivariant
Hermitian-Einstein metrics
of $(E,\delbar_E)$ satisfying the following.
\begin{itemize}
\item
 $h_1$ and $h_2$ are mutually bounded.
\item
 $\Lambda F(h_1)=\Lambda F(h_2)$.
\end{itemize}
Then, there exists a $G$-invariant holomorphic decomposition 
$E=\bigoplus_{i=1}^m E_i$
and a tuple 
$(c_1,\ldots,c_m)\in\real_{>0}^m$
such that the following holds:
\begin{itemize}
\item
 The decomposition $E=\bigoplus_{i=1}^m E_i$
 is orthogonal with respect to
 both $h_i$ $(i=1,2)$.
\item
 $h_{1|E_i}=c_i\cdot h_{2|E_i}$.
\end{itemize}
\end{prop}
\pf
Let $b$ be the automorphism of $E$
determined by $h_1=h_2b$,
which is $G$-invariant.
According to \cite[Lemma 3.1]{Simpson88},
the following holds:
\[
 \Delta_X\Tr(b)
=
-\bigl|
 \delbar_E(b)b^{-1/2}
 \bigr|_{h_2,g_X}^2\leq 0.
\]
By Assumption \ref{assumption;17.11.29.40},
we obtain that
$\Delta_X\Tr(b)=0$.
Hence, we obtain
$\bigl|
 \delbar_E(b)b^{-1/2}
 \bigr|_{h_2,g_X}=0$.
It implies that 
$\delbar_E(b)=0$.
Because $b$ is self-adjoint with respect to
$h_i$ $(i=1,2)$,
we also obtain $\del_{E,h_i}(b)=0$,
i.e.,
$b$ is flat with respect to the Chern connections of
$(E,\delbar_E,h_i)$.
In particular, the eigenvalues of $b$ are constant.
Let $E=\bigoplus_{i=1}^m E_i$
denote the eigen decomposition of $b$.
It satisfies the condition desired in the lemma.
\hfill\qed

\subsection{Statements}

\label{subsection;17.11.27.4}

\subsubsection{Existence of Hermitian-Einstein metrics and 
the analytic stability condition}
\label{subsection;18.1.19.20}

Let $(E,\delbar_E)$ be a $G$-equivariant holomorphic vector bundle on $X$.
Let $h_0$ be a $G$-invariant Hermitian metric of $E$
such that 
$|\Lambda F(h_0)|_{h_0}\leq B\varphi_X$ for $B>0$.
We shall prove the following theorem in
\S\ref{subsection;17.11.15.1}--\ref{subsection;17.11.27.1}.

\begin{thm}
\label{thm;17.8.8.11}
If $(E,\delbar_E,h_0)$ is analytically stable
with respect to the $G$-action,
there exists a $G$-invariant Hermitian-Einstein metric $h$
of $(E,\delbar_E)$ satisfying the following conditions.
\begin{itemize}
\item
$\det(h)=\det(h_0)$,
 where $\det(h)$ and $\det(h_0)$ denote
 the Hermitian metrics on 
 the determinant line bundle $\det(E)$
 induced by $h$ and $h_0$, respectively.
\item Let $b$ be the automorphism of $E$
determined by $h=h_0b$.
Then, $|b|_{h_0}$ and $|b^{-1}|_{h_0}$ are bounded,
and $\int_X|\delbar b|_{h_0,g_X}^2<\infty$.
\end{itemize}
\end{thm}

\subsubsection{Complement on curvature decay}
\label{subsection;17.8.10.50}

Let $Y$ be a $G$-invariant closed end of $X$
satisfying the following conditions.
We study the behaviour of the Hermitian-Einstein metric 
in Theorem \ref{thm;17.8.8.11} on $Y$.

\begin{condition}
\label{condition;18.12.23.10}
\mbox{{}}
\begin{itemize}
\item
The curvature of $g_{X|Y}$ is bounded.
\item
There exist $r_0>0$ and a compact subset
$Y_0\subset Y$ such that
the injective radius of any point $P\in Y\setminus Y_0$ in $Y$
is larger than $r_0$.
\end{itemize}
\end{condition}

Let $(E,\delbar_E,h_0)$ 
and $h$ be as in 
Theorem \ref{thm;17.8.8.11}.
We shall prove the following proposition
in \S\ref{subsection;17.11.30.1}.
\begin{prop}
\label{prop;17.8.10.40}
Assume Condition {\rm\ref{condition;18.12.23.10}}
on $(Y,g_{X|Y})$.
Moreover, we assume the following additional condition
on $(E,\delbar_E,h_0)_{|Y}$.
\begin{itemize}
\item
For any $\epsilon>0$,
there exists a compact subset $K(\epsilon)$ of $Y$
such that
\[
|F(h_0)^{\bot}|_{h_0,g_X}\leq \epsilon,
\quad
 \bigl|
 \delbar_{E}\Lambda F(h_0)^{\bot}
 \bigr|_{h_0,g_X}\leq\epsilon
\]
on $Y\setminus K(\epsilon)$.
Note that it also implies
$\bigl|
 \del_{E,h_0}\Lambda F(h_0)^{\bot}
 \bigr|_{h_0,g_X}\leq\epsilon$
on $Y\setminus K(\epsilon)$
\end{itemize}
Then, for any $\epsilon>0$,
there exists a compact subset $K(\epsilon)\subset Y$
such that
$|F(h)^{\bot}|_{h,g_X}\leq\epsilon$ on $Y\setminus K(\epsilon)$.
\end{prop}

The following proposition is similar
and proved in \S\ref{subsection;17.11.30.1}.
\begin{prop}
\label{prop;18.1.20.30}
Assume Condition {\rm\ref{condition;18.12.23.10}}
on $(Y,g_{X|Y})$.
Moreover, we assume the following additional condition
on $(E,\delbar_E,h_0)_{|Y}$.
\begin{itemize}
\item
$|F(h_0)^{\bot}|_{h_0,g_X}$
and
$\bigl|
 \delbar_{E}\Lambda F(h_0)^{\bot}
 \bigr|_{h_0,g_X}$
are bounded on $Y$.
Note that it also implies 
the boundedness of
$\bigl|
 \del_{E,h_0}\Lambda F(h_0)^{\bot}
 \bigr|_{h_0,g_X}$ on $Y$.
\end{itemize}
Then, 
$|F(h)^{\bot}|_{h,g_X}$  is bounded on $Y$.
\end{prop}

\subsubsection{Complement on the $L^2$-property of the curvature}
\label{subsection;17.11.29.11}

We give some sufficient conditions
for the Hermitian-Einstein metric in Theorem {\rm \ref{thm;17.8.8.11}}
to be $L^2$.
We shall prove the following proposition
in \S\ref{subsection;17.11.30.2}.
\begin{prop}
\label{prop;17.8.12.1}
Let $(E,\delbar_E,h_0)$ be as in 
Theorem {\rm\ref{thm;17.8.8.11}}.
We assume
$\int_X\bigl|F(h_0)^{\bot}\bigr|^2_{h_0,g_X}<\infty$.
Moreover, we assume that
there exists a $G$-invariant  exhaustion function $\varrho$ on $X$
satisfying the following conditions.
\begin{itemize}
\item
$\del\delbar\varrho$ is bounded.
\item
$\lim_{t\to\infty}
 \frac{1}{t}
 \int_{\{\varrho\leq t\}}
 \bigl|
 F(h_0)^{\bot}\cdot\del\delbar\varrho
 \bigr|_{h_0,g_X}
 =0$.
\item
 $\lim_{t\to\infty}
 \frac{1}{t}
  \Bigl(
 \int_{\{\varrho\leq t\}}
 \bigl|
 F(h_0)^{\bot}\cdot\delbar\varrho
 \bigr|^2_{h_0,g_X}
 \Bigr)^{1/2}=0$.
\end{itemize}
Let $h$ be the $G$-invariant Hermitian-Einstein metric
of $(E,\delbar_E)$ as in Theorem {\rm \ref{thm;17.8.8.11}}.
Then, the following holds:
\begin{equation}
 \label{eq;18.1.21.1}
 \int_X\Tr\bigl((F(h)^{\bot})^2\bigr)\cdot\omega_X^{n-2}
=\int_X\Tr\bigl((F(h_0)^{\bot})^2\bigr)\cdot\omega_X^{n-2}.
\end{equation}
In particular,
$F(h)^{\bot}$ is $L^2$.
\end{prop}

\begin{cor}
\label{cor;17.11.29.20}
Suppose that there exists a $G$-invariant 
exhaustion function $\varrho$ on $X$
such that
(i) $\del\delbar\varrho$ is $L^2$ and bounded,
(ii) $\del\varrho$ is bounded.
Let $(E,\delbar_E)$ be a holomorphic vector bundle on $X$
with a Hermitian metric $h_0$
such that 
(a) $|\Lambda F(h_0)|_{h_0}\leq B\varphi_X$ for some $B>0$,
(b) $\Tr F(h_0)=0$,
(c) $F(h_0)$ is $L^2$.
If $(E,\delbar_E,h_0)$ is analytically stable
with respect to the $G$-action,
then there exists a $G$-invariant Hermitian metric $h$ of $E$
such that
(i) $h$ and $h_0$ are mutually bounded,
(ii) $\det(h)=\det(h_0)$,
(iii) $\Lambda F(h)=0$,
(iv) $F(h)$ and $\delbar_E(h h_0^{-1})$ are $L^2$.
Moreover, 
the equality {\rm(\ref{eq;18.1.21.1})} holds.
\end{cor}
\pf
Let $h$ be the $G$-invariant
Hermitian-Einstein metric for $(E,\delbar_E)$
in Theorem {\rm \ref{thm;17.8.8.11}}.
Because $\Tr F(h_0)$ and $\det(h)=\det(h_0)$,
we obtain $\Tr F(h)=0$.
In particular, $\Lambda F(h)=0$ holds.
Because the assumptions 
in Proposition \ref{prop;17.8.12.1}
are satisfied for $(E,\delbar_E,h_0)$ with $\varrho$,
we obtain that $F(h)$ is $L^2$.
\hfill\qed

\begin{rem}
Proposition {\rm\ref{prop;17.8.12.1}}
and Corollary {\rm\ref{cor;17.11.29.20}}
are variants of 
{\rm\cite[Proposition 3.5, Lemma 7.4]{Simpson88}}.
\hfill\qed
\end{rem}

\subsubsection{Uniqueness}

Let us give a sufficient condition
for the uniqueness of metrics with the properties
in Theorem \ref{thm;17.8.8.11}.
We shall prove the following proposition
in \S\ref{subsection;17.11.30.4}.

\begin{prop}
\label{prop;17.11.30.3}
Suppose that there exists an exhaustion function
$\varrho_1:X\lrarr \real_{>0}$ 
such that
$\delbar\log\varrho_1$ is $L^2$ on $X$.
Let $(E,\delbar_E,h_0)$ be an analytically stable bundle
on $X$ with respect to the $G$-action
as in {\rm\S\ref{subsection;18.1.19.20}}.
Suppose that $h_i$ $(i=1,2)$ are 
$G$-invariant Hermitian-Einstein metrics
such that
(i) $\det(h_i)=\det(h_0)$,
(ii) $h_i$ and $h_0$ are mutually bounded.
Then, $h_1=h_2$ holds.
\end{prop}

Let us mention why the uniqueness of Hermitian-Einstein metrics
does not follow immediately from Proposition \ref{prop;17.11.30.5}.
It is not clear to the author whether
the analytic stability of $(E,\delbar_E,h_0)$
implies the analytic stability of $(E,\delbar_E,h)$
without any additional condition,
where $h$ is a Hermitian-Einstein metric as 
in Theorem \ref{thm;17.8.8.11}.
This issue does not appear
in the case where $X$ is compact.
Namely, if $X$ is compact,
the analytic stability condition of $(E,\delbar_E,h_0)$
is equivalent to the stability condition of $(E,\delbar_E)$,
and the latter is equivalent to
the analytic stability condition of $(E,\delbar_E,h)$.
As a result,
the uniqueness of Hermitian-Einstein metrics
is an easy consequence of Proposition \ref{prop;17.11.30.5}
in the compact case.
However, in the non-compact case,
we need an additional work to prove 
the analytic stability of $(E,\delbar_E,h)$
from 
the analytic stability of $(E,\delbar_E,h_0)$.

\begin{rem}
It would be instructive to state explicitly that
a Hermitian-Einstein metric
in {\rm\cite[Theorem 1]{Simpson88}} is unique
without any additional assumption.
Indeed, 
suppose that 
there exist Hermitian-Einstein metrics $h_i$ $(i=1,2)$
under the assumption in {\rm\cite{Simpson88}},
such that 
(i) $h_i$ are mutually bounded with $h_0$,
(ii) $\det(h_i)=\det(h_0)$.
By an argument in Proposition {\rm\ref{prop;17.11.30.5}},
we obtain the holomorphic decomposition
$(E,\delbar_E)=\bigoplus_{i=1}^m(E_i,\delbar_{E_i})$,
which are orthogonal with respect to both $h_i$ $(i=1,2)$.
For $\veca=(a_1,\ldots,a_m)\in\real_{>0}^m$ with
$\prod a_i^{\rank E_i}=1$,
we consider the automorphism
$F_{\veca}=
 \bigoplus_{i=1}^m a_i\id_{E_i}$
of $E$.
Let $s_{\veca}$ be determined by
$F_{\veca}^{\ast}h_1=h_0e^{s_{\veca}}$.
According to {\rm\cite[Proposition 5.3]{Simpson88}},
there exist positive constants $C_i$ $(i=1,2)$
such that
$\sup_X|s_{\veca}|_{h_0}
\leq
C_1+C_2M(h_0,F_{\veca}^{\ast}h_1)$.
Here, 
$M(h_0,F_{\veca}^{\ast}h_1)$
denotes the Donaldson functional
in {\rm\cite[\S5]{Simpson88}}.
By {\rm\cite[Proposition 5.1]{Simpson88}},
$M(h_0,F_{\veca}^{\ast}h_1)
=M(h_0,h_1)+M(h_1,F_{\veca}^{\ast}h_1)$
holds.
Because $F_{\veca}^{\ast}h_1$ are Hermitian-Einstein
with $\det(F_{\veca}^{\ast}h_1)=\det(h_1)$,
$\delbar_E(F_{\veca}^{\ast}(h_1)h_1^{-1})=0$ holds
as observed in the proof of Proposition {\rm\ref{prop;17.11.30.5}}.
Hence, we obtain $M(h_1,F_{\veca}^{\ast}h_1)=0$
by the definition of
Donaldson functional in {\rm\cite[\S5]{Simpson88}}.
We obtain 
$\sup_X|s_{\veca}|_{h_0}
\leq
C_1+C_2M(h_0,h_1)$
for any $\veca$.
If $m\geq 2$,
there exists a sequence $\veca_i=(a_{i,1},\ldots,a_{i,m})$
such that $a_{i,1}\to\infty$,
which implies
$\sup_X|s_{\veca_i}|_{h_0}\to\infty$.
Hence, we obtain $m=1$.
\hfill\qed
\end{rem}

\subsection{Review of the Dirichlet Problem for Hermitian-Einstein metrics}

Let $(Z,g_Z)$ be a connected compact K\"ahler manifold
with a non-empty smooth boundary $\del Z$
equipped with an action of a compact Lie group $G$.
For simplicity, we assume that $Z$ is embedded to
a K\"ahler manifold $(Z',g_{Z'})$ equipped with a $G$-action
such  that
$Z\setminus \del Z$
is a $G$-invariant open subset of $Z'$.

Let $(E',\delbar_{E'})$ be a $G$-equivariant holomorphic vector bundle
on a neighbourhood of $Z$.
Let $(E,\delbar_E)$ be the restriction of $(E',\delbar_{E'})$
to $Z$.
We say that a Hermitian metric $h$ of $E$ is $C^{\infty}$
if for each $k\in\seisuu_{>0}$
there exists a $C^k$-metric $h'_k$ of $E'$
such that $h'_{k|Z}=h'$.
A section $s$ of $E$ is called $C^{\infty}$
if for each $k\in\seisuu_{>0}$
there exists a $C^k$-section $s'_k$ of $E'$
such that $s'_{k|Z}=s$.
In the following, Hermitian metrics 
and sections are $C^{\infty}$ unless otherwise specified.

Let us recall an important theorem of Donaldson.
\begin{prop}[Donaldson \cite{Donaldson-boundary-value}]
Let $h_{E,\del Z}$ be 
a $G$-invariant $C^{\infty}$-Hermitian metric
of $E_{|\del Z}$.
Then, there exists a unique $G$-invariant
$C^{\infty}$-Hermitian metric $h_E$
of $E$ satisfying the following condition.
\begin{itemize}
\item
 $\Lambda F(h_E)=0$.
\item
 $h_{E|\del Z}=h_{E,\del Z}$.
\hfill\qed
\end{itemize}
\end{prop}
Note that the $G$-invariance of $h_E$ follows from
the uniqueness.

Let $h_0$ be any 
$G$-equivariant $C^{\infty}$-Hermitian metric of $(E,\delbar_E)$.
\begin{cor}
\label{cor;17.11.28.10}
There exits the $G$-invariant Hermitian-Einstein metric $h_E$ of $E$
such that
(i) $\det(h_E)=\det(h_0)$ on $Z$,
(ii) $h_{E|\del Z}=h_{0|\del Z}$.
\end{cor}
\pf
There exists the $G$-invariant Hermitian-Einstein metric $h_{1}$ of $E$
such that
(i) $\Lambda F(h_1)=0$,
(ii) $h_{1|\del Z}=h_{0|\del Z}$.
We obtain the $\real$-valued function $a$ determined by
$\det(h_0)=\det(h_1)\cdot e^{-a}$.
Set 
$h_E:=h_1e^{-a/\rank E}$.
By the construction,
we obtain
$\det(h_E)=\det(h_1 e^{-a/\rank E})=\det(h_1)e^{-a}=\det(h_0)$.
Because $h_{1|\del Z}=h_{0|\del Z}$,
we obtain $a_{|\del Z}=0$,
which implies $h_{E|\del Z}=h_{0|\del Z}$.
Moreover, we obtain
$F(h_E)=F(h_1)-(\delbar\del a)(\rank E)^{-1}\id_E$,
and hence
$F(h_E)^{\bot}=F(h_1)^{\bot}$.
It implies 
the Hermitian-Einstein condition
$\Lambda F(h_E)^{\bot}=0$ for $h_E$.
Hence, $h_E$ satisfies the desired conditions.
\hfill\qed

\subsubsection{Non-increasing property of the Donaldson functional
along the heat equation}
\label{subsection;17.8.10.21}

Let us recall the construction of
the Donaldson functional in this context
by following \cite{Simpson88}.

\paragraph{Preliminary}
Let $U$ be any finite dimensional complex vector space
with a Hermitian metric $h_U$.
Let $b$ be an endomorphism of $U$
which is self-adjoint with respect to $h_U$.
Let $e_1,\ldots,e_r$ be an orthonormal base of $U$
such that $b(e_i)=\lambda_ie_i$.
Let $e_1^{\lor},\ldots,e_r^{\lor}$
denote the dual base of $U^{\lor}$.
For any function $f:\real\lrarr\real$,
let $f(b)$ denote the endomorphism of $U$
determined by 
$f(b)(e_i)=f(\lambda_i)e_i$.
It is self-adjoint with respect to $h_U$.
For any function
$\Phi:\real\times\real\lrarr \real$,
let $\Phi(b)$ denote the endomorphism of $\End(U)$
determined as
$\Phi(b)(e_i^{\lor}\otimes e_j):=
 \Phi(\lambda_i,\lambda_j)e_i^{\lor}\otimes e_j$.
The endomorphism $\Phi(b)$ is self-adjoint
with respect to the induced metric on $\End(U)$.
If $f$ (resp. $\Phi$) is positive valued,
then $f(b)$ (resp. $\Phi(b)$) is positive definite.

Let $E$ be a vector bundle on $Z$
with a $C^{\infty}$-Hermitian metric $h_0$.
Let $b$ be a $C^{\infty}$-section of $\End(E)$
which is self-adjoint with respect to $h_0$.
For any $C^{\infty}$-function $f:\real\lrarr\real$,
we obtain a $C^{\infty}$-section $f(b)$ of $\End(E)$
which is self-adjoint with respect to $h_0$.
For any $C^{\infty}$-function 
$\Phi:\real\times\real\lrarr\real$,
we obtain a $C^{\infty}$-endomorphism
$\Phi(b)$ of $\End(E)$,
which is self-adjoint with respect to $h_0$.

\paragraph{Donaldson functional}

Let $(E,\delbar_E)$ be a holomorphic vector bundle
on $Z$ as above.
Let $h_0$ be a $C^{\infty}$-Hermitian metric of $E$.
Let $\nbigp$ denote the space of $C^{\infty}$-Hermitian metrics
of $E$ such that $h_{|\del Z}=h_{0|\del Z}$.

We take $h_1,h_2\in\nbigp$.
Let $s$ be the endomorphism of $E$
which is self-adjoint with respect to $h_i$ $(i=1,2)$
determined by the condition
$h_2=h_1e^{s}$.
By the construction,
$s_{|\del Z}=0$ holds.
We set
\[
 \Psi(\lambda_1,\lambda_2):=
 \frac{e^{\lambda_2-\lambda_1}-(\lambda_2-\lambda_1)-1}{(\lambda_2-\lambda_1)^2}.
\]
Let $\htilde_1$ be the metric on
$\End(E)\otimes\Omega^{p,q}$
induced by $h_1$ and $g_Z$.
(In the case $p=q=0$, we use the notation $h_1$.)
Then, we put
\begin{equation}
\label{eq;18.12.24.53}
 M(h_1,h_2):=
 \sqrt{-1}\int_Z\Tr(s\Lambda F(h_1))
+\int_Z
 \htilde_1\bigl(
 \Psi(s)(\delbar_Es),\delbar_Es
 \bigr).
\end{equation}

\begin{prop}
\label{prop;17.8.10.2}
For any $h_i\in\nbigp$ $(i=1,2,3)$,
$M(h_1,h_3)
=M(h_1,h_2)+M(h_2,h_3)$
holds.
\end{prop}
\pf
We just apply the argument of Donaldson 
\cite{Donaldson-surface,Donaldson-infinite}
and Simpson \cite{Simpson88}
under the Dirichlet condition.
We give only an outline.

We take a sufficiently large $p>0$.
Let $\nbigpbar$ denote the space of 
$L_2^p$-Hermitian metrics $h$ of $E$
such that $h_{|\del Z}=h_{0|\del Z}$.
It is naturally a Banach manifold.
(See \cite{Hamilton-book} for
Sobolev spaces on manifolds with boundary.)
For each $h\in\nbigpbar$,
the tangent space of $\nbigpbar$ at $h$
is naturally identified with the space $S_h$ of
$L_2^p$-sections $a$ of $\End(E)$ 
which are self-adjoint with respect to $h$
such that $a_{|\del Z}=0$.
Let $U_h:=\{a\in S_h\,|\,\sup|a|_h\leq 1/2\}$.
We define the map
$\Upsilon_h:U_h\lrarr \nbigpbar$ by
$\Upsilon_h(b)=h\cdot (\id+b)$,
which induces an isomorphism of $U_h$
and a neighbourhood of $h$ in $\nbigpbar$.

For each $h\in\nbigpbar$,
we define the linear map $\Phi_h:S_h\lrarr \real$
by
\[
 \Phi_h(b)
=\sqrt{-1}\int_Z\Tr\bigl(b\Lambda F(h)\bigr).
\]
Thus, we obtain a differential $1$-form $\Phi$ on $\nbigpbar$.

\begin{lem}
$\Phi$ is a closed $1$-form.
\end{lem}
\pf
It is enough to prove 
$(d\Phi)_{h_1}=0$ at each $h_1\in \nbigpbar$.
We take $u,v\in S_{h_1}$.
They naturally gives a vector field on 
the linear space $S_{h_1}$.
By $\Upsilon_{h_1}$,
they induce vector fields on
$\Upsilon_{h_1}(U_{h_1})$,
which are denoted by
$\underline{u}$ and $\underline{v}$.
We have 
$[\underline{u},\underline{v}]=0$.
For $b\in U_{h_1}$,
we obtain
$\underline{u}_{h_1\cdot (\id+b)}=
 (\id+b)^{-1}u\in S_{h_1\cdot (\id+b)}$.
The following holds:
\[
 \Phi_{h_1(\id+\epsilon v)}(\underline{u}_{h_1(\id+\epsilon v)})
=\sqrt{-1}\int_Z
 \Lambda\Tr\Bigl(
 (\id+\epsilon v)^{-1}u F(h_1(\id+\epsilon v))
 \Bigr).
\]
Hence, we obtain
\[
 \Phi_{h_1(\id+\epsilon v)}(\underline{u}_{h_1(\id+\epsilon v)})
-\Phi_{h_1}(\underline{u}_{h_1})
=
 \epsilon
 \sqrt{-1}\int_Z
 \Lambda\Tr
\Bigl(
(-vu)F(h_1)
+u\delbar\del_{h_1} v
+O(\epsilon)
\Bigr).
\]
We obtain the following:
\begin{equation}
 \label{eq;17.8.11.10}
\Bigl(
 \underline{v}
 \Phi_{h}(\underline{u}_{h})
\Bigr)_{h=h_1}
-\Bigl(
 \underline{u}
 \Phi_h(\underline{v}_h)
\Bigr)_{h=h_1}
=\sqrt{-1}\int_Z\Lambda\Tr
 \Bigl(
(-vu+uv)F(h_1)
+u\delbar\del_{h_1} v
-v\delbar\del_{h_1}u
\Bigr).
\end{equation}

Because $u_{|\del Z}=v_{|\del Z}=0$,
we obtain the following by the Stokes formula:
\begin{equation}
\int_Z\Lambda\Tr\bigl(v\delbar\del_{h_1}u\bigr)
=-\int_Z\Lambda\Tr(\delbar v\del_{h_1}u)
=\int_Z\Lambda\Tr(\del_{h_1}u\delbar v)
=-\int_Z\Lambda\Tr(u\del_{h_1}\delbar v)
\end{equation}
Note the relation
$(\delbar\del_{h_1}+\del_{h_1}\delbar)v=[F(h_1),v]$.
Hence, we obtain
\begin{multline}
 \int_{Z}
 \Lambda\Tr\bigl(
 u\delbar\del_{h_1}v
-v\delbar\del_{h_1}u
 \bigr)
=\int_Z\Lambda\Tr
 \bigl(
 u\delbar\del_{h_1}v+
 u\del_{h_1}\delbar v
 \bigr)
=\int_Z\Lambda\Tr
 \bigl(
 u[F(h_1)v]
 \bigr) \\
=\int_Z\Lambda\Tr\Bigl(
 vuF(h_1)-uv F(h_1)
 \Bigr).
\end{multline}
Hence, the right hand side of (\ref{eq;17.8.11.10})
is $0$.
Thus, we obtain the claim of the lemma.
\hfill\qed

\vspace{.1in}

For $h_1,h_2\in \nbigpbar$,
we take a path $\gamma$ from $h_1$ to $h_2$
in $\nbigpbar$,
and we set
$M^D(h_1,h_2):=\int_{\gamma}\Phi$.
We clearly have
$M^D(h_1,h_3)=M^D(h_1,h_2)+M^D(h_2,h_3)$
for $h_i\in\nbigpbar$ $(i=1,2,3)$.
For $s\in S_{h_1}$,
we obtain the following:
\begin{equation}
\frac{d^2}{dt^2}
 M^D(h_1,h_1e^{ts})_{|t=t_0}
=\frac{d}{dt_1}
 \Bigl(
 \frac{d}{dt}
  M^D(h_1,h_1e^{(t_1+t)s})_{|t=0}
 \Bigr)_{|t_1=t_0}
=\frac{d}{dt_1}
\Bigl(
 \frac{d}{dt}
 M^D(h_1e^{t_1s},h_1e^{(t_1+t)s})_{|t=0}
\Bigr)_{|t_1=t_0}.
\end{equation}
By the construction of $M^D$,
the following holds:
\[
 \frac{d}{dt}
 M^D(h_1e^{t_1s},h_1e^{(t_1+t)s})_{|t=0}
=\sqrt{-1}\int_Z\Tr(s\Lambda F(h_1e^{t_1s})).
\]
Because
$F(h_1e^{(t_0+t)s})
=F(h_1e^{t_0s})+\delbar\bigl(
e^{-ts}\del_{h_1e^{t_0s}}(e^{ts})
\bigr)$,
the following holds:
\[
 \frac{d}{dt_1}
 F(h_1e^{t_1s})_{|t_1=t_0}
=\frac{d}{dt}
 F(h_1e^{(t_0+t)s})_{|t=0}
=\delbar\del_{h_1e^{t_0s}}(s).
\]
Hence, we obtain
\begin{equation}
\label{eq;17.8.11.41}
 \frac{d^2}{dt^2}
 M^D(h_1,h_1e^{ts})_{|t=t_0}
=\sqrt{-1}\int_Z
 \Lambda \Tr(s\delbar\del_{h_1e^{t_0s}}s).
\end{equation}

As given in the proof of \cite[Proposition 5.1]{Simpson88},
we have the following:
\begin{equation}
\label{eq;17.8.11.20}
 \frac{d^2}{dt^2}
 M(h_1,h_1e^{ts})_{|t=t_0}
=\int_Z
 \htilde_1\bigl(
 \Psi_{1,t_0}(s)\delbar s,\delbar s\bigr).
\end{equation}
Here, $\Psi_{1,t_0}(y_1,y_2)=e^{t_0(y_2-y_1)}$.
Hence, the right hand side of (\ref{eq;17.8.11.20}) is
\begin{equation}
\label{eq;17.8.11.40}
 -\sqrt{-1}\int_Z
 \Lambda
 \Tr\Bigl(
 \delbar s\del_{h_1e^{t_0s}}s
 \Bigr).
\end{equation}
Because $s_{|\del Z}=0$,
we obtain the following equality
for any $0\leq t\leq 1$,
from (\ref{eq;17.8.11.41}),
(\ref{eq;17.8.11.40})
and the Stokes formula:
\[
 \frac{d^2}{dt^2}
 M^D(h_1,h_1e^{ts})
=\frac{d^2}{dt^2}
 M(h_1,h_1e^{ts}).
\]
We clearly have
$M(h_1,h_1e^{ts})_{|t=0}
=0
=M^D(h_1,h_1e^{ts})_{|t=0}$
and 
\[
 \frac{d}{dt}
  M(h_1,h_1e^{ts})_{|t=0}
=\sqrt{-1}
 \int_Z\Tr\bigl(s\Lambda F(h_1)\bigr)
=\frac{d}{dt}
 M^D(h_1,h_1e^{ts})_{|t=0}.
\]
Hence, we obtain
$M(h_1,h_1e^s)=M^D(h_1,h_1e^s)$.
In particular,
we obtain
$M(h_1,h_3)=M(h_1,h_2)+M(h_2,h_3)$
for any $h_i\in\nbigp$ $(i=1,2,3)$.
\hfill\qed

\subsubsection{Preliminary for $C^0$-bound}

Let us consider the heat equation
associated to the Hermitian-Einstein condition:
\begin{equation}
\label{eq;17.8.10.11}
 h_t^{-1}\frac{dh_t}{dt}
=-\sqrt{-1}\Lambda F(h_t)^{\bot}.
\end{equation}
According to Simpson \cite[Corollary 6.5]{Simpson88},
there exists a unique solution $h_t$ of the heat equation
satisfying 
$h_{t|\{0\}\times Z}=h_0$
and $h_{t|\{t\}\times\del Z}=h_{0|\del Z}$.
The following proposition is useful for our study.

\begin{prop}
\label{prop;17.8.10.22}
Let $h_E$ be the Hermitian-Einstein metric of $(E,\delbar_E)$
such that
$\det(h_E)=\det(h_0)$
and $h_{E|\del Z}=h_{0|\del Z}$.
Then, 
$M(h_0,h_E)\leq 0$
holds.
\end{prop}
\pf
Note that if we regard $h_t$ as 
a Hermitian metric of the pull back of
$E$ on $\real_{\geq 0}\times Z$
by the projection $\real_{\geq 0}\times Z\lrarr Z$,
then $h_t$ is $C^{\infty}$
on $(\real_{\geq 0}\times Z)\setminus(\{0\}\times\del Z)$.
(See \cite[Part IV, \S8]{Hamilton-book}
for the regularity outside of the corner $\{0\}\times\del Z$,
which is available in our case
as remarked in the proof of \cite[Corollary 6.5]{Simpson88}.)
We also recall that $h_t$ is $L_2^p$ for any large $p$
around any point of $\{0\}\times \del Z$,
(See the construction of 
a short-time $L_2^p$-solution of a non-linear heat equation
in \cite[\S IV-11]{Hamilton-book}.)
It implies that $h_t$ is $C^0$
even around any point of $\{0\}\times\del Z$,
according to the Sobolev embedding
in \cite[\S II-15]{Hamilton-book}.

For $t_1,t_2\in\real_{\geq 0}$,
let $b_{t_1,t_2}$ be the automorphism of $E$
determined by $h_{t_2}=h_{t_1}b_{t_1,t_2}$,
which is self-adjoint and positive definite
with respect to both of
$h_{t_1}$ and $h_{t_2}$.
Let $s_{t_1,t_2}$ be the endomorphism of $E$,
which is self-adjoint with respect to $h_{t_i}$ $(i=1,2)$
and determined by $b_{t_1,t_2}=e^{s_{t_1,t_2}}$.
Let $\Delta_Z:=\sqrt{-1}\Lambda\delbar_Z\del_Z$.
By using
\cite[Lemma 3.1]{Simpson88}
and $\Tr\Lambda F(h_{t_1})=\Tr\Lambda F(h_{t_2})$,
we obtain the following equality:
\[
 \Delta_Z \Tr(b_{t_1,t_2})
=\sqrt{-1}\Tr\Bigl( (b_{t_1,t_2}-\id)
 \bigl(\Lambda F(h_{t_2})-\Lambda F(h_{t_1})
 \bigr)
 \Bigr)
-\bigl|
 \delbar_E(b_{t_1,t_2})b_{t_1,t_2}^{-1/2}
 \bigr|_{\htilde_{t_1}}^{2}.
\]
Because $\Tr(b_{t_1,t_2})\geq \rank E$ on $Z$,
and because $\Tr(b_{t_1,t_2})=\rank E$ on the boundary $\del Z$,
we obtain $\del_{\nu}\Tr(b_{t_1,t_2})\leq 0$ 
at $\del Z$,
where $\del_{\nu}$ denote the outward unit normal vector field
of $\del Z\subset Z$.
By Green's formula (Lemma \ref{lem;17.11.27.3}),
we obtain
\[
 \int_Z\Delta_Z\Tr(b_{t_1,t_2})
=-\frac{1}{2}\int_{\del Z}
 \del_{\nu}\Tr(b_{t_1,t_2})\geq 0.
\]
Hence, we obtain
\begin{equation}
\label{eq;18.12.24.31}
 \int_Z\bigl|
 \delbar_E(b_{t_1,t_2})b_{t_1,t_2}^{-1/2}
 \bigr|^2_{\htilde_{t_0}}
\leq
 \int_Z
 \bigl|
 (b_{t_1,t_2}-\id)
 \bigr|_{h_{t_1}}
\cdot
 \bigl|
 (\Lambda F(h_{t_2})-\Lambda F(h_{t_1}))
 \bigr|_{h_{t_1}}.
\end{equation}

\begin{lem}
\label{lem;18.12.24.60}
$M(h_0,h_t)$ is a continuous function of $t\in\real_{\geq 0}$.
Moreover, on $\real_{>0}$,
$M(h_0,h_t)$ is $C^{\infty}$-class,
and the following equality holds:
\begin{equation}
\label{eq;18.12.24.30}
 \frac{d}{dt}M(h_0,h_t)
=-\int_{Z}\bigl|
 \Lambda F(h_t)^{\bot}
 \bigr|_{h_t}^2.
\end{equation}
\end{lem}
\pf
We closely follow the argument of Simpson
in \cite[Lemma 7.1]{Simpson88}.
Indeed, because $Z$ is compact, the proof is easier.
Because $h_t$ is $C^{\infty}$ on $\real_{>0}\times Z$,
$M(h_0,h_t)$ is $C^{\infty}$ for $t\in\real_{>0}$.
Let us prove (\ref{eq;18.12.24.30}).
Because of Proposition \ref{prop;17.8.10.2},
it is enough to prove the following for any $t_0>0$.
\[
 \frac{d}{dt}M(h_{t_0},h_t)_{|t=t_0}
=-\int_Z\bigl|
 \Lambda F(h_{t_0})^{\bot}
 \bigr|^2_{h_{t_0}}.
\]
We put 
$I_{\epsilon}:=\{t\in\real\,|\,t_0-\epsilon<t<t_0+\epsilon\}$
for a small number $\epsilon>0$.
Because $h_t$ is $C^{\infty}$ 
on $I_{\epsilon}\times Z$,
there exists $C_1>0$ such that
$|s_{t_0,t}|_{h_{t_0}}\leq C_1|t-t_0|$
in $I_{\epsilon}\times Z$.
By the heat equation,
we have
$\lim_{t\to t_0}\frac{1}{t-t_0}s_{t_0,t}
=-\sqrt{-1}\Lambda F(h_{t_0})^{\bot}$,
and hence
\[
 \lim_{t\to t_0}
 \frac{\sqrt{-1}}{t-t_0}\int_{Z}
 \Tr(s_{t_0,t}\Lambda F(h_t))
=\int_Z
 \Tr\bigl(\Lambda F(h_{t_0})^{\bot}\cdot 
 \Lambda F(h_{t_0})
 \bigr)
=-\int_Z
 \bigl|
 \Lambda F(h_{t_0})^{\bot}
 \bigr|_{h_{t_0}}^2.
\]
We may assume 
$|b_{t_0,t}-\id|<1/2$ 
for $t\in I_{\epsilon}$.
There exists $C_2>0$
such that
$|\delbar s_{t_0,t}|_{\htilde_{t_0}}
\leq
 C_2|\delbar b_{t_0,t}|_{\htilde_{t_0}}$
for $t\in I_{\epsilon}$.
Hence,
there exists 
$C_3>0$
such that
$\htilde_0(\Psi(s_{t_0,t})\delbar_E s_{t_0,t},\delbar_Es_{t_0,t})
\leq
 C_3
 \bigl|
 \delbar_E(b_{t_0,t})b_{t_0,t}^{-1/2}
 \bigr|_{\htilde_{t_0}}^2$
for $t\in I_{\epsilon}$.
Hence, 
according to the formula (\ref{eq;18.12.24.53}),
it is enough to prove 
\begin{equation}
\label{eq;18.12.24.52}
\lim_{t\to t_0}
 \frac{1}{t-t_0}
 \int_Z
 \bigl|
 \delbar_E(b_t)b_{t_0,t}^{-1/2}
 \bigr|_{\htilde_0}^2
=0.
\end{equation}
Because $b_{t_0,t}$ is $C^{\infty}$
around $\{t_0\}\times Z$,
$(t-t_0)^{-1}(b_{t_0,t}-\id)$ is 
bounded on $(I_{\epsilon}\setminus\{t_0\})\times Z$.
Therefore, 
we obtain (\ref{eq;18.12.24.52})
by using (\ref{eq;18.12.24.31}).
Thus, we obtain (\ref{eq;18.12.24.30}).

\vspace{.1in}
Let us prove the continuity of $M(h_0,h_{t})$ at $t=0$.
Set $b_t:=b_{0,t}$ and $s_t:=s_{0,t}$.
By the continuity of $h_t$,
there exists $T_0>0$
such that $|b_t-\id|_{h_0}\leq 1/2$
if $0\leq t\leq T_0$.
According to \cite[Lemma 6.2]{Simpson88},
$\bigl|\Lambda F(h_t)\bigr|$ is uniformly bounded.
At any point of $Z\setminus \del Z$,
$\lim_{t\to 0}\Lambda F(h_t)=\Lambda F(h_0)$ holds.
Hence, we obtain the following from
(\ref{eq;18.12.24.31}):
\begin{equation}
\label{eq;18.12.24.50}
\lim_{t\to 0}
 \int_Z
 \bigl|
 \delbar_E(b_t)b_t^{-1/2}
 \bigr|_{\htilde_0}^2
=0.
\end{equation}
There exists $C_5>0$
such that
$|\delbar s_t|_{\htilde_0}
\leq
 C_5|\delbar b_t|_{\htilde_0}$
for $0\leq t\leq T_0$.
Hence,
there exists 
$C_6>0$
such that
$\htilde_0(\Psi(s_t)\delbar_E s_t,\delbar_Es_t)
\leq
 C_6
 \bigl|
 \delbar_E(b_t)b_t^{-1/2}
 \bigr|_{\htilde_0}^2$
for any $0\leq t\leq T_0$.
Note that $s_t$ is continuous
and $\lim_{t\to 0} s_t=0$.
Then, we obtain 
$\lim_{t\to 0}M(h_0,h_t)=0=M(h_0,h_0)$
by the formula (\ref{eq;18.12.24.53}).
\hfill\qed

\vspace{.1in}
Let us complete the proof of Proposition \ref{prop;17.8.10.22}.
By Lemma \ref{lem;18.12.24.60},
$M(h_0,h_t)$ is continuous and non-increasing.
Because of $M(h_0,h_0)=0$,
we obtain $M(h_0,h_t)\leq 0$ for any $t$.
According to Donaldson \cite{Donaldson-boundary-value},
a subsequence
$\{h_{t_i}\}$ converges to $h_E$ in $C^{\infty}$
as $t_i\to\infty$.
Hence, we obtain
$M(h_0,h_E)\leq 0$.
\hfill\qed

\begin{rem}
We shall explain a refinement of
Lemma {\rm\ref{lem;18.12.24.60}}
in {\rm\S\ref{subsection;18.12.24.100}}.
\hfill\qed
\end{rem}

We recall the well known general formula
called Green's formula.
\begin{lem}
\label{lem;17.11.27.3}
Let $\del_{\nu}$ denote the outward unit normal vector field
at $\del Z$.
Let $\Delta_{g_Z}$ denote the Laplacian
associated to the K\"{a}hler metric $g_Z$.
Note that 
$\Delta_{g_Z}=2\Delta_Z=2\sqrt{-1}\Lambda\delbar_Z\del_Z$.
Then, the following equality holds
for any $C^{\infty}$-functions $f_i$ $(i=1,2)$ on $Z$:
\begin{equation}
\label{eq;17.11.27.2}
 \int_Z f_1\Delta_{g_Z}f_2\dvol_Z
=\int_Zf_2\Delta_{g_Z} f_1\dvol_Z
-\int_{\del Z}
 f_1\del_{\nu}f_2\dvol_{\del Z}
+\int_{\del Z}
 f_2\del_{\nu}f_1\dvol_{\del Z}.
\end{equation}
(For example,
see {\rm\cite[Problem 14-8]{Lee-manifolds}}.)
\hfill\qed
\end{lem}

\subsubsection{Appendix: Complement for the Donaldson functional 
along the heat flow}
\label{subsection;18.12.24.100}

Let us prove that the Donaldson functional
$M(h_0,h_t)$ is right differentiable even at $t=0$
(Proposition \ref{prop;17.8.10.12}).
This was originally given for the proof of Proposition 
\ref{prop;17.8.10.22}
although we replace it with a simpler statement
(Lemma \ref{lem;18.12.24.60}).
The argument of the proof is 
contained in \cite{Donaldson-surface, Simpson88}.
We include it for completeness.
We use the notation in the proof of 
Proposition \ref{prop;17.8.10.22}.

Let us recall some basic inequalities
for positive definite self-adjoint automorphisms
which was introduced by 
Donaldson \cite{Donaldson-surface}
to the study of the heat flow 
associated to Hermitian-Einstein condition.

\begin{lem}
\label{lem;18.12.24.1}
Let $U$ be a finite dimensional $\cnum$-vector space
with a Hermitian metric $h_U$.
Let $\nbigh_+(U,h_U)$ denote the set of
automorphisms $f$ of $U$
such that 
(i) $f$ are self-adjoint and positive definite
 with respect to $h_U$,
(ii) $\det(f)=1$.
For any $0<R$,
we set
$\nbigh_+(U,h_U,R):=
\bigl\{
 f\in\nbigh_+(U,h_U,R)\,\big|\,
 |f-\id|_{h_U}\leq R
\bigr\}$.
Then, the following holds.
\begin{itemize}
\item
For any $f\in \nbigh_+(U,h_U)$,
the inequality
$\Tr(f)\geq\dim_{\cnum}U$ holds.
The equality
$\Tr(f)=\dim_{\cnum}U$ holds if and only if 
$f=\id_U$.
\item
For any $0<R$, there exists a positive number $C(R,\dim_{\cnum}U)$
depending only on $R$ and $\dim_{\cnum}U$
such that
the following holds for any
$f\in\nbigh_+(U,h_U,R)$:
\[
 C(R,\dim_{\cnum}U)^{-1}|f-\id_U|^2_{h_U}\leq
 \Tr(f)+\Tr(f^{-1})-2\dim_{\cnum}U
\leq 
 C(R,\dim_{\cnum}U)|f-\id_U|^2_{h_U}.
\]
\end{itemize}
\end{lem}
\pf
The first claim follows from
the basic inequality
$r^{-1}\sum_{i=1}^r a_i\geq
 \prod a_i^{1/r}$ for any $(a_i)\in\real_{>0}^r$.
The second claim follows from
$a+a^{-1}-2=a^{-1}(a-1)^{2}$ for any $a\in\real_{>0}$.
\hfill\qed

\vspace{.1in}

Let us study the behaviour of $|b_t-\id|$ as $t\to 0$
by following \cite{Donaldson-surface} and \cite{Simpson88}.
We set 
\[
 \nbigf(b_t):=\Tr(b_t)+\Tr(b_t^{-1})-2\rank(E).
\]
Because $b_t$ is continuous,
Lemma \ref{lem;18.12.24.1} implies that
there exist a constant $C_{10}>1$
such that the following holds for $0\leq t\leq T_0$:
\[
 C_{10}^{-1}|b_t-\id|_{h_0}^2
\leq
 \nbigf(b_t)
\leq
 C_{10}|b_t-\id|_{h_0}^2.
\]
According to \cite[Proposition 6.2]{Simpson88},
there exists $C_{11}>0$ such that
$\sup|\Lambda F(h_t)|_{h_0}\leq C_{11}$
for any $0\leq t\leq T_0$.

\begin{lem}
\label{lem;18.12.24.10}
There exists $C_{12}>0$ such that
\[
 (\del_t+\Delta_Z)\nbigf(b_t)
\leq
 C_{12}\nbigf(b_t)^{1/2}.
\]
Indeed, we may put $C_{12}=2C_{10}C_{11}$.
\end{lem}
\pf
We obtain the following equality for $b_t$
from the heat equation for $h_t$
(see \cite[\S6]{Simpson88}):
\[
 \del_tb_t
+\sqrt{-1}\,\delbar\del_{h_0} b_t
=
-\sqrt{-1}b_t\Lambda F(h_0)^{\bot}
+\sqrt{-1}\Lambda\delbar(b_t)b_t^{-1}\del_{h_0}(b_t).
\]
Note that
$\Tr\bigl(\del_tb_t+\sqrt{-1}\,\delbar\del_{h_0}b_t\bigr)
=(\del_t+\Delta_Z)\Tr(b_t)$.
We have the following inequality:
\[
\Bigl|
 \Tr\bigl(b_t \Lambda F(h_0)^{\bot}\bigr)
\Bigr|
=
\Bigl|
 \Tr\bigl((b_t-\id) \Lambda F(h_0)^{\bot}\bigr)
\Bigr|
\leq
 \bigl|(b_t-\id)\bigr|_{h_0}
\cdot
 \bigl|\Lambda F(h_0)^{\bot}\bigr|_{h_0}
\leq
 C_{10}C_{11}\nbigf(b_t)^{1/2}.
\]
According to \cite[Lemma 3.1]{Simpson88},
the following holds:
\[
 \sqrt{-1}\Lambda
 \Tr\bigl(
 \delbar(b_t)\cdot b_t^{-1}\cdot \del(b_t)
 \bigr)
=-\bigl|
 \delbar(b_t)b_t^{-1/2}
 \bigr|^2_{\htilde_0}\leq 0.
\]
Hence, we obtain
$(\del_t+\Delta_Z)\Tr(b_t)
\leq
 C_{10}C_{11}\nbigf(b_t)^{1/2}$.

Similarly,
we obtain the following equality for $b_t^{-1}$ from the heat equation
for $h_t$:
\[
 \del_t(b_t^{-1})
+\sqrt{-1}\Lambda\delbar\del_{h_0}(b_t^{-1})
=-\sqrt{-1}\Lambda F(h_0)^{\bot}\cdot b_t^{-1}
-\sqrt{-1}\Lambda \del_{h_0}(b_t^{-1})\cdot b_t\cdot\delbar(b_t^{-1}).
\]
Note that 
$-\sqrt{-1}\Lambda \Tr\bigl(
 \del_{h_0}(b_t^{-1})\cdot b_t\cdot\delbar(b_t^{-1})
 \bigr)
=-\bigl|
 b_t^{1/2}\delbar(b_t^{-1})
 \bigr|_{\htilde_0}^2\leq 0$.
Hence, 
we obtain
$(\del_t+\Delta_Z)\Tr(b_t^{-1})
\leq
  C_{10}C_{11}\nbigf(b_t)^{1/2}$.
Thus, we obtain the claim of 
Lemma \ref{lem;18.12.24.10}.
\hfill\qed

\begin{prop}
\label{prop;18.12.24.20}
There exists $C_{13}>1$ such that 
the following holds for any $0\leq t\leq T_0$:
\[
 0\leq \nbigf(b_t) \leq C_{13}t^2.
\]
As a result, we obtain
$|b_t-\id|_{h_0}\leq
 C_{10}C_{13}^{1/2}t$
for any $0\leq t\leq T_0$.
Indeed, we may put $C_{13}=C_{10}^2C_{11}^2$.
\end{prop}
\pf
Let $B_0$ be a constant
such that
$\nbigf(b_t)\leq B_0$
for any $0\leq t\leq T_0$.
For any $n\geq 1$,
we put
\[
 B_n:=
 C_{12}^{2-2^{-n+1}}
 B_0^{2^{-(n-1)}}
 \prod_{j=0}^{n-1}
 \left(
 \frac{1}{2-2^{-j}}
 \right)^{2^{-(n-1-j)}}
= C_{12}^{2-2^{-n+1}}
 B_0^{2^{-(n-1)}}
 2^{-(2-2^{-n+1})}
 \prod_{j=0}^{n-1}
 \left(
 \frac{1}{1-2^{-j-1}}
 \right)^{2^{-(n-1-j)}}.
\]
Let us prove that
$\nbigf(b_t)\leq B_nt^{2-2^{-(n-1)}}$
by an induction on $n$.
It holds in the case $n=0$
by our choice of $B_0$.
Suppose 
$\nbigf(b_t)\leq B_nt^{2-2^{-(n-1)}}$ holds.
We obtain
$(\del_t+\Delta_Z)\nbigf(b_t)
\leq
 C_{12}B_n^{1/2}t^{1-2^{-n}}$.
Hence, we obtain
\[
(\del_t+\Delta_Z)\Bigl(
 \nbigf(b_t)-(2-2^{-n})^{-1}C_{12}B_n^{1/2}t^{2-2^{-n}}
 \Bigr)
\leq 0.
\]
Because $\nbigf(b_t)=0$
on $(\{0\}\times Z)\cup\bigl(
 \{0\leq t\leq T_0\}\times\del Z\bigr)$,
we have
$\nbigf(b_t)-(2-2^{-n})C_{12}B_n^{1/2}t^{2-2^{-n}}\leq 0$
on $(\{0\}\times Z)\cup\{0\leq t\leq T_0\}\times\del Z$.
Hence we obtain
$\nbigf(b_t)-(2-2^{-n})C_{12}B_n^{1/2}t^{2-2^{-n}}\leq 0$
on $\{0\leq t\leq T_0\}\times Z$
by the maximum principle
\cite[IV-2]{Hamilton-book}.
It is easy to check
$B_{n+1}=(2-2^{-n})C_{12}B_n^{1/2}$,
and the induction can proceed.

We have the convergence $\lim_{n\to\infty}C_{12}^{2-2^{n+1}}=C_{12}^2$,
$\lim_{n\to\infty}B_0^{2^{-(n-1)}}=1$
and 
$\lim_{n\to\infty}2^{-(2-2^{-n+1})}=2^{-2}$.
The following holds:
\begin{equation}
0\leq
 \log
  \prod_{j=0}^{n-1}
 \left(
 \frac{1}{1-2^{-j-1}}
 \right)^{2^{-(n-1-j)}}
=-\sum_{j=0}^{n-1} 2^{-(n-1-j)}\log(1-2^{-j-1})
\leq
 n2^{-n}.
\end{equation}
Hence, we obtain the convergence:
\[
\lim_{n\to\infty}
  \log
  \prod_{j=0}^{n-1}
 \left(
 \frac{1}{1-2^{-j-1}}
 \right)^{2^{-(n-1-j)}}
=0.
\]
In all,
we obtain
$\nbigf(b_t)\leq 2^{-2}C_{12}^2t^2$.
\hfill\qed

\begin{rem}
Note that 
the estimate in Proposition 
{\rm\ref{prop;18.12.24.20}}
follows from
only 
the boundary condition of $h_t$
and the estimates for the boundedness of
$|b_t-\id|_{h_0}$
and 
$|\Lambda F(h_t)|_{h_0}$.
\hfill\qed
\end{rem}

\begin{prop}
\label{prop;17.8.10.12}
$M(h_0,h_t)$ is right differentiable at $t=0$,
and the following equality holds:
\begin{equation}
\label{eq;18.12.24.61}
 \lim_{t\to 0}\frac{1}{t}M(h_0,h_t)
=-\int_{Z}\bigl|
 \Lambda F(h_t)^{\bot}
 \bigr|_{h_t}^2.
\end{equation}
\end{prop}
\pf
We closely follow the argument of Simpson
in \cite[Lemma 7.1]{Simpson88}.
There exists $C_{14}>1$
such that the following holds for any $0\leq t\leq T_0$:
\[
 C_{14}^{-1}|s_t|^2_{h_0}
\leq
 \nbigf(b_t)
\leq
 C_{14}|s_t|^2_{h_0}.
\]
Hence, $|s_t|_{h_0}\leq C_{15}t$
for some $C_{15}>0$.
By the heat equation,
we have
$\lim_{t\to 0}\frac{1}{t}s_t=-\sqrt{-1}\Lambda F(h_0)^{\bot}$.
Because $h_t$ is $C^{\infty}$ on 
$(\real_{\geq 0}\times Z)\setminus (\{0\}\times\del Z)$,
we obtain $\lim_{t\to 0}\Lambda F(h_t)=\Lambda F(h_0)$
at any point of $Z\setminus\del Z$.
Because $\bigl|\Lambda F(h_t)\bigr|_{h_t}$ is uniformly bounded
on $\{0\leq t\leq T_0\}\times(Z\setminus \del Z)$,
we obtain
\[
 \lim_{t\to 0}
 \frac{\sqrt{-1}}{t}\int_{Z}
 \Tr(s_t\Lambda F(h_t))
=\int_Z
 \Tr\bigl(\Lambda F(h_0)^{\bot}\cdot 
 \Lambda F(h_0)
 \bigr)
=-\int_Z
 \bigl|
 \Lambda F(h_0)^{\bot}
 \bigr|_{h_0}^2.
\]
There exists 
$C_{16}>0$
such that
$\htilde_0(\Psi(s_t)\delbar_E s_t,\delbar_Es_t)
\leq
 C_{16}
 \bigl|
 \delbar_E(b_t)b_t^{-1/2}
 \bigr|_{\htilde_0}^2$
for any $0\leq t\leq T_0$.
By Proposition \ref{prop;18.12.24.20},
$t^{-1}(b_t-\id)$
is uniformly bounded
on $\{0\leq t\leq T_0\}\times (Z\setminus \del Z)$.
By using (\ref{eq;18.12.24.31}),
we obtain
\begin{equation}
\label{eq;17.8.10.10}
\lim_{t\to 0}
 \frac{1}{t}
 \int_Z
 \bigl|
 \delbar_E(b_t)b_t^{-1/2}
 \bigr|_{\htilde_0}^2
=0.
\end{equation}
Then, by the formula (\ref{eq;18.12.24.53}),
we obtain the right differentiability of
$M(h_0,h_t)$ at $t=0$,
and the formula (\ref{eq;18.12.24.61}).
\hfill\qed

\subsection{Preliminary for $C^1$-bound}

\subsubsection{Kodaira identity}

Let $(Y,g_Y)$ be any K\"ahler manifold.
Let $(V,\delbar_V)$ be any holomorphic vector bundle on $Y$
with a Hermitian metric $h$.
Let $A^{p,q}(V)$ denote the space of $C^{\infty}$-sections of
$V\otimes\Omega^{p,q}$.
Let $\del_{V,h}+\delbar_V$ denote the Chern connection of $V$.
The differential operators
$A^{r,s}(V)\lrarr A^{r+1,s}(V)$ induced by $\del_{V,h}$
are also denoted by $\del_{V,h}$.

Take a non-negative integer $p$.
We set $\Vtilde:=V\otimes\Omega_Y^{p,0}$.
It is naturally a holomorphic vector bundle,
and equipped with the induced Hermitian metric $\htilde$.
(In the case $p=0$, we also use the notation $h$.)
Let $\del_{\Vtilde,\htilde}+\delbar_{\Vtilde}$
denote the Chern connection of 
$(\Vtilde,\delbar_{\Vtilde},\htilde)$.
The differential operators
$A^{r,s}(\Vtilde)\lrarr A^{r+1,s}(\Vtilde)$
are also denoted by 
$\del_{\Vtilde,\htilde}$.

The wedge product 
induces an isomorphism
$\Omega^{p,0}_Y\otimes\Omega^{0,q}_Y
\simeq
 \Omega^{p,q}_Y$
given by
$\rho^{p,0}\otimes\rho^{0,q}\longmapsto
 \rho^{p,0}\wedge\rho^{0,q}$.
It induces an isomorphism
$A^{0,q}(\Vtilde)\simeq A^{p,q}(V)$.
We identify
$A^{0,q}(\Vtilde)$ 
and $A^{p,q}(V)$
by the isomorphism.
Note that
we have
$\delbar_{\Vtilde}(\rho)=(-1)^p\delbar_V(\rho)$
for any $\rho\in A^{0,q}(\Vtilde)$
under this identification.

For any $\tau\in A^{p,1}(V)=A^{0,1}(\Vtilde)$,
we obtain the  following elements:
\begin{itemize}
\item
The operator
$\Lambda:A^{1,1}(\Vtilde)\lrarr A^{0,0}(\Vtilde)$
and the element $\del_{\Vtilde,\htilde}(\tau)
 \in A^{1,1}(\Vtilde)$
induce
$\Lambda\del_{\Vtilde,\htilde}(\tau)\in A^{0,0}(\Vtilde)$.
\item
The operator
$\Lambda:A^{p+1,1}(V)\lrarr A^{p,0}(V)$
and the element
$\del_{V,h}(\tau)$
induce
$\Lambda\del_{V,h}(\tau)
\in A^{p,0}(V)$.
\item
The operator
$\del_{V,h}:A^{p-1,0}(V)\lrarr A^{p,0}(V)$
and the element
$\Lambda\tau\in A^{p-1,0}(V)$
induce
$\del_{V,h}\Lambda\tau\in A^{p,0}(V)$.
\end{itemize}

\begin{lem}
\label{lem;17.8.9.1}
For any $\tau\in A^{p,1}(V)=A^{0,1}(\Vtilde)$,
we have
the following equality in 
$A^{p,0}(V)=A^{0,0}(\Vtilde)$:
\[
 -\sqrt{-1}
 \Lambda\del_{\Vtilde,\htilde}(\tau)
+(-1)^{p+1}\sqrt{-1}\del_{V,h}\Lambda\tau
+(-1)^p\sqrt{-1}\Lambda\del_{V,h}(\tau)=0.
\]
\end{lem}
\pf
It is enough to prove the equality
at each point $P\in Y$.
We take a holomorphic coordinate system
$(z_1,\ldots,z_n)$ around $P$
such that
(i) $z_i(P)=0$,
(ii) $g=\sum dz_i\,d\zbar_i+O\bigl(\sum|z_i|^2\bigr)$.
For any tuple
$I=(i_1,\ldots,i_r)\in\{1,\ldots,n\}^r$,
we set
$dz_I=dz_{i_1}\wedge\cdots \wedge dz_{i_r}$.
For any positive integer $\ell$,
let $\gbigs_{\ell}$ denote the $\ell$-th symmetric group.
For any $\sigma\in\gbigs_{\ell}$,
let $\sgn(\sigma)$ denote the signature of $\sigma$.
We refer \cite[\S3.2]{Kobayashi-vector-bundle}
as a general reference for the following computations.

We have the following expression
of $\tau\in A^{p,1}(V)$:
\[
 \tau=\frac{1}{p!}
 \sum_{|I|=p}\sum_j \tau_{I,\juebar}\,dz_I\wedge d\zbar_j.
\]
We set $\eta:=\del_{V,h}\tau\in A^{p+1,1}(V)$.
We have the expression
\[
 \eta=\frac{1}{(p+1)!}
 \sum_{|J|=p+1}\sum_j \eta_{J,\juebar}\,dz_J\wedge d\zbar_j.
\]
We assume that
$\tau_{I,\juebar}=\sgn(\sigma)\tau_{\sigma(I),\juebar}$
(resp.
$\eta_{I,\juebar}=\sgn(\sigma)\eta_{\sigma(I),\juebar}$)
for any $\sigma\in \gbigs_p$
(resp. $\sigma\in\gbigs_{p+1}$).
We have the relation
\[
 \frac{1}{p!}
 \sum_k\sum_{|I|=p}\sum_j
 \nabla_k\tau_{I,\juebar}\,dz_k\wedge dz_I\wedge d\zbar_j
=\frac{1}{(p+1)!}
 \sum_{|J|=p+1}\sum_j
 \eta_{J,\juebar}\,dz_J\wedge d\zbar_j.
\]
We obtain the following relation
for any $(a_0,\ldots,a_p)\in\{1,\ldots,n\}^{p+1}$:
\[
 \sum_{i=0}^{p}\nabla_{a_i}
 \tau_{(a_0,\ldots,a_{i-1},a_{i+1},\ldots,a_p),\juebar}
 dz_{a_i}\wedge dz_{a_0}\wedge\cdots\wedge
 dz_{a_{i-1}}\wedge dz_{a_{i+1}}\wedge \cdots
 \wedge dz_{a_p}\wedge d\zbar_j
\\
=\eta_{(a_0,\ldots,a_{p}),\juebar}\,
 dz_{a_0}\wedge \cdots \wedge dz_{a_p}\wedge d\zbar_j.
\]
Hence, we obtain the following equality
for any $(a_1,\ldots,a_p)\in \{1,\ldots,n\}^p$
with $a_i\neq j$:
\begin{equation}
 \label{eq;17.11.27.30}
 \eta_{(j,a_1,\ldots,a_p),\juebar}
=
 \nabla_j\tau_{(a_1,\ldots,a_p),\juebar}
+\sum_{i=1}^p(-1)^i
 \nabla_{a_i}
 \tau_{(j,a_1,\ldots,a_{i-1},a_{i+1},\ldots,a_{p}),\juebar}.
\end{equation}

We have the following expression
of $\del_{\Vtilde,\htilde}\tau\in A^{1,1}(\Vtilde)$:
\[
 \del_{\Vtilde,\htilde}\tau
=\frac{1}{p!}
 \sum_{k,j,I}
 \nabla_k
 \bigl(
 \tau_{I,\juebar}
 \otimes dz_I
 \bigr)
 dz_k\wedge d\zbar_j.
\]
Then, $-\sqrt{-1}\Lambda\del_{\Vtilde,\htilde}\tau
\in A^{0,0}(\Vtilde)=A^{p,0}(V)$
is expressed as follows
(see \cite[\S3.2 (3.2.17)]{Kobayashi-vector-bundle}
for the formula of $\Lambda$):
\begin{multline}
\label{eq;17.11.27.21}
 -\sqrt{-1}\Lambda\del_{\Vtilde,\htilde}\tau
=-\frac{1}{p!}
 \sum_{j,I}\nabla_j\tau_{I,\juebar}\,dz_I
+A_0
 \\
=-\frac{1}{p!}
 \sum_j\sum_{\substack{|I|=p\\j\not\in I}}
 \nabla_j\tau_{I,\juebar}dz_I
-\frac{1}{(p-1)!}\sum_j\sum_{\substack{|K|=p-1\\j\not\in K}}
 \nabla_j\tau_{(j,K),\juebar}\,dz_j\wedge dz_K
+A_0.
\end{multline}
Here, $A_0(0,\ldots,0)=0$.
We have the following expression of
$\Lambda\eta$:
\begin{equation}
 \label{eq;17.11.27.20}
 \Lambda\eta
=\frac{1}{p!} (-1)^{p+1}\sqrt{-1}
 \sum_j\sum_{\substack{|I|=p\\j\not\in I}}
 \eta_{(j,I),\juebar}dz_I
+A_1.
\end{equation}
Here, $A_1=O\bigl(\sum|z_i|^2\bigr)$.
We also have
\[
 \Lambda\tau
=\frac{1}{(p-1)!}
 (-1)^p\sqrt{-1}
 \sum_{j}\sum_{\substack{|K|=p-1\\j\not\in K}}
 \tau_{(j,K),\juebar}dz_K
+A_2.
\]
Here, $A_2=O\bigl(\sum|z_i|^2\bigr)$.
Hence, we obtain
\begin{multline}
\label{eq;17.11.27.22}
\del_{V,h}\Lambda\tau
=\frac{1}{(p-1)!}(-1)^p\sqrt{-1}
 \sum_{k}\sum_j\sum_{\substack{|K|=p-1\\ j,k\not\in K}}
 \nabla_{k}
 \tau_{(j,K),\juebar}dz_k\wedge dz_K
+A_3
 \\
=\frac{1}{(p-1)!}(-1)^p\sqrt{-1}
 \sum_j\sum_{\substack{|K|=p-1\\ j\not\in K}}
 \nabla_{j}
 \tau_{(j,K),\juebar}dz_j\wedge dz_K
 \\
+\frac{1}{(p-1)!}(-1)^p\sqrt{-1}
 \sum_{k}\sum_{j\neq k}
 \sum_{\substack{|K|=p-1\\ j,k\not\in K}}
 \nabla_{k}
 \tau_{(j,K),\juebar}dz_k\wedge dz_K
+A_3.
\end{multline}
Here, $A_3(0,\ldots,0)=0$.
By  (\ref{eq;17.11.27.30}),
(\ref{eq;17.11.27.21}) and (\ref{eq;17.11.27.22}),
we obtain the following:
\begin{multline}
 -\sqrt{-1}\Lambda\del_{\Vtilde,\htilde}\tau
+(-1)^{p+1}\sqrt{-1}\del_{V,h}\Lambda\tau
=\\
-\frac{1}{p!}
 \sum_j\sum_{\substack{|I|=p\\ j\not\in I}}
 \nabla_j\tau_{I,\juebar}dz_I
+\frac{1}{(p-1)!}
 \sum_{k}\sum_{j\neq k}
 \sum_{\substack{|K|=p-1\\ j,k\not\in K}}
 \nabla_{k}
 \tau_{(j,K),\juebar}dz_k\wedge dz_K
+A_4
\\
=-\frac{1}{p!}
 \sum_j\sum_{\substack{|I|=p\\ j\not\in I}}
 \eta_{(j,I),\juebar}dz_I+A_4.
\end{multline}
Here, $A_4(0,\ldots,0)=0$.
By (\ref{eq;17.11.27.20}),
we obtain that 
$-\sqrt{-1}\Lambda\del_{\Vtilde,\htilde}\tau
+(-1)^{p+1}\sqrt{-1}\del_{V,h}\Lambda\tau
+(-1)^{p}\sqrt{-1}\Lambda\eta$
is $0$ at $P$.
Thus, we obtain the claim of Lemma \ref{lem;17.8.9.1}.
\hfill\qed

\subsubsection{Hermitian-Einstein metrics and other metrics}

Suppose that $h$ is a Hermitian-Einstein metric of $(V,\delbar_V)$.
Let $h_0$ be any Hermitian metric of $V$
such that $\det(h_0)=\det(h)$.
Let $b_1$ be the automorphism of $V$
determined by $h_0=h\,b_1$.
Set $\Vtilde:=V\otimes\Omega^{1,0}_Y$.
Let $\htilde$ denote the Hermitian metric of $\Vtilde$
induced by $h$ and $g_Y$.

\begin{lem}
\label{lem;17.10.8.1}
We have the following equality:
\[
-\sqrt{-1}\Lambda\del_{\Vtilde,\htilde}\delbar_{\Vtilde}
 (b_1^{-1}\del_{V,h}b_1)
=\sqrt{-1}
 \del_{V,h_0}\Lambda F(h_0)^{\bot}
-\sqrt{-1}\bigl[
 b_1^{-1}\del_{V,h}b_1,\Lambda F(h_0)^{\bot}
 \bigr]
+\sqrt{-1}\Lambda\bigl[b_1^{-1}\del_{V,h}b_1,F(h_0)^{\bot}\bigr].
\]
\end{lem}
\pf
By applying Lemma \ref{lem;17.8.9.1} to 
$\delbar_{\Vtilde}(b_1^{-1}\del_{V,h}b_1)
=-\delbar_{V}(b_1^{-1}\del_{V,h}b_1)$,
we obtain the following equality:
\[
 -\sqrt{-1}
\Lambda\del_{\Vtilde,\htilde}
 \delbar_{\Vtilde}(b_1^{-1}\del_{V,h}b_1)
=\sqrt{-1}\del_{V,h}\Lambda\bigl(
 \delbar_V(b_1^{-1}\del_{V,h}b_1)\bigr)
-\sqrt{-1}\Lambda
 \del_{V,h}\delbar_V(b_1^{-1}\del_{V,h}b_1).
\]
Because $h_0=hb_1$,
we have 
$\del_{V,h_0}=\del_{V,h}+b_1^{-1}\del_{V,h}b_1$
and $F(h_0)-F(h)=\delbar_V (b_1^{-1}\del_{V,h}b_1)$.
By the Hermitian-Einstein condition $\Lambda F(h)^{\bot}=0$,
and by $\det(h)=\det(h_0)$,
we obtain
$\Lambda\bigl(F(h_0)-F(h)\bigr)
=\Lambda F(h_0)^{\bot}$,
and hence
\begin{multline}
\sqrt{-1}\del_{V,h}\Lambda\bigl(
 \delbar_V(b_1^{-1}\del_{V,h}b_1)\bigr)
=\sqrt{-1}(\del_{V,h_0}-b_1^{-1}\del_{V,h}b_1)
 \Lambda F(h_0)^{\bot}
 \\
=\sqrt{-1}\del_{V,h_0}\Lambda F(h_0)^{\bot}
-\sqrt{-1}\bigl[b_1^{-1}\del_{V,h}b_1,\Lambda F(h_0)^{\bot}\bigr].
\end{multline}
By the Bianchi identity,
we have 
$\del_{V,h_0}F(h_0)=0$
and $\del_{V,h}F(h)=0$.
Hence, we obtain the following:
\begin{multline}
-\sqrt{-1}\Lambda \del_{V,h}\delbar_V(b_1^{-1}\del_{V,h}b_1)
=-\sqrt{-1}\Lambda\del_{V,h}\bigl(F(h_0)-F(h)\bigr)
=-\sqrt{-1}\Lambda(\del_{V,h_0}-b_1^{-1}\del_{V,h}b_1)F(h_0)
 \\
=\sqrt{-1}\Lambda\bigl[b_1^{-1}\del_{V,h}b_1,F(h_0)\bigr]
=\sqrt{-1}\Lambda\bigl[b_1^{-1}\del_{V,h}b_1,F(h_0)^{\bot}\bigr].
\end{multline}
Then, we obtain the claim of the lemma.
\hfill\qed

\vspace{.1in}
We regard $g_Y$ as a Hermitian metric of 
the holomorphic vector bundle $\Omega_Y^{1,0}$.
Let $F(g_Y)$ denote the curvature of 
the Chern connection of $(\Omega_Y^{1,0},g_Y)$.
In particular,
we obtain $\Lambda F(g_Y)$ as a section of
$\End(\Omega_Y^{1,0})$.
The induced section of
$\End(\Vtilde)$
is denoted as $\id_V\otimes\Lambda F(g_Y)$.
Set $\Delta_Y:=-\sqrt{-1}\del_Y\delbar_Y=\sqrt{-1}\,\delbar_Y\del_Y$.

\begin{lem}
\label{lem;17.10.8.2}
We have the following inequality:
\begin{multline}
 \Delta_Y\bigl|b_1^{-1}\del_{V,h}b_1\bigr|^2_{\htilde}
\leq
 2\Re \htilde\bigl(
-\sqrt{-1}\del_{\Vtilde,\htilde}\delbar_{\Vtilde}
(b_1^{-1}\del_{V,h}b_1),\,
 b_1^{-1}\del_{V,h}b_1
 \bigr)
 \\
+\htilde\bigl(
 b_1^{-1}\del_{V,h}b_1,\,
 \sqrt{-1}\bigl[\id_V\otimes \Lambda F(g_Y),b_1^{-1}\del_{V,h}b_1\bigr]
 \bigr).
\end{multline}
\end{lem}
\pf
Let $F(\htilde)$ denote the curvature of the Chern connection
of $(\Vtilde,\delbar_{\Vtilde},\htilde)$.
Recall
\[
 F(\htilde)=F(h)\otimes\id_{\Omega_Y^{1,0}}
+\id_V\otimes F(g_Y).
\]
Because $\Lambda F(h)^{\bot}=0$,
the adjoint action of $\Lambda F(\htilde)$
on $\End(\Vtilde)$
is equal to the adjoint action of
$\id_V\otimes\Lambda F(g_Y)$.
Hence, 
for any section $s$ of $\End(\Vtilde)$,
the following holds:
\[
 \sqrt{-1}\Lambda\bigl(
 \del_{\Vtilde,\htilde}\delbar_{\Vtilde}
+\delbar_{\Vtilde}\del_{\Vtilde,\htilde}
 \bigr)s
=\sqrt{-1}
 \bigl[
 \id_V\otimes\Lambda F(g_Y),s
 \bigr].
\]

The following holds:
\begin{multline}
-\sqrt{-1}\Lambda
 \del\delbar \htilde(b_1^{-1}\del_hb_1,b_1^{-1}\del_hb_1)
=\\
\htilde\Bigl(
-\sqrt{-1}\Lambda\del_{\Vtilde,\htilde}\delbar_{\Vtilde}
 (b_1^{-1}\del_hb_1),
 b_1^{-1}\del_hb_1
 \Bigr)
+\htilde\Bigl(
  b_1^{-1}\del_hb_1,
 \sqrt{-1}\Lambda\delbar_{\Vtilde}\del_{\Vtilde,\htilde}
 (b_1^{-1}\del_hb_1)
 \Bigr)
\\
+\sqrt{-1}\Lambda
 \htilde\bigl(\delbar_{\Vtilde}(b_1^{-1}\del_hb_1),
 \delbar_{\Vtilde}(b_1^{-1}\del_hb_1)\bigr)
-\sqrt{-1}\Lambda
 \htilde\bigl(
 \del_{\Vtilde,\htilde}(b_1^{-1}\del_hb_1),
 \del_{\Vtilde,\htilde}(b_1^{-1}\del_hb_1)
 \bigr).
\end{multline}
Because
$
 \sqrt{-1}\,\delbar_{\Vtilde}\del_{\Vtilde,\htilde}
 (b_1^{-1}\del_hb_1)
=-\sqrt{-1}\del_{\Vtilde,\htilde}\delbar_{\Vtilde}
 (b_1^{-1}\del_hb_1)
+\bigl[
 \sqrt{-1}\id_V\otimes F(g_Y),
 b_1^{-1}\del_hb_1
 \bigr]$
and 
\[
 \sqrt{-1}\Lambda
 \htilde\bigl(\delbar_{\Vtilde}(b_1^{-1}\del_hb_1),\,
 \delbar_{\Vtilde}(b_1^{-1}\del_hb_1)\bigr)
\leq 0,
\quad
-\sqrt{-1}\Lambda
 \htilde\bigl(
 \del_{\Vtilde,\htilde}(b_1^{-1}\del_hb_1),
 \del_{\Vtilde,\htilde}(b_1^{-1}\del_hb_1)
 \bigr)\leq 0,
\]
we obtain the claim of the lemma.
\hfill\qed

\subsection{Sequence of 
Hermitian-Einstein metrics
and $C^0$-bound}
\label{subsection;17.11.15.1}

We return to the setting in \S\ref{subsection;17.11.27.4}.
If $X$ is compact,
the claim of Theorem \ref{thm;17.8.8.11}
is contained in the result of Simpson
\cite[Theorem 1]{Simpson88}.
(If $X$ is compact and $G$ is trivial,
it is a result of Donaldson \cite{Donaldson-infinite}
and Uhlenbeck-Yau \cite{Uhlenbeck-Yau}.)
Hence, we shall prove Theorem \ref{thm;17.8.8.11}
under the assumption
where $X$ is non-compact.
Recall that $X$ is assumed to be connected.

\subsubsection{Sequence of $G$-invariant regions
with smooth boundary}

\begin{lem}
\label{lem;18.1.18.2}
There exists a $G$-invariant $C^{\infty}$-function $f$ on $X$
which is exhaustive.
\end{lem}
\pf
We take an increasing sequence of compact subsets
$K_i\subset X$ $(i=1,2,\ldots)$
such that $\bigcup K_i=X$.
We set $\Ktilde_1:=\bigcup_{g\in G}\kappa_g(K_1)$.
Because $\Ktilde_1$
is the image of the map $G\times K_1\lrarr X$,
$\Ktilde_1$ is compact and contains $K_1$.
We take a relatively compact open neighbourhood
$U_1$ of $\Ktilde_1\cup K_2$.
Let $\overline{U_1}$ be the closure of $U_1$ in $X$,
which is compact.
We set $\Ktilde_2:=\bigcup_{g\in G}\kappa_g(\overline{U_1})$.
Inductively,
we construct a sequence of compact subsets
$\Ktilde_i\subset X$ $(i=1,2,\ldots)$
as follows.
Suppose that we have already constructed
$\Ktilde_i$.
We take a relatively compact open neighbourhood 
$U_i$ of $\Ktilde_i\cup K_{i+1}$.
Let $\overline{U_i}$ be the closure of $U_i$ in $X$,
which is compact.
We set 
$\Ktilde_{i+1}:=
 \bigcup_{g\in G}\kappa_g(\overline{U_i})$.
Thus, we obtain an increasing sequence of
$G$-invariant compact subsets $\Ktilde_i\subset X$
$(i=1,2\ldots)$
such that $\Ktilde_{i+1}$ contains
an open neighbourhood of $\Ktilde_i$
and that $\bigcup\Ktilde_{i}=X$.
Let $L_i$ denote the closure of 
$X\setminus\Ktilde_i$.

There exists a $C^{\infty}$-function 
$\psi_i:X\lrarr\{0\leq t\leq 1\}$
such that
$\psi_i(P)=0$ $(P\in \Ktilde_{i-1})$
and 
$\psi_i(P)=1$ $(P\in L_{i})$.
We set
$\psitilde_i:=\int_{G}\kappa_g^{\ast}(\psi_i)$,
where the integral is taken with respect to 
the normalized bi-invariant measure of $G$.
Then, $\psitilde_i:X\lrarr \{0\leq t\leq 1\}$ 
are $C^{\infty}$-functions
such that 
$\psitilde_i(P)=0$ $(P\in \Ktilde_{i-1})$
and 
$\psitilde_i(P)=1$ $(P\in L_{i})$.

We set $f:=\sum_{i=1}^{\infty}\psitilde_i$.
On $\Ktilde_{i+1}\cap L_i$,
we have $\psitilde_j=0$ $(j\geq i+1)$
and $\psitilde_j=1$ $(j<i)$.
Hence, we obtain
$i\leq f\leq i+1$ on $\Ktilde_{i+1}\cap L_i$.
Hence, $f$ is a function
as desired in Lemma \ref{lem;18.1.18.2}.
\hfill\qed

\begin{lem}
\label{lem;18.1.18.3}
There exists an increasing sequence of compact subsets  $X_i$ $(i=1,2,\ldots)$
in $X$ with $\bigcup X_i=X$
satisfying the following condition.
\begin{itemize}
\item
Each $X_i$ is a $G$-invariant
submanifold with non-empty smooth boundary $\del X_i$
such that $X_i\setminus \del X_i$ is an open subset of $X$.
Moreover each connected component of $X_i$
has non-empty boundary.
\end{itemize}
\end{lem}
\pf
We take a $G$-invariant $C^{\infty}$-function 
$f:X\lrarr \real_{\geq 0}$
which is exhaustive
as in Lemma \ref{lem;18.1.18.2}.
According to Sard's theorem,
there exists a sequence $a_i\to\infty$
such that each $a_i$ is not a critical value of $f$.
It is enough to put $X_i:=f^{-1}(\{t\in\real\,|\,t\leq a_i\})$.
Note that if the boundary of a connected component of $X_i$
is empty, the component should be equal to $X$ because
$X$ is assumed to be connected.
But, it contradicts with our assumption
that $X$ is non-compact.
Hence, the boundary of any connected component of $X_i$
is not empty.
Moreover, because $\del X_i=f^{-1}(a_i)$,
and because $a_i$ is not a critical value of $f$,
$\del X_i$ is smooth.
\hfill\qed

\subsubsection{Sequence of Hermitian-Einstein metrics}

We take a sequence of compact subsets $X_i\subset X$ 
as in Lemma \ref{lem;18.1.18.3}.
We set
$(E_i,\delbar_{E_i},h_{0,i}):=
 (E,\delbar_E,h_0)_{|X_i}$.
According to the theorem of Donaldson
(see Corollary \ref{cor;17.11.28.10}),
there exists a unique $G$-invariant Hermitian-Einstein metric $h_i$ 
of $(E_i,\delbar_{E_i})$
such that
$h_{i|\del X_i}=h_{0,i|\del X_i}$
and $\det(h_{i})=\det(h_{0,i})$.
Let $s_i$ be the endomorphism of $E_i$
determined by $h_i=h_{0,i}e^{s_i}$
which is self-adjoint with respect to both $h_i$ and $h_{0,i}$.
Note that $\Tr(s_i)=0$.

We shall prove the following boundedness of the sequence,
which is the counterpart of
\cite[Proposition 5.3]{Simpson88}.
\begin{prop}
\label{prop;17.8.8.10}
Suppose that
$(E,\delbar_E,h_0)$ is analytically stable.
Then, there exists a positive constant $C_1$
such that
$\sup_{X_i}|s_i|_{h_{0,i}}\leq C_1$
for any $i$.
\end{prop}

We explain only an outline of the proof.

\subsubsection{Comparison of the sup-norms and the $L^1$-norms}

Set $r:=\rank E$.
According to \cite[Lemma 3.1]{Simpson88},
we have the following inequality
on $X_i$:
\[
 \Delta_{X_i}\bigl(
 \log\bigl(
 \Tr(e^{s_i})/r
 \bigr)
 \bigr)
\leq
 \bigl|
 \Lambda F(h_{0,i})
 \bigr|_{h_{0,i}}.
\]
We extend
$\log\bigl(
 \Tr(e^{s_i})/r
 \bigr)
 \bigr)$
and 
$\bigl|
 \Lambda F(h_{0,i})
 \bigr|_{h_{0,i}}$
to the functions 
$\log\bigl(
 \Tr(e^{s_i})/r
 \bigr)
 \bigr)^{\sim}$
and 
$\bigl|
 \Lambda F(h_{0,i})
 \bigr|_{h_{0,i}}^{\sim}$
on $X$
by setting $0$ outside $X_i$.

\begin{lem}
The following inequality holds as distributions on $X$:
\[
 \Delta_{X}\bigl(
 \log\bigl(
 \Tr(e^{s_i})/r
 \bigr)
 \bigr)^{\sim}
\leq
 \bigl|
 \Lambda F(h_{0,i})
 \bigr|_{h_{0,i}}^{\sim}.
\]
\end{lem}
\pf
Let $\del_{\nu,i}$ be the outward unit normal vector field
at $\del X_i$.
Note that $\log\bigl(\Tr(e^{s_i})/r\bigr)=0$
on $\del X_i$,
and that $\log\bigl(\Tr(e^{s_i})/r\bigr)\geq 0$
on $X_i$.
Hence, we obtain
\[
 \del_{\nu,i}
\bigl(
  \log\bigl(\Tr(e^{s_i})/r\bigr)
\bigr)
\leq 0.
\]

Let $\psi$ be any $\real_{\geq 0}$-valued test function on $X$.
By Green's formula (Lemma \ref{lem;17.11.27.3}),
we obtain the following:
\begin{multline}
 \int_{X_i}
 \log\bigl(\Tr(e^{s_i})/r\bigr)
 \cdot
 \Delta_{X_i}(\psi_{|X_i})
=\int_{X_i}
 \Delta_{X_i}
 \log\bigl(\Tr(e^{s_i})/r\bigr)\cdot\psi
 \\
-\frac{1}{2}\int_{\del X_i}
  \log\bigl(\Tr(e^{s_i})/r\bigr)\cdot
\del_{\nu,i}\psi
+\frac{1}{2}\int_{\del X_i}
 \psi\cdot
 \del_{\nu,i}
  \log\bigl(\Tr(e^{s_i})/r\bigr).
\end{multline}
We have
$\int_{\del X_i}
  \log\bigl(\Tr(e^{s_i})/r\bigr)\cdot
\del_{\nu,i}\psi=0$
and 
$\int_{\del X_i}
 \psi\cdot
 \del_{\nu,i}
  \log\bigl(\Tr(e^{s_i})/r\bigr)\leq 0$.
Hence, we obtain
\[
 \int_{X_i}
 \log\bigl(\Tr(e^{s_i})/r\bigr)
 \cdot
 \Delta_{X_i}(\psi_{|X_i})
\leq
 \int_{X_i}
 \Delta_{X_i}
 \log\bigl(\Tr(e^{s_i})/r\bigr)\cdot\psi
\leq
 \int_{X_i}
 \bigl|\Lambda F(h_0)
 \bigr|_{h_0}\cdot \psi.
\]
It implies the claim of the lemma.
\hfill\qed

\vspace{.1in}
By the assumption on $(X,g_X)$
and on $h_0$,
there exist constants $C_{10},C_{11}>0$
such that
the following holds for any $i$:
\[
 \sup_{X_i}\log
 \bigl(\Tr(e^{s_i})/r\bigr)
\leq
 C_{10}
+C_{11}
 \int_{X_i}
 \log\bigl(\Tr(e^{s_i})/r\bigr)
 \cdot
 \varphi_{X|X_i}.
\]
Hence, there exist positive constants $C_{12},C_{13}$
such that the following holds for any $i$:
\begin{multline}
\label{eq;17.8.10.20}
 \sup_{X_i}|s_i|_{h_0}\leq
 r^{1/2}
 \sup_{X_i} \log\bigl(\Tr(e^{s_i})\bigr)
\leq
 r^{1/2}
 \Bigl(
 \log r
+C_{10}
+C_{11}\int_{X_i}
 \log\bigl(\Tr(e^{s_i})/r\bigr)\cdot\varphi_{X|X_i}
 \Bigr)
 \\
\leq
 C_{12}
+C_{13}\int_{X_i}|s_i|_{h_0}\cdot\varphi_{X|X_i}.
\end{multline}

\subsubsection{Proof of Proposition \ref{prop;17.8.8.10}}

Suppose that there exists a sequence $s_i$
such that
$\sup_{X_i} |s_i|_{h_0}\to\infty$ $(i\to\infty)$.
By the estimate (\ref{eq;17.8.10.20}),
we obtain
$\ell_i:=
 \int_{X_i}|s_i|_{h_0}\varphi_X\to\infty$
as $i\to\infty$.
We set
$u_i:=\ell_i^{-1}s_i$.
They are $G$-invariant sections of 
$\End(E_i)$ on $X_i$
which are self-adjoint with respect to $h_i$.

\begin{lem}
There exists a $G$-invariant $L_1^2$-section $u_{\infty}$ of $\End(E)$
on $X$
such that the following holds:
\begin{itemize}
\item $u_{\infty}\neq 0$.
\item
 A subsequence of
 $\{u_i\}$ is weakly convergent to $u_{\infty}$
 in $L_1^2$ on any compact subset of $X$.
\item
Let $\Phi:\real\times\real\lrarr\openopen{0}{\infty}$
be a $C^{\infty}$-function such that
 $\Phi(y_1,y_2)<(y_1-y_2)^{-1}$
 if $y_1>y_2$.
Then, the following holds:
\begin{equation}
\label{eq;17.8.10.30}
 \sqrt{-1}\int_X\Tr(u_{\infty}\Lambda F(h_0))
+\int_X
 \htilde_0\bigl(
 \Phi(u_{\infty})(\delbar_E u_{\infty}),
 \delbar_E u_{\infty}
 \bigr)\leq 0.
\end{equation}
(See {\rm\S\ref{subsection;17.8.10.21}} 
for the notation $\Phi(u_{\infty})(\delbar_E u_{\infty})$.)
Here, $\htilde_0$
denotes the metric of 
$\End(E)\otimes\Omega^{p,q}$
induced by $h_0$ and $g_X$.
(In the case $p=q=0$, we also use the notation $h_0$.)
\end{itemize}
\end{lem}
\pf
We closely follow the argument of \cite[Lemma 5.4]{Simpson88}.
By Proposition \ref{prop;17.8.10.22},
we have $M(h_{0,i},h_{i})\leq 0$,
and hence
\[
 \sqrt{-1}\int_{X_i}
 \Tr(u_i\Lambda F(h_0))
+
\int_{X_i}
 \htilde_0
 \bigl(
 \ell_i\Psi(\ell_i u_i)(\delbar_E u_i),
 \delbar_E u_i
 \bigr)
\leq 0.
\]
Note that there exists $C_1>0$
such that
$\sup|u_i|<C_1$ for any $i$.
We remark that
$\ell\Psi(\ell\lambda_1,\ell\lambda_2)$
is monotonously increasing for $\ell$,
and convergent to 
$(\lambda_1-\lambda_2)^{-1}$ 
(if $\lambda_1>\lambda_2$)
or 
$\infty$ 
(if $\lambda_1\leq \lambda_2$).
Hence, we obtain
\begin{equation}
\label{eq;17.8.10.50}
 \sqrt{-1}\int_{X_i}
 \Tr(u_i\Lambda F(h_0))
+
\int_{X_i}
 \htilde_0
 \bigl(
 \Phi(u_i)(\delbar_E u_i),
 \delbar_E u_i
 \bigr)
\leq 0.
\end{equation}
There exists $C_2>0$ such that
$\int_{X_i}
 \bigl|u_i\cdot
 \Lambda F(h_0)\bigr|_{h_{0,i}}\leq C_2$
for any $i$.
By (\ref{eq;17.8.10.50}),
there exists $C_3>0$
such that 
the following holds for any $i$:
\[
 \int_{X_i}\bigl|
 \delbar_E u_i
 \bigr|^2_{\htilde_{0}}
<C_3.
\]
Hence, for any compact subset $K$ of $X$,
we obtain that
$u_i$ are bounded in $L_1^2$ on $K$.
By going to a subsequence,
we may assume that
$\{u_i\}$ is weakly convergent in $L_1^2$
on any compact subset of $X$.
Let $u_{\infty}$ denote the weak limit.
Because $|u_i|_{h_{0,i}}\leq C_1$,
we obtain $|u_{\infty}|_{h_0}\leq C_1$.
Because $u_i$ are $G$-invariant,
$u_{\infty}$ is also $G$-invariant.

For any compact subset $Z$ of $X$,
we have
$\int_Z|u_i|_{h_{0,i}}\varphi_X
 \to
 \int_Z|u_{\infty}|_{h_0}\varphi_X$.
Recall $\int_{X_i}|u_i|_{h_{0,i}}=1$ for any $i$
by the construction.
Because $|u_i|_{h_0}\leq C_1$ for any $i$,
there exists a compact subset $Z_0$ such that 
$\int_{X_i\setminus Z_0}|u_i|_{h_{0,i}}\varphi_X\leq 1/2$
for any $i$.
Then, we obtain that
$\int_{Z_0}|u_i|_{h_{0,i}}\varphi_X\geq 1/2$
for any $i$,
which implies
$\int_{Z_0}|u_{\infty}|_{h_0}\varphi_X\geq 1/2$.
In particular,
$u_{\infty}\neq 0$.

For any compact subset $Z$ of $X$,
we have the convergence
\[
 \lim_{i\to\infty}
 \int_{Z}
 \Tr\bigl(u_i\Lambda F(h_0)\bigr)
=\int_{Z}
 \Tr\bigl(u_{\infty}\Lambda F(h_0)\bigr).
\]
Note that $|u_{\infty}|_{h_0}\leq C_1$
and $|u_i|_{h_{0,i}}\leq C_1$ for any $i$.
Hence, for any $\epsilon>0$,
there exists a compact subset $Z(\epsilon)$
such that 
\[
 \int_{X_i\setminus Z(\epsilon)}
 \bigl|
  \Tr\bigl(u_i\Lambda F(h_0)\bigr)
 \bigr|<\epsilon,
\quad
 \int_{X\setminus Z(\epsilon)}
  \bigl|
  \Tr\bigl(u_{\infty}\Lambda F(h_0)\bigr)
 \bigr|<\epsilon.
\]
Hence, we obtain the convergence
\begin{equation}
\label{eq;17.11.27.10}
\lim_{i\to\infty}
 \int_{X_i}
 \Tr\bigl(u_i\Lambda F(h_0)\bigr)
=\int_{X}
 \Tr\bigl(u_{\infty}\Lambda F(h_0)\bigr).
\end{equation}

Take any $\epsilon>0$.
Because of (\ref{eq;17.11.27.10}),
there exists $i_1$
such that the following holds for any $i\geq i_1$:
\[
 \sqrt{-1}
 \int_X\Tr\bigl(u_{\infty}\Lambda F(h_0)\bigr)
+\int_{X_i}
\bigl|
 \Phi^{1/2}(u_i)(\delbar u_i)
  \bigr|_{\htilde_{0,i}}^2
\leq
 \epsilon.
\]
If $i\leq j$,
we have
$\int_{X_i}
 \bigl|
 \Phi^{1/2}(u_j)(\delbar u_j)
 \bigr|_{\htilde_{0,j}}^2
\leq
 \int_{X_j}
 \bigl|
 \Phi^{1/2}(u_j)(\delbar u_j)
 \bigr|_{\htilde_{0,j}}^2$.
Hence, we have the following for any $i<j$:
\[
 \sqrt{-1}
 \int_X\Tr(u_{\infty}\Lambda F(h_0))
+\int_{X_i}
 \bigl|
 \Phi^{1/2}(u_j)(\delbar u_j)
  \bigr|^2_{\htilde_{0,j}}
\leq
 \epsilon.
\]
On the compact space $X_i$,
by applying the argument of Simpson
in the proof of \cite[Lemma 5.4]{Simpson88},
we obtain 
\[
 \sqrt{-1}
 \int_X\Tr(u_{\infty}\Lambda F(h_0))
+\int_{X_i}
 \bigl|
 \Phi^{1/2}(u_{\infty})(\delbar u_{\infty})
 \bigr|_{\htilde_{0}}^2
\leq
 2\epsilon.
\]
Hence, we obtain 
\[
 \sqrt{-1}
 \int_X\Tr(u_{\infty}\Lambda F(h_0))
+\int_X\bigl|
 \Phi^{1/2}(u_{\infty})(\delbar u_{\infty})
 \bigr|^2_{\htilde_0}
\leq
 2\epsilon.
\]
Because this holds for any $\epsilon>0$,
we obtain the desired inequality
(\ref{eq;17.8.10.30}).
\hfill\qed

\vspace{.1in}

We obtain the following lemma
by the argument of Simpson
\cite[Lemma 5.5, Lemma 5.6]{Simpson88}.

\begin{lem}\mbox{{}}
\label{lem;17.8.10.31}
\begin{itemize}
\item
 The eigenvalues of $u_{\infty}$
 are constant.
 We denote them by $\lambda_1,\lambda_2,\ldots,\lambda_{\rank E}$.
\item
Let $\Phi:\real\times\real\lrarr\openopen{0}{\infty}$
be a $C^{\infty}$-function
such that $\Phi(\lambda_i,\lambda_j)=0$
if $\lambda_i>\lambda_j$.
Then,
$\Phi(u_{\infty})(\delbar u_{\infty})=0$.
\hfill\qed
\end{itemize}
\end{lem}

Let $\gamma$ be an open interval between the eigenvalues.
We take a $C^{\infty}$-function
$p_{\gamma}:\real\lrarr \closedopen{0}{\infty}$
such that
$p_{\gamma}(\lambda_i)=1$ if $\lambda_i<\gamma$
and $p_{\gamma}(\lambda_i)=0$ if $\lambda_i>\gamma$.
Set $\pi_{\gamma}:=p_{\gamma}(u_{\infty})$.
By the construction,
we have $\pi_{\gamma}^2=\pi_{\gamma}$ 
and $(\pi_{\gamma})^{\dagger}_{h_0}=\pi_{\gamma}$.
Moreover $\pi_{\gamma}$ is $G$-invariant.

Because
$\delbar\pi_{\gamma}
=dp_{\gamma}(u_{\infty})(\delbar u_{\infty})$,
$\pi_{\gamma}$ is a locally $L_1^2$-section of $\End(E)$,
where
\[
dp_{\gamma}(y_1,y_2):=
 (y_1-y_2)^{-1}(p_{\gamma}(y_1)-p_{\gamma}(y_2)).
\]
We also have 
$\int\bigl|\delbar\pi_{\gamma}\bigr|_{\htilde_0}^2<\infty$.
By setting $\Phi_{\gamma}:=(1-p_{\gamma}(y_2))dp_{\gamma}(y_1,y_2)$.
Because 
$\Phi_{\gamma}(u_{\infty})(\delbar u_{\infty})
=(1-\pi_{\gamma})\delbar\pi_{\gamma}$,
we have $(1-\pi_{\gamma})\delbar\pi_{\gamma}=0$
by Lemma \ref{lem;17.8.10.31}.

According to \cite[\S7]{Uhlenbeck-Yau}
(see also \cite[Proposition 5.8]{Simpson88}
and a more recent work 
\cite[Theorem 0.1.1]{Popovici}\footnote{this reference was informed
by one of the reviewers}),
$\pi_{\gamma}$ determines 
a saturated $\nbigo_X$-submodule $V_{\gamma}$ of $E$
such that $\pi_{\gamma}$ is the orthogonal projection 
of $E$ onto $V_{\gamma}$ outside $Z(V_{\gamma})$.
Because $\pi_{\gamma}$ is $G$-invariant,
we obtain that $V_{\gamma}$ is also $G$-invariant.

The rest of the proof of Proposition \ref{prop;17.8.8.10}
is completely the same
as the proof of \cite[Proposition 5.3]{Simpson88},
which we do not repeat.
\hfill\qed

\subsection{Proof of Theorem \ref{thm;17.8.8.11}}
\label{subsection;17.11.27.1}

By Proposition \ref{prop;17.8.8.10},
we have $\sup_{X_i}|s_i|_{h_{0,i}}<C_1$
for any $i$.
Set $b_i:=e^{s_i}$ on $X_i$.
There exists $C_{20}>0$ such that
the following holds for any $i$:
\[
 \sup_{X_i}
 \bigl(
|b_i|_{h_{0,i}}+|b_i^{-1}|_{h_{0,i}}
\bigr)<C_{20}.
\]
By \cite[Lemma 3.1]{Simpson88},
we obtain the following equality on $X_i$:
\begin{equation}
\label{eq;18.12.26.1}
 \sqrt{-1}\Lambda
 \delbar_{E_i}
 \del_{E_i,h_{0,i}}(b_i)=
-b_i\sqrt{-1}\Lambda F(h_{0,i})^{\bot}
+\sqrt{-1}\Lambda \delbar_{E_{i}}(b_i)
 b_i^{-1}\del_{E_{i},h_{0,i}}(b_i).
\end{equation}
Hence, we obtain the following equality on $X_i$:
\[
 \Delta_{X_i}\Tr(b_i)
 =-\Tr\bigl(
 b_i\sqrt{-1}\Lambda F(h_{0,i})^{\bot}
 \bigr)
-\bigl|
 \delbar_{E_i}(b_i)b_i^{-1/2}
 \bigr|^2_{\htilde_{0,i}}.
\]
We have $\Tr(b_i)-\rank E\geq 0$ on $X_i$
and $\Tr(b_i)-\rank E=0$ on $\del X_i$.
By Green's formula (Lemma \ref{lem;17.11.27.3})
and the inequality $\del_{\nu,i}\Tr(b_i)\leq 0$ on $\del X_i$,
we also have the following:
\[
 \int_{X_i}\Delta_{X_i}\bigl(\Tr(b_i)-\rank E\bigr)
=-\frac{1}{2}\int_{\del X_i}
 \del_{\nu,i}\bigl(
 \Tr(b_i)-\rank E
 \bigr)
\geq 0.
\]
Hence, we obtain the following:
\[
0\leq
 -\int_{X_i}
 \Tr\bigl(b_i\sqrt{-1}\Lambda F(h_{0,i})^{\bot}\bigr)
-\int_{X_i}
 \bigl| 
 \delbar(b_i)b_i^{-1/2}
 \bigr|^2_{\htilde_{0,i}}.
\]
There exist positive constants $C_{21},C_{22}$
such that the following holds for any $i$:
\[
\int_{X_i}
 \bigl|
\delbar(b_i)b_i^{-1/2}
 \bigr|^2_{\htilde_{0,i}}
\leq
 C_{21}\int_X\bigl|\Lambda F(h_0)^{\bot}\bigr|_{h_0}
\leq C_{22}.
\]
Hence, there exists $C_{23}>0$
such that the following holds
for any $i$:
\begin{equation}
\label{eq;17.11.14.2}
 \int_{X_i}\bigl|
 \delbar b_i
 \bigr|^2_{\htilde_{0,i}}
\leq C_{23} 
\end{equation}

Let us prove that the $C^1$-norm of $b_i$
are locally bounded 
by using the argument in the proof of
\cite[Theorem 2.10]{Ni1}.
Because $b_i$ is self-adjoint with respect to 
$h_{0,i}$,
$|\delbar b_i|_{\htilde_{0,i}}$
and
$|\del_{h_{0,i}}b_i|_{\htilde_{0,i}}$
are equal.
Hence, it is enough to study
$|\del_{h_{0,i}} b_i|_{\htilde_{0,i}}$.
Let $P$ be any point of $X$.
Take a holomorphic coordinate neighbourhood 
$(X_P;z_1,\ldots,z_n)$ around $P$ in $X$.
There exists $i_1(P)$ such that $X_P\subset X_i$
for any $i\geq i_1(P)$.
Let $X_P'$ be a relatively compact neighbourhood
of $P$ in $X_P$.

\begin{lem}
\label{lem;18.12.29.1}
$|\del_{h_{0,i}} b_i|_{\htilde_{0,i}}$ are bounded
on $X_P'$.
\end{lem}
\pf
Let $\htilde_i$ denote the metric
of $E_i\otimes\Omega^{1,0}_{X_i}$
induced by $h_i$ and $g_{X_i}$.
We set $b_{1,i}:=b_i^{-1}$,
i.e., 
$b_{1,i}$ is determined by $h_{0,i}=h_ib_{1,i}$.
Because
\[
 b_{1,i}^{-1}\del_{E_i,h_i}b_{1,i}
=\del_{E_i,h_{0,i}}-\del_{E_i,h_i}
=-b_i^{-1}\del_{E_i,h_{0,i}}b_i,
\]
and because $|b_{1,i}|_{h_i}$ are uniformly bounded,
it is enough to prove that
$\bigl|b_{1,i}^{-1}\del_{E_i,h_i}(b_{1,i})\bigr|_{\htilde_{i}}$
are uniformly bounded.
From (\ref{eq;17.11.14.2}),
there exist $C_{24}>0$
such that
\begin{equation}
\label{eq;18.12.25.10}
 \int_{X_i}\bigl|
 b_{1,i}^{-1}\del_{E,h_i}(b_{1,i})
 \bigr|^2_{\htilde_i}
<C_{24}.
\end{equation}
By Lemma \ref{lem;17.10.8.1}
and Lemma \ref{lem;17.10.8.2},
there exists $C_{P,k}>0$ $(k=1,2)$
such that the following holds
for any $i\geq i_1(P)$ on $X_P$:
\[
 \Delta_X\bigl|
 b_{1,i}^{-1}\del_{E_i,h_i}(b_{1,i})
 \bigr|_{\htilde_i}^2
\leq
 C_{P,1}\bigl|
 b_{1,i}^{-1}\del_{E_i,h_i}(b_{1,i})
 \bigr|_{\htilde_i}
+C_{P,2}
\bigl|
 b_{1,i}^{-1}\del_{E_i,h_i}(b_{1,i})
 \bigr|_{\htilde_i}^2.
\]
Hence, there exists $C_{P,k}>0$ $(k=3,4)$
such that the following holds for any $i\geq i_1(P)$
on $X_P$
\[
 \Delta_X\bigl|
 b_{1,i}^{-1}\del_{E_i,h_i}(b_{1,i})
 \bigr|_{\htilde_i}^2
-C_{P,4}
\bigl|
 b_{1,i}^{-1}\del_{E_i,h_i}(b_{1,i})
 \bigr|_{\htilde_i}^2
\leq
 C_{P,3}.
\]
Then, by using \cite[Theorem 9.20]{Gilbarg-Trudinger}
with the $L^2$-boundedness (\ref{eq;18.12.25.10}),
we obtain the uniform boundedness of
$\bigl|b_{1,i}^{-1}\del_{E_i,h_i}b_{1,i}\bigr|_{\htilde_i}^2$.
\hfill\qed

\begin{rem}
\label{rem;19.1.2.1}
In the proof of Lemma {\rm\ref{lem;18.12.29.1}},
if the curvature of $(X,g_X)$ is flat,
and if $F(h_0)=0$,
then 
$\bigl|
b_{1,i}^{-1} \del_{E_i,h_{i}}b_{1,i}
 \bigr|_{\htilde_{0,i}}^2$
is subharmonic,
as proved in {\rm\cite{Donaldson-boundary-value}}.
More generally,
if $F(h_0)=0$
and if $\sqrt{-1}\Lambda F(g_X)$ is non-positive,
then 
$\bigl|
b_{1,i}^{-1} \del_{E_i,h_{i}}b_{1,i}
 \bigr|_{\htilde_{0,i}}^2$
is subharmonic
which follows from Lemma {\rm\ref{lem;17.10.8.1}}
and Lemma {\rm\ref{lem;17.10.8.2}}.
\hfill\qed
\end{rem}

On any compact subset $K$ of $X$,
the equation (\ref{eq;18.12.26.1}),
the $C^0$-boundedness of $b_i$ and $b_i^{-1}$
and the $C^1$-boundedness of $b_i$
imply the $L_2^p$-boundedness of
$b_i$ for any $p$.
Going to a subsequence,
we may assume that $\{b_i\}$ is weakly convergent in $L_2^p$
for any $p\geq 1$ on any compact subset of $X$.
Then, we obtain the limit $h=\lim h_i$,
which is a Hermitian-Einstein metric of $(E,\delbar_E)$
with the desired property.
Thus, we obtain Theorem \ref{thm;17.8.8.11}.
\hfill\qed

\subsection{Proof of Proposition \ref{prop;17.8.10.40}
and Proposition \ref{prop;18.1.20.30}}
\label{subsection;17.11.30.1}

\subsubsection{Proof of Proposition \ref{prop;17.8.10.40}}
\label{subsection;18.1.20.31}

We use the notation in \S\ref{subsection;17.8.10.50}.
We know that
$h$ and $h_0$ are mutually bounded.
Let $b$ be determined by
$h=h_0\cdot b$
as in Theorem \ref{thm;17.8.8.11}.
Let $b_1:=b^{-1}$,
i.e.,
$b_1$ is determined by $h_0=h\cdot b$.
Because
\begin{equation}
\label{eq;17.10.8.3}
 b_1^{-1}\del_{E,h}b_1=\del_{E,h_0}-\del_{E,h}
=-b^{-1}\del_{E,h_0}b,
\end{equation}
we obtain
$\int_Y\bigl|b_1^{-1}\del_{E,h}b_1\bigr|_{\htilde}^2<\infty$.
For any $\epsilon>0$,
there exists a compact subset $K_{\epsilon}\subset Y$
such that the following holds.
\begin{itemize}
\item
$\int_{Y\setminus K_{\epsilon}}
 |b_1^{-1}\del_{E,h}b_1|_{\htilde}^2<\epsilon$.
\item
 $\bigl|\del_{E,h_0} \Lambda F(h_0)^{\bot}\bigr|_{\htilde}<\epsilon$
 on $Y\setminus K_{\epsilon}$.
\item
 $\bigl| F(h_0)^{\bot}\bigr|_{\htilde}<\epsilon$
 on $Y\setminus K_{\epsilon}$.
\end{itemize}

By Lemma \ref{lem;17.10.8.1}
and Lemma \ref{lem;17.10.8.2},
there exists $C_1>0$,
which is independent of $(\epsilon,K_{\epsilon})$,
such that the following holds
on $Y\setminus K_{\epsilon}$:
\[
 \Delta_X\bigl|
 b_1^{-1}\del_{E,h}b_1
 \bigr|_{\htilde}^2
\leq
 C_1\epsilon
 |b_1^{-1}\del_{E,h}b_1|_{\htilde}
+C_1(\epsilon+1)
|b_1^{-1}\del_{E,h}b_1|_{\htilde}^2
\leq
 \frac{C_1\epsilon}{2}
+C_1\bigl(3\epsilon+1\bigr)
 |b_1^{-1}\del_{E,h}b_1|_{\htilde}^2.
\]
Hence, there exists a constant $C_2>0$
which is independent of $(\epsilon,K_{\epsilon})$,
such that the following holds
on $Y\setminus K_{\epsilon}$:
\[
 \bigl(
 \Delta_X-C_2
 \bigr)|b_1^{-1}\del_{E,h}b_1|^2_{\htilde}
\leq C_2\epsilon.
\]
By using \cite[Theorem 9.20]{Gilbarg-Trudinger},
we obtain the following.
\begin{lem}
For any $\epsilon>0$,
there exists a compact subset
$K_{\epsilon}'\subset Y$
such that
the following holds
on $Y\setminus K_{\epsilon}'$:
\[
 \sup_{Y\setminus K_{\epsilon}'}
 \bigl|
 b_1^{-1}\del_{E,h}b_1
 \bigr|_{\htilde}^2
\leq
  \epsilon.
\]
\hfill\qed
\end{lem}

By (\ref{eq;17.10.8.3}),
for any $\epsilon>0$,
there exists a compact subset
$K''_{\epsilon}\subset Y$
such that
$\bigl|
 b^{-1}\del_{E,h_0}b
 \bigr|_{\htilde_0}\leq \epsilon$
on $Y\setminus K''_{\epsilon}$.
According to \cite[Lemma 3.1]{Simpson88},
the following equality holds:
\begin{equation}
\label{eq;17.10.13.300}
\sqrt{-1}\Lambda\delbar_{E}\del_{E,h_0}(b)
=-b\sqrt{-1}\Lambda F(h_0)^{\bot}
+\sqrt{-1}\Lambda
 \delbar_E(b)b^{-1}\del_{E,h_0}(b).
\end{equation}
Hence, for any $\epsilon>0$,
there exists a compact subset $K^{(3)}_{\epsilon}\subset Y$
such that
$\bigl|
\sqrt{-1}\Lambda\delbar_E\del_{E,h_0}(b)
 \bigr|_{h_0}\leq \epsilon$
on $Y\setminus K^{(3)}_{\epsilon}$.

Take a large $p$.
We obtain that
the $L_2^p$-norm  of $b-\id$ on
the disc with radius $r_0$ centered at $P\in Y$
goes to $0$ as $P$ goes to $\infty$.
By using (\ref{eq;17.10.13.300}),
we obtain that 
the $L_3^p$-norm of $b-\id$ on
the disc with radius $r_0$ centered at $P\in Y$
goes to $0$ as $P$ goes to $\infty$.
Then, we obtain the claim of Proposition \ref{prop;17.8.10.40}.
\hfill\qed

\subsubsection{Proof of Proposition \ref{prop;18.1.20.30}}

We use the notation in \S\ref{subsection;18.1.20.31}.
There exists $C_{10}>0$ such that 
the following holds:
\begin{itemize}
\item
$\int_Y|b_1^{-1}\del_{E,h}b_1|_{\htilde}^2\leq C_{10}$.
\item
$\bigl|
  \del_{E,h_0}\Lambda F(h_0)^{\bot}
 \bigr|_{\htilde}\leq C_{10}$
and 
$\bigl|
 F(h_0)^{\bot}
 \bigr|_{\htilde}\leq C_{10}$ hold on $Y$.
\end{itemize}
By Lemma \ref{lem;17.10.8.1}
and Lemma \ref{lem;17.10.8.2},
there exists $C_{11}>0$
such that the following holds on $Y$:
\[
 \Delta_X\bigl|
 b_1^{-1}\del_{E,h}b_1
 \bigr|_{\htilde}^2
\leq
 C_{11}
 |b_1^{-1}\del_{E,h}b_1|_{\htilde}
+2C_{11}|b_1^{-1}\del_{E,h}b_1|_{\htilde}^2
\leq
 \frac{C_{11}}{2}
+\frac{5C_{11}}{2}
 |b_1^{-1}\del_{E,h}b_1|_{\htilde}^2.
\]
Again, by using \cite[Theorem 9.20]{Gilbarg-Trudinger},
we obtain that there exists a constant $C_{12}>0$
such that
\[
 \sup_Y
 |b_1^{-1}\del_{E,h}b_1|_{\htilde}^2
\leq
 C_{12}.
\]
Hence, there exists $C_{13}>0$ such that
$|b^{-1}\del_{E,h_0}b|_{\htilde_0}^2
\leq
 C_{13}$ on $Y$.
By using (\ref{eq;17.10.13.300}),
we obtain that 
there exists a compact subset $K\subset Y$
such that 
the $L_3^p$-norms of $b-\id$ 
on the disc with radius $r_0$ centered at $P\in Y\setminus K$
are bounded.
Hence we obtain the boundedness of $b$.
Thus, the proof of Proposition \ref{prop;18.1.20.30}
is completed.
\hfill\qed

\subsection{Proof of Proposition \ref{prop;17.8.12.1}}
\label{subsection;17.11.30.2}

We closely follow the argument of Simpson
in the proof of \cite[Proposition 3.5, Lemma 7.4]{Simpson88}.
We shall prove the following inequalities:
\begin{equation}
\label{eq;17.11.28.20}
\int_{X}
 \Tr\bigl((F(h)^{\bot})^2\bigr)\omega_X^{\dim X-2}
\leq
 \int_{X}
 \Tr\bigl((F(h_{0})^{\bot})^2\bigr)\omega_X^{\dim X-2},
\end{equation}
\begin{equation}
\label{eq;17.11.28.21}
\int_{X} 
\Tr\bigl((F(h_0)^{\bot})^2\bigr)\omega_X^{\dim X-2}
\leq
\int_{X} 
\Tr\bigl((F(h)^{\bot})^2\bigr)\omega_X^{\dim X-2}.
\end{equation}

We may assume that 
$X_i=\{P\in X\,|\,\varrho(P)\leq a_i\}$
for a sequence $a_i>0$.
On $X_i$,
we set $f_i:=1-a_i^{-1}\varrho$.
We have $f_i>0$ on $X_i\setminus \del X_i$,
and $f_i=0$ on $\del X_i$.
Let $\omega_{X_i}$ denote the K\"ahler form of $X_i$.
In the following,
when we are given a Hermitian metric $k$ on $E$,
the induced metrics $E\otimes\Omega^{p,q}$
are also denoted by $k$
for simplicity of the description.

\subsubsection{Proof of (\ref{eq;17.11.28.20})}

Let $\nbigp_i$ denote the space of $C^{\infty}$ Hermitian metrics $k$
of $E_i$ such that $k_{|\del X_i}=h_{0|\del X_i}$.
Take any $k_1,k_2\in\nbigp_i$.
Let $s$ be the endomorphism $s$ of $\End(E_i)$
determined by $k_2=k_1\cdot e^{s}$
which is self-adjoint with respect to both of $k_1$ and $k_2$.
We put
$\tau_i:=-2\sqrt{-1}\del\delbar f_i
=2\sqrt{-1}a_i^{-1}\del\delbar \varrho$ on $X_i$.
Let $n:=\dim X$.
We set
\[
 M_i(k_1,k_2):=
 \sqrt{-1}\int_{X_i}
 \Tr\bigl(
 sF(k_1)
 \bigr)\tau_i\omega_{X_i}^{n-2}
-\sqrt{-1}
 \int_{X_i}
 \Tr\Bigl(
 \Psi(s)(\delbar s)\cdot
\del_{k_1}s
 \Bigr)
 \tau_i\omega_{X_i}^{n-2}.
\]
The following is an analogue of \cite[Lemma 7.2]{Simpson88},
and the proof is also similar.
\begin{lem}
We have the following equality:
\begin{equation}
\label{eq;18.12.25.20}
 M_i(k_1,k_2)=
\int_{X_i}
 f_i\bigl(
 \Tr(F(k_1)^2)
-\Tr(F(k_2)^2)
 \bigr)\omega_{X_i}^{n-2}.
\end{equation}
\end{lem}
\pf
We give only an indication.
By the argument in the proof of Proposition \ref{prop;17.8.10.2}.
we can prove $M_i(k_1,k_2)=M_i(k_1,k_3)+M_i(k_3,k_2)$
for any $k_j\in\nbigp$ $(j=1,2,3)$.
Moreover, we can obtain the following
for any $k\in\nbigp_i$
and any endomorphism $s$ of $E_i$
which is self-adjoint with respect to $k$.
\[
 \frac{d}{dt}M_i(k,ke^{ts})_{|t=0}=
 \sqrt{-1}\int_{X_i}\Tr(s F(k))\tau_i\omega_{X_i}^{n-2}.
\]
Let $M_i'(k_1,k_2)$ denote the right hand side of
(\ref{eq;18.12.25.20}).
For any $k\in \nbigp$
and any endomorphism $s$ of $E_i$
which is self-adjoint with respect to $k$,
the following holds:
\[
 \frac{d}{dt}M_i'(k,ke^{ts})_{|t=0}
=-2\int_{X_i}f_i\Tr\bigl(\delbar\del_{E_i,k}(s)F(k)\bigr)
 \tau_i\omega_{X_i}^{n-2}.
\]
Because of $f_{i|\del X_i}=0$ 
and $s_{|\del X_i=0}$,
we obtain the following by the Stokes formula:
\[
-2\int_{X_i}f_i\Tr\bigl(\delbar\del_{E_i,k}(s)F(k)\bigr)
 \tau_i\omega_{X_i}^{n-2}
=2\int_{X_i}\delbar(f_i)
 \Tr\bigl(\del_{E_i,k}(s)F(k)\bigr)\omega_{X_i}^{n-2}
=2\int_{X_i}\del\delbar(f_i)
 \Tr\bigl(sF(k)\bigr)\omega_{X_i}^{n-2}.
\]
Hence, we obtain $\frac{d}{dt}M_i(k,ke^{ts})_{|t=0}
=\frac{d}{dt}M_i'(k,ke^{ts})_{|t=0}$,
and the claim of the lemma.
\hfill\qed

\vspace{.1in}

Let $h_{0,i}$, $h_i$ and $s_i$
be as in \S\ref{subsection;17.11.15.1}.
\begin{lem}
There exists a constant $C_1>0$
such that the following holds for any $i$:
\begin{equation}
 \label{eq;17.8.11.100}
 M_i(h_{0,i},h_i)
\geq
-\frac{C_1}{a_i}
 \int_{X_i}
 \bigl|
 \bigl(
 F(h_{0})^{\bot}\cdot\del\delbar\varrho
 \bigr)_{|X_i}
 \bigr|_{h_{0,i}}
-\frac{C_1}{a_i}
 \int_{X_i}
 \bigl|
 \delbar_Es_i
 \bigr|_{h_{0,i}}^2.
\end{equation}
\end{lem}
\pf
We consider 
$M_i(h_{0,i},h_{0,i}e^{ys_i})$ for $0\leq y\leq 1$.
We have
\[
 \frac{d^2}{dy^2}
 M_i(h_{0,i},h_{0,i}e^{ys_i})
=\sqrt{-1}
 \int_{X_i}
 \Tr\bigl(
 s_i\delbar\del_{h_{0,i}e^{ys_i}}(s_i)
 \bigr)\tau_i\omega_{X_i}^{n-2}
=-\sqrt{-1}
 \int_{X_i}
 \Tr\bigl(
 \delbar_{E_i}(s_i)\cdot
 \del_{h_{0,i}e^{ys_i}}(s_i)
 \bigr)
 \tau_i\omega_{X_i}^{n-2}.
\]
Take any point $P\in X_i\setminus \del X_i$.
We take a holomorphic coordinate system
$(z_1,\ldots,z_n)$ around $P$ such that
(i) $z_i(P)=0$,
(ii) $g_{X|P}=(\sum dz_i\,d\zbar_i)_{|P}$,
(iii) $\tau_{i|P}=\sum b_i\,(dz_i\,d\zbar_i)_{|P}$
for some $(b_1,\ldots,b_n)\in\cnum^n$.
There exists $C_2>0$ 
which is independent of $i$ and $P$,
such that $|b_p|<C_2/a_i$.
We express $\delbar_{E_i}(s_i)$
and $\del_{h_{0,i}e^{ys_i}}(s_i)$
by using local sections $A_p$ 
of $\End(E_i)$:
\[
 \delbar_{E_i}(s_i)
=\sum A_p\,d\zbar_p,
\quad
 \del_{h_{0,i}e^{ys_i}}(s_i)
=\sum (A_p)^{\dagger}_{h_{0,i}e^{ys_i}}\,
 dz_p.
\]
At $P$,
we obtain
\[
 \Tr\bigl(
 \delbar_{E_i}(s_i)
 \del_{h_{0,i}e^{ys_i}}(s_i)
 \bigr)\cdot
 \tau_i
 \omega_{X_i}^{n-2}
=\sum_{p\neq q}
 \bigl|
 A_p
 \bigr|^2
 _{h_{0,i}e^{ys_i}}
 \cdot
  b_q
 (dz_p\,d\zbar_p)
 (dz_q\,d\zbar_q)
 \omega_X^{n-2}
=C_{3}
\sum_{p\neq q}
 \bigl|
 A_p
 \bigr|^2
 _{h_{0,i}e^{ys_i}}
 \cdot
  b_q\cdot
\omega_X^{n}.
\]
Here, $C_3$ depends only on 
the dimension $n$.
Because $0\leq y\leq 1$,
and because of Proposition \ref{prop;17.8.8.10},
there exist positive constants $C_5$ and $C_6$,
which are independent of $i$ and $P$,
such that the following holds:
\[
\Bigl|
\sum_{p\neq q}
 \bigl|
 A_p
 \bigr|^2
 _{h_{0,i}e^{ys_i}}
 \cdot
  b_q
\Bigr|
\leq 
\frac{C_5}{a_i}
 \bigl|
 \delbar_{E_i}(s_i)
 \bigr|^2_{h_{0,i}e^{ys_i}}
\leq
\frac{C_5}{{a_i}}\Bigl(
  \bigl|
 \delbar_{E_i}(s_i)
 \bigr|^2_{h_{0,i}}
+ \bigl|
 \delbar_{E_i}(s_i)
 \bigr|^2_{h_{0,i}e^{s_i}}
 \Bigr)
\leq
\frac{C_6}{a_i}\bigl|
 \delbar_{E_i}s_i
 \bigr|^2_{h_{0,i}}.
\]
Hence, we obtain the following:
\[
 \frac{d^2}{dy^2}
 M_{i}(h_{0,i},h_{0,i}e^{ys_i})
\geq
 -\frac{C_6}{a_i}
 \int_{X_i}
 \bigl|
 \delbar_{E_i}s_i
 \bigr|_{h_{0,i}}^2.
\]
Because $\Tr(s_i)=0$,
we obtain
\[
 \frac{d}{dy}
 M_i(h_{0,i},h_{0,i}e^{ys_i})_{|y=0}
=\sqrt{-1}\int_{X_i}
 \Tr(s_iF(h_{0,i}))
 \tau_i\omega_X^{n-2}
=\sqrt{-1}\int_{X_i}
 \Tr(s_iF(h_{0,i})^{\bot})
 \tau_i\omega_X^{n-2}.
\]
Hence, there exists a positive constant $C_7$,
which is independent of $i$,
such that the following holds:
\[
 \frac{d}{dy}
 M_{i}(h_{0,i},h_{0,i}e^{ys_i})
\geq
-\frac{C_7}{a_i}
 \int_{X_i}
 |s_i|_{h_{0,i}}
 \cdot
 \bigl|
 F(h_{0,i})^{\bot}\del\delbar\varrho
 \bigr|_{h_{0,i}}
-\frac{C_7}{a_i}
 \int_{X_i}
 \bigl|\delbar_{E_i}s_i\bigr|_{h_{0,i}}^2.
\]
Then, we obtain the inequality
(\ref{eq;17.8.11.100}).
\hfill\qed

\vspace{.1in}

Hence, we have a positive constant $C_{10}$, $C_{11}$
which is independent of $i$
such that the following holds:
\[
 \int_{X_i}
 f_i
 \Tr\bigl(F(h_i)^2\bigr)
\leq
 \int_{X_i}
 f_i
 \Tr\Bigl(
 F(h_{0,i})^2
 \Bigr)
+\frac{C_{10}}{a_i}
+\frac{C_{11}}{a_i}
 \int_{X_i}
 \bigl|
 F(h_0)^{\bot}\del\delbar\varrho
 \bigr|_{h_{0,i}}.
\]
Because
$\Tr F(h_{0,i})=\Tr F(h_i)$ on $X_i$,
we have
\[
 \int_{X_i}
 \Tr\bigl(f_i(F(h_i)^{\bot})^2\bigr)
\leq
 \int_{X_i}
 \Tr\Bigl(
 f_i(F(h_{0,i})^{\bot})^2
 \Bigr)
+\frac{C_{10}}{a_i}
+\frac{C_{11}}{a_i}
 \int_{X_i}
 \bigl|
 F(h_0)^{\bot}\del\delbar\varrho
 \bigr|_{h_{0,i}}.
\]
By the Hermitian-Einstein condition for $h_j$
and the non-negativity $f_j\geq 0$,
we have the following for any $i\leq j$:
\begin{multline}
 \int_{X_i}
 f_j
 \Tr\bigl((F(h_j)^{\bot})^2\bigr)
\leq
 \int_{X_j}
 f_j
 \Tr\bigl((F(h_j)^{\bot})^2\bigr)
 \\
\leq
 \int_{X_j}
 f_j
 \Tr\Bigl(
 (F(h_{0,j})^{\bot})^2
 \Bigr)
+\frac{C_{10}}{a_j}
+\frac{C_{11}}{a_j}
 \int_{X_j}
 \bigl|
 F(h_0)^{\bot}\del\delbar\varrho
 \bigr|_{h_{0,j}}.
\end{multline}
We have the convergences
\[
 \lim_{j\to\infty}
 \int_{X_i}
 f_j
 \Tr\bigl((F(h_j)^{\bot})^2\bigr)
=\int_{X_i}
 \Tr\bigl((F(h)^{\bot})^2\bigr),
\]
\[
 \lim_{j\to\infty}
 \int_{X_j}
 f_j
 \Tr\bigl((F(h_{0,j})^{\bot})^2\bigr)
=\int_{X}
 \Tr\bigl((F(h_{0})^{\bot})^2\bigr).
\]
Hence, we obtain
\[
\int_{X_i}
 \Tr\bigl((F(h)^{\bot})^2\bigr)
\leq
 \int_{X}
 \Tr\bigl((F(h_{0})^{\bot})^2\bigr).
\]
By taking the limit $i\to\infty$,
we obtain (\ref{eq;17.11.28.20}).

\subsubsection{Proof of (\ref{eq;17.11.28.21})}

We extend the function $f_i$ on $X_i$
to the function $\ftilde_i$ on $X$
by setting $0$ outside $X_i$.
We obtain the current
$\tautilde_i:=-2\sqrt{-1}\del\delbar f_i$ on $X$.

Let $k_1,k_2$ be $C^{\infty}$ Hermitian metrics of $E$
such that $\det(k_i)=\det(h_0)$.
We have the endomorphism $s$ of $E$
determined by $k_2=k_1e^s$.
We set 
\[
 \Mtilde_i(k_1,k_2):=
 \sqrt{-1}
\int_{X}\Tr\bigl(s F(k_1)\bigr)\tautilde_i\omega_X^{n-2}
-\sqrt{-1}
 \int_X\Tr\Bigl(
 \Psi(s)(\delbar s)\cdot \del_{k_1}s
 \Bigr)\cdot
 \tautilde_i\omega_X^{n-2}.
\]
As proved in \cite[Lemma 7.2]{Simpson88},
the following holds:
\[
 \Mtilde_i(k_1,k_2)=
\int_{X}\ftilde_i
 \Tr\Bigl(
 F(k_1)^2
-F(k_2)^2
 \Bigr)\omega_X^{n-2}
=
 \int_{X}\ftilde_i
 \Tr\Bigl(
 (F(k_1)^{\bot})^2
-(F(k_2)^{\bot})^2
 \Bigr)\omega_X^{n-2}.
\]

\begin{lem}
There exists a positive constant $C_{20}$
such that the following holds for any $i$:
\begin{equation}
 \label{eq;17.8.11.101}
 \Mtilde_i(h_{0},h)
\leq
\frac{C_{20}}{a_i}
 \Bigl(
\int_{X_i}
 \bigl|
 F(h_0)^{\bot}\cdot \del\varrho
 \bigr|^2_{h_0}
 \Bigr)^{1/2}
\Bigl(
\int_{X_i}
 \bigl|
 \delbar_Es
 \bigr|^2_{h_0}
 \Bigr)^{1/2}
+\frac{C_{20}}{a_i}
 \int_{X}
 \bigl|
 \delbar_Es
 \bigr|_{h_{0}}^2.
\end{equation}
Here, $s$ is determined by
$h=h_0e^s$.
\end{lem}
\pf
Recall $\Tr(s)=0$.
The following holds:
\begin{multline}
\Bigl|
\int_{X}\Tr(s F(h_{0}))
 \tautilde_i\omega_X^{n-2}
\Bigr|
=
\Bigl|
\int_{X}\Tr(s F(h_{0})^{\bot})
 \tautilde_i\omega_X^{n-2}
\Bigr|
=2\Bigl|
 \int_{X_i} 
 \del\delbar\Tr(s F(h_0)^{\bot})
 (1-a_i^{-1}\varrho)\omega_X^{n-2}
 \Bigr|
 \\
=2\Bigl|
 \frac{1}{a_i}
 \int_{X_i}
 \Tr\bigl(\delbar s
\cdot F(h_0)^{\bot}\cdot \del\varrho\bigr)
 \omega_X^{n-2}
 \Bigr|.
\end{multline}
Hence, we obtain
\[
\Bigl|
\int_{X}\Tr(s F(h_{0}))
 \tautilde_i\omega_X^{n-2}
\Bigr|
\leq
 \frac{2}{a_i}
 \Bigl(
 \int_{X_i}
 \bigl|\delbar s\bigr|_{h_0}^2
 \Bigr)^{1/2}
 \Bigl(
 \int_{X_i}
 \bigl|F(h_0)^{\bot}\del\varrho
 \bigr|_{h_0}^2
 \Bigr)^{1/2}.
\]
By definition,
we obtain
\begin{multline}
-\sqrt{-1}
 \int_X\Tr\bigl(
 \Psi(s)\delbar s\del_{h_0}s
 \bigr)
 \tautilde_i\omega_X^{n-2}
=-\sqrt{-1}\int_{X_i}
 -\delbar\del\Tr(\Psi(s)\delbar s\del_{h_0}s)
 (-2\sqrt{-1}f_i)\omega_X^{n-2}
\\
 =\sqrt{-1}\int_{\del X_i}
 \Tr(\Psi(s)\delbar s\del_{h_0}s)
 \bigl(-2\sqrt{-1}(\delbar f_i)\bigr)\omega_X^{n-2}
-\sqrt{-1}\int_{X_i}
 \Tr(\Psi(s)\delbar s\del_{h_0}s)
 (-2\sqrt{-1}\del\delbar f_i)\omega_X^{n-2}.
\end{multline}

\begin{lem}
We have
\[
 \sqrt{-1}\int_{\del X_i}
 \Tr(\Psi(s)\delbar s\del_{h_0}s)
 \bigl(-2\sqrt{-1}(\delbar f_i)\bigr)\omega_X^{n-2}
\leq 0.
\]
\end{lem}
\pf
It is enough to prove that
the integrand is positive at each $P\in \del X_i$.
Note that $X_i=\{f_i\geq 0\}=\{\varrho-a_i\leq 0\}$.
We take a holomorphic coordinate neighbourhood
$(z_1,\ldots,z_n)=(x_1+\sqrt{-1}y_1,\ldots,x_n+\sqrt{-1}y_n)$
around $P$
such that
(i) $z_i(P)=0$,
(ii) $g_{X|P}=\sum dz_p\,d\zbar_p$,
(iii) $(d\varrho)_P=b\,dx_1$ for $b>0$.
We may also assume the following:
\[
 \sqrt{-1}\Tr\bigl(\Psi(s)\delbar s\,\del_{h_0}s\bigr)
=-\sum c_p\,\sqrt{-1}dz_p\,d\zbar_p
\]
for some $c_p\geq 0$ $(p=1,\ldots,n)$.
We have 
$(-2\sqrt{-1}(\delbar f_i))_P
=-a_i^{-1}(-2\sqrt{-1}\,\delbar \varrho)_P
=a_i^{-1}b(\sqrt{-1}d\zbar_1)_P$.
As a form on the tangent space
$T_P\del X_i$,
we have
\begin{multline}
 \Bigl(
 \sqrt{-1} 
 \Tr(\Psi(s)\delbar s\del_{h_0}s)
 (-2\sqrt{-1}(\delbar f_i))\omega^{n-2}
\Bigr)_P
= \\
 -\Bigl(
 \sum_{p=2}^n c_p \sqrt{-1}dz_p\,d\zbar_p
  \Bigr)_{P}
\cdot
 \Bigl(
 a_i^{-1}b dy_1
 \Bigr)_{P}
 \cdot
 \Bigl(
 \sum_{i=2}^n \frac{\sqrt{-1}}{2}dz\,d\zbar_i
 \Bigr)^{n-2}_{P}.
\end{multline}
Then, the claim of the lemma follows.
\hfill\qed

\vspace{.1in}

Because $s$ and $\del\delbar\varrho$ are bounded,
there exists a constant $C_{21}>0$ such that
the following holds for any $i$:
\[
\Bigl|
\sqrt{-1}\int_{X_i}
 \Tr(\Psi(s)\delbar s\del_{h_0}s)
 (-2\sqrt{-1}\del\delbar f_i)\omega^{n-2}
\Bigr|
\leq
 \frac{C_{21}}{a_i}
 \int_{X}
 \bigl|\delbar s\bigr|_{h_0}^2.
\]
Hence, we obtain
\[
 \sqrt{-1}
 \int_X\Tr\bigl(
 \Psi(s)\delbar s\del_{h_0}s
 \bigr)
 \tautilde_i\omega^{n-2}
\leq
\frac{C_{21}}{a_i}
 \int_{X}
 \bigl|\delbar s\bigr|_{h_0}^2.
\]
Thus, we obtain (\ref{eq;17.8.11.101}).
\hfill\qed

\vspace{.1in}
There exist positive constants $C_{30}$ and $C_{31}$
such that the following holds for any $i$:
\[
 \int_{X_i}
 f_i
 \Tr\bigl(
 (F(h_{0})^{\bot})^2
 \bigr)
\leq
 \int_{X_i}
 f_i\Tr\bigl((F(h)^{\bot})^2\bigr)
+\frac{C_{30}}{a_i}
+\frac{C_{31}}{a_i}
 \Bigl(
\int_{X_i}
 \bigl|
 F(h_0)^{\bot}\del\varrho
 \bigr|_{h_0}^2
\Bigr)^{1/2}.
\]
By the theorem of Lebesgue,
we have the following convergence:
\[
\lim_{i\to\infty}
 \int_{X_i}
 f_i
 \Tr\bigl(
 (F(h_{0})^{\bot})^2
 \bigr)
=
 \int_X
 \Tr\bigl(
 (F(h_{0})^{\bot})^2
 \bigr).
\]
Note that $\Tr\bigl((F(h)^{\bot})^2\bigr)$
is positive by the Hermitian-Einstein condition,
and $f_i$ is monotonously increasing for $i$.
Hence, we have the following convergence:
\[
\lim_{i\to\infty}
 \int_{X_i}
 f_i\Tr\bigl((F(h)^{\bot})^2\bigr)
=
 \int_X\Tr\bigl(
(F(h)^{\bot})^2\bigr).
\]
Hence, we obtain (\ref{eq;17.11.28.21}).
Thus the proof of Proposition \ref{prop;17.8.12.1}
is completed.
\hfill\qed

\subsection{Proof of Proposition \ref{prop;17.11.30.3}}
\label{subsection;17.11.30.4}

Let $h_1$ and $h_2$ be metrics
as in Proposition \ref{prop;17.11.30.3}.
For the proof of the proposition,
we may assume that $h_1$ is as in 
Theorem \ref{thm;17.8.8.11}.
Namely, let $b_1$ be determined by 
$h_1=h_0b_1$,
then $\delbar b_1$ is $L^2$.

By Proposition \ref{prop;17.11.30.5},
we have the decomposition
$(E,\delbar_E)=\bigoplus_{j=1}^m (E_j,\delbar_{E_j})$
such that 
(i) the decomposition is 
orthogonal with respect to $h_i$ $(i=1,2)$,
(ii) $h_{2|E_j}=a_j\cdot h_{1|E_j}$ for some $a_j>0$.
Let $\pi_i$ denote the projection onto $E_i$
with respect to the decomposition.
Let $\pi_{i,h_0}^{\dagger}$
denote the adjoint of $\pi_i$
with respect to $h_0$.
Note that $\pi_i$ are bounded
with respect to $h_0$,
because $h_0$ and $h_i$ are mutually bounded.

\begin{lem}
\label{lem;17.11.30.10}
$\del_{E,h_0}\pi_i$
and 
$\delbar_E\pi_{i,h_0}^{\dagger}$
are $L^2$.
\end{lem}
\pf
Because the holomorphic decomposition 
$E=\bigoplus E_i$
is orthogonal with respect to $h_1$,
we have
$\del_{E,h_1}\pi_i=0$.
We have 
$\del_{E,h_1}=\del_{E,h_0}+b_1^{-1}\del_{E,h_0}b_1$.
Because
$\del_{E,h_0}\pi_{i,h_0}
=-[b_1^{-1}\del_{E,h_0}b_1,\pi_{i,h_0}]$,
we obtain that 
$\del_{E,h_0}\pi_{i,h_0}$ is $L^2$.
We also obtain that 
$\delbar_E\pi_{i,h_0}^{\dagger}$
are $L^2$.
\hfill\qed

\vspace{.1in}

We consider the Hermitian metric 
$h_3$ obtained as the direct sum of $h_{0|E_i}$.
\begin{lem}
\label{lem;17.11.30.11}
$h_3$ and $h_0$ are mutually bounded.
\end{lem}
\pf
Because $h_0$ and $h_1$ are mutually bounded,
$h_{0|E_i}$ and $h_{1|E_i}$ are mutually bounded.
Because $h_1=\bigoplus h_{1|E_i}$,
we obtain that
$h_1$ and $h_3$ are mutually bounded.
Then, we obtain the claim of the lemma.
\hfill\qed

\vspace{.1in}

Let $b_3$ be determined by $h_3=h_0b_3$.
We have
\[
 b_3=\sum_{j=1}^m \pi_{j,h_0}^{\dagger}\circ\pi_j.
\]

\begin{lem}
$b_3^{-1}\delbar_Eb_3$
and $b_{3}^{-1}\del_{E,h_0}b_3$
are $L^2$ with respect to $h_0$ and $g_X$.
\end{lem}
\pf
By Lemma \ref{lem;17.11.30.10},
$\delbar_Eb_3$ and $\del_{E,h_0}b_3$ are  $L^2$
with respect to $h_0$ and $g_X$.
By Lemma \ref{lem;17.11.30.11},
$b_3$ and $b_3^{-1}$ is also bounded with respect to $h_0$.
Then, we obtain the claim of the lemma.
\hfill\qed

\begin{lem}
$\Lambda\delbar\Tr(b_3^{-1}\del_{E,h_0}b_3)$
is $L^1$.
\end{lem}
\pf
We have 
$\delbar_E (b_3^{-1}\del_{E,h_0}b_3)
=-b_3^{-1}(\delbar_E b_3)\cdot b_3^{-1}\del_{E,h_0}b_3
+b_3^{-1}\delbar_E\del_{E,h_0}b_3$.
We also have the following:
\[
 \delbar_E\del_{E,h_0}b_3
=\sum \delbar_E\pi^{\dagger}_{j,h_0}\circ\del_{E,h_0}\pi_j
+\sum \pi^{\dagger}_{j,h_0}\circ [F(h_0),\pi_j].
\]
By the assumption
$|\Lambda F(h_0)|_{h_0}\leq B\varphi_X$,
$|\Lambda F(h_0)|_{h_0}$ is $L^1$.
Then, we obtain the claim of the lemma.
\hfill\qed

\begin{lem}
We have
$\int_X\Tr\bigl(
 \delbar_E(b_3^{-1}\del_{E,h_0}b_3)
 \bigr)\omega_X^{\dim X-1}=0$.
\end{lem}
\pf
We take a $C^{\infty}$-function
$\mu:\real\lrarr\real_{\geq 0}$
such that
$\mu(t)=0$ if $t\geq 2$
and $\mu(t)=1$ if $t\leq 1$.
Let $\varrho_1$ be an exhaustion function as in 
Proposition \ref{prop;17.11.30.3}.
We set
$\chi_N:=\mu\bigl(N^{-1}\varrho_1\bigr)$
on $X$.
Because
$\Tr\bigl(\delbar_E(b_3^{-1}\del_{E,h_0}b_3)\bigr)\omega_X^{\dim X-1}$
is $L^1$,
it is enough to prove
\begin{equation}
 \label{eq;17.11.30.30}
 \lim_{N\to\infty}
 \int_X\chi_N\cdot
 \Tr\bigl(\delbar_E(b_3^{-1}\del_{E,h_0}b_3)\bigr)
 \omega_X^{\dim X-1}=0.
\end{equation}
We have the following:
\begin{multline}
 \int_X\chi_N\cdot
 \Tr\bigl(\delbar_E(b_3^{-1}\del_{E,h_0}b_3)\bigr)
 \omega_X^{\dim X-1}
=-\int_X\delbar \chi_N\cdot
 \Tr\bigl(b_3^{-1}\del_{E,h_0}b_3\bigr)
 \omega_X^{\dim X-1}
 \\
=-\int_X
 \mu'\bigl(N^{-1}\varrho_1\bigr)
 N^{-1}\delbar\varrho_1\cdot
  \Tr\bigl(b_3^{-1}\del_{E,h_0}b_3\bigr)
 \omega_X^{\dim X-1}.
\end{multline}
Note that if
$\mu'\bigl(N^{-1}\varrho_1\bigr)\neq 0$,
we have
$N\leq \varrho_1\leq 2N$.
Hence, we have
$\bigl|N^{-1}\delbar\varrho_1\bigr|_{g_X}
\leq
 2\bigl|\varrho_1^{-1}\delbar\varrho_1\bigr|_{g_X}
=2\bigl|\delbar\log\varrho_1\bigr|_{g_X}$.
Then, we obtain 
(\ref{eq;17.11.30.30})
by the theorem of Lebesgue.
\hfill\qed

\vspace{.1in}
Because
$\Tr F(h_3)=\Tr F(h_0)+\delbar\Tr(b_3^{-1}\del_{E,h_0}b_3)$,
we obtain
\[
\int_X\Lambda\Tr F(h_0)
=\int_X\Lambda\Tr F(h_3)
=\sum_{i=1}^m
 \int_X\Lambda\Tr F(h_{0|E_i}).
\]
We also have $\rank E=\sum_{i=1}^m \rank E_i$.
Then, it is standard that
there exists $i_0$ such that
$\deg(E,h_0)/\rank E\leq \deg(E_{i_0},h_0)/\rank E_{i_0}$.
By the analytic stability of $(E,\delbar_E,h_0)$,
we obtain that $m=1$.
It implies $h_1=h_2$.
\hfill\qed

\section{Examples}
\label{section;17.11.14.1}

\subsection{Preliminary}

\subsubsection{Ahlfors type lemma}

Take $R>0$ and $C_0>0$.
We use the standard coordinate system
$\vecx=(x_1,\ldots,x_n)\in\real^n$.
We set $r(\vecx):=\sqrt{\sum_{i=1}^n x_i^2}$.
We set $U(R):=\bigl\{\vecx\in\real^n\,\big|\,r(\vecx)\geq R\bigr\}$.
We set
$\Delta:=-\sum\del_{x_i}^2$.

Let $C_0$ be a positive constant.
Let $g:U(R)\lrarr \real_{\geq 0}$ be 
a function such that
\[
 \Delta g\leq -C_0g.
\]
We also assume that $g=O(r(\vecx)^N)$ for some $N>0$.
\begin{lem}
\label{lem;17.10.7.300}
There exists $\epsilon_1>0$, depending only on $C_0$,
such that
$g(\vecx)=O\bigl(
 \exp(-\epsilon_1r(\vecx))
 \bigr)$.
\end{lem}
\pf
For any $a\in\real$,
we have
\[
 -(\del_r^2+(n-1)r^{-1}\del_r)
 \exp(ar)
=\bigl(
 -a^2-(n-1)r^{-1}a
 \bigr)
 \exp(ar).
\]
Hence, 
there exist $\epsilon_1>0$ and $R_1\geq R$
such that the following holds
on $\bigl\{\vecx\in\real^n\,\big|\,r(\vecx)\geq R_1\bigr\}$:
\[
\Delta
 \exp\bigl(-\epsilon_1 r(\vecx)\bigr)
\geq
 -C_0\exp\bigl(-\epsilon_1 r(\vecx)\bigr),
\quad\quad
\Delta \exp\bigl(\epsilon_1r(\vecx)\bigr)
\geq
 -C_0\exp\bigl(\epsilon_1r(\vecx)\bigr).
\]
We take $C_1>0$ such that
$g(\vecx)<C_1\exp\bigl(-\epsilon_1r(\vecx)\bigr)$
on $\{r(\vecx)=R_1\}$.
For any $\delta>0$,
we set
\[
 F_{\delta}(\vecx):=
 C_1\exp\bigl(-\epsilon_1r(\vecx)\bigr)
+\delta\exp\bigl(\epsilon_1r(\vecx)\bigr).
\]
We have $g(\vecx)<F_{\delta}(\vecx)$ on $\{r(\vecx)=R_1\}$.
We also have
$\Delta(g-F_{\delta})
\leq
 -C_0(g-F_{\delta})$.
We set
\[
 Z(\delta):=\Bigl\{
 \vecx\in\real^n\,\Big|\,
 r(\vecx)\geq R_1,\,\,
 g(\vecx)>F_{\delta}(\vecx)
 \Bigr\}.
\]
Because $g(\vecx)=O(r(\vecx)^N)$
and $g(\vecx)<F_{\delta}(\vecx)$ on $\{r(\vecx)=R_1\}$,
$Z(\delta)$ is relatively compact in $\{r(\vecx)>R_1\}$.
On $Z(\delta)$,
the following holds:
\[
 \Delta(g-F_{\delta})<0.
\]
On $\del Z(\delta)$,
we have $g-F_{\delta}=0$.
Hence, if $Z(\delta)\neq\emptyset$,
we obtain
$g-F_{\delta}\leq 0$ on $Z(\delta)$,
which contradicts with the construction of $Z(\delta)$.
Hence, we have $Z(\delta)=\emptyset$.
Namely, we have
$g\leq F_{\delta}$ on $\{r(\vecx)>R_1\}$
for any $\delta>0$.
We obtain that
$g\leq C_1\exp(-\epsilon_1r(\vecx))$.
\hfill\qed

\subsubsection{An estimate on $\real$}

Let $\varphi_{\real}$ be a positive $C^{\infty}$-function on $\real$
such that 
$\varphi_{\real}(t)=e^{-\delta |t|}$  $(|t|\geq 1)$.

\begin{lem}
\label{lem;17.10.16.10}
There exist positive constants $C_i$ $(i=0,1)$
such that the following holds.
\begin{itemize}
\item
Let $g$ be any bounded function 
$\real\lrarr\closedopen{0}{\infty}$
such that $-\del_t^2g\leq B\varphi_{\real}$ for some constant
$B>0$.
Then, $\sup g\leq C_0B+C_1\int\varphi_{\real} g$ holds.
\end{itemize}
\end{lem}
\pf
On $t\geq 1$,
we set
$F:=g+\delta^{-2}Be^{-\delta t}$.
Then, 
we have
$\del_t^2F\geq 0$ on $t\geq 1$.
Because $F$ is bounded,
we obtain that $\del_t F\leq 0$ on $t\geq 1$,
and hence $F(t)\leq F(1)$ $(t\geq 1)$.
It implies that
$g(t)\leq g(1)+\delta^{-2}Be^{-\delta}$ for $t\geq 1$.
Similarly,
we obtain that
$g(t)\leq g(-1)+\delta^{-2}Be^{-\delta}$ for $t\leq -1$.

We set $\psi(t):=\int_0^tds\int_0^{\tau}\varphi(\tau)\,d\tau$.
There exists a positive constant
$C_{10}$ such that
$|\psi(t)|\leq C_{10}$ for any $|t|\leq 3$.
We set $F:=g-B\psi(t)+BC_{10}$.
Then, $-\del_t^2F\leq 0$ and $F\geq 0$ hold.
By the convexity,
$F(t)\leq 2^{-1}(F(t+\tau)+F(t-\tau))$ hold
for any $-1\leq t\leq 1$ and $0\leq \tau\leq 1$.
Hence, for $-1\leq t\leq 1$, we obtain
\[
 F(t)\leq \int_0^12^{-1}(F(t+\tau)+F(t-\tau))d\tau
=\frac{1}{2}\int_{-1}^1F(t+\tau)\,d\tau
\leq \frac{1}{2}\int_{-2}^2F(s)\,ds.
\]
It implies the following for $-1\leq t\leq 1$:
\[
 g(t)\leq \frac{1}{2}\int_{-2}^2
 \bigl(
 g(s)-B\psi(t)+BC_{10}
 \bigr)\,ds
+B\psi(t)-BC_{10}.
\]
Hence, there exist $C_i>0$ $(i=3,4)$
such that 
$g(t)\leq C_3B+C_4\int_{-2}^2g\varphi$
for $-1\leq t\leq 1$.
Then, we obtain 
the claim of the lemma.
\hfill\qed

\subsubsection{Inequality for distributions}

Let $U$ be a neighbourhood of $(0,0,0)$ in $\real^3$.
Set $\Utilde:=S^1\times U$
and $W:=S^1\times\{(0,0,0)\}$.
We set $\Utilde^{\ast}:=\Utilde\setminus W$.
We regard $S^1=\real/\seisuu$.
Let $t$ be the standard coordinate of $\real$,
which induces local coordinates on $S^1$.
Let $(x,y,z)$ be the standard coordinate system of $\real^3$.
We set
$\Delta:=-(\del_t^2+\del_x^2+\del_y^2+\del_z^2)$.
Let $\psi$ be a bounded positive function on $\Utilde$.
The following lemma is implicitly implied in the proof of Proposition 2.2 of 
\cite{Simpson88}.
\begin{lem}
\label{lem;17.10.16.21}
Let $g:\Utilde^{\ast}\lrarr\real_{\geq 0}$ 
be a bounded function such that 
$\Delta g\leq \psi$ on $\Utilde^{\ast}$.
Then, 
$\Delta g\leq \psi$ holds
as distributions on $\Utilde$.
\end{lem}
\pf
In the proof of \cite[Proposition 2.2]{Simpson88},
it is studied in the case where $W$ is a complex submanifold.
The argument can work even in the case of real submanifolds
whose codimension is larger than $2$.
\hfill\qed

\subsection{Instantons}

\subsubsection{Doubly periodic instantons}

Let $\Gamma$ be a lattice in $\cnum$.
Let us consider the action of $\Gamma$
on $\cnum_z\times\cnum_w$
given by $\chi(z,w)=(z+\chi,w)$.
Set $X:=(\cnum_z\times\cnum_w)/\Gamma$
with the K\"ahler metric
$g_X:=dz\,d\zbar+dw\,d\wbar$.
Let $\dvol_X$ denote the associated volume form.
We set
$\dvol_{z}:=\frac{\sqrt{-1}}{2}dz\,d\zbar$
and 
$\dvol_{w}:=\frac{\sqrt{-1}}{2}dw\,d\wbar$.
Clearly, $\dvol_X=\dvol_z\,\dvol_w$ holds
on $X$.
We set
$\varphi_X:=(1+|w|^{2})^{-\delta-1}$
for a $\delta>0$.

\begin{prop}
\label{prop;17.10.15.30}
$(X,g_X,\varphi_X)$ satisfies 
Assumption {\rm\ref{assumption;17.11.29.40}}.
\end{prop}
\pf
Set $\varphi_{\cnum}:=(1+|w|^2)^{-1-\delta}$
on $\cnum$.
Let us regard $\cnum$ with the metric 
$\varphi_{\cnum}\,dw\,d\wbar$
as a K\"ahler manifold.
The Laplacian is given by $-\varphi_{\cnum}^{-1}\del_{w}\del_{\wbar}$.
Recall the following,
which is a special case of \cite[Proposition 2.4]{Simpson88}.
\begin{lem}
\label{lem;17.10.16.1}
There exists an increasing function
$a_1:\closedopen{0}{\infty}\lrarr \closedopen{0}{\infty}$
with $a_1(0)=0$ and $a_1(x)=x$ $(x\geq 1)$,
such that if $g$ is a bounded function $\cnum\lrarr \closedopen{0}{\infty}$
satisfying $-\varphi_{\cnum}^{-1}\del_w\del_{\wbar}g\leq B$,
then 
\[
 \sup g\leq C_1(B)a_1\Bigl(\int_{\cnum}g\varphi_{\cnum}\dvol_w\Bigr)
\]
holds,
where $C_1(B)$ is a positive constant depending on $B$.
Moreover, if $-\del_{\wbar}\del_wf\leq 0$
then $f$ is constant.
\hfill\qed
\end{lem}

Let $f$ be an $\real_{\geq 0}$-valued bounded
function on $X$ such that
$-(\del_z\del_{\zbar}+\del_w\del_{\wbar})f
 \leq B\varphi_X$ for $B>0$.
Set $T:=\cnum_z/\Gamma$.
Let $\vol(T)$ denote the volume of $T$
with respect to $\dvol_z$.
For any $w\in\cnum$,
we set 
$F(w):=\int_{T\times\{w\}}f\,\dvol_z$.
We set $B':=B\cdot\vol(T)$.
Then, we obtain
$-\del_w\del_{\wbar}F(w)\leq B'\varphi_{\cnum}$.
By Lemma \ref{lem;17.10.16.1},
we obtain 
\[
 \sup_{w\in\cnum} F(w)
\leq
 C_1(B')
a_1\Bigl(
 \int_{\cnum} F\varphi_{\cnum}\,\dvol_w
 \Bigr)=
C_1(B')
 a_1\Bigl(
 \int_{T\times\cnum}f\varphi_X\dvol_X
 \Bigr).
\]
For any $w_0\in\cnum_w$ and $r_0>0$,
we set $B_{r_0}(w_0):=
 \bigl\{w\in\cnum\,\big|\,|w-w_0|<r_0\bigr\}$.
Let $\vol(B_{r_0}(w_0))$ denote 
the volume of $B_{r_0}(w_0)$
with respect to $\dvol_w$.
We have 
$\Delta_X f\leq B\varphi_X$ for $B>0$
on $T\times B_{2}(w_0)$.
We also have
\[
 \int_{T\times B_2(w_0)} f\dvol_X
\leq
 \int_{B_2(w_0)}F\dvol_w
\leq
 C_1(B')
 {\rm vol\,}(B_2(w_0))
 a_1\Bigl(
 \int_{T\times\cnum}f\varphi_X\dvol_X
 \Bigr).
\]
Then, we there exist positive constants
$C_i$ $(i=2,3)$
such that
\[
 \sup_{w\in B_{1}(w_0)}f(w)
\leq
 C_2B+
 C_3C_1(B\vol(T))
 a_1\Bigl(
 \int_{T\times \cnum}f\varphi_X\dvol_X
 \Bigr).
\]
(See the proof of \cite[Proposition 2.1]{Simpson88}.)

Suppose that $f$ satisfies the stronger condition
$-(\del_z\del_{\zbar}+\del_w\del_{\wbar})f\leq 0$
on $X$.
Let us prove that $f$ is constant.
We set $F(w):=\int_{T\times\{w\}}f$ as above.
Then, 
$F$ is a bounded function on $\cnum_w$
satisfying $-\del_w\del_{\wbar}F\leq 0$.
Hence, by Lemma \ref{lem;17.10.16.1},
$F$ is constant.

We have the following:
\[
 -(\del_z\del_{\zbar}+\del_w\del_{\wbar})|f|^2
=
 -2(\del_z\del_{\zbar}+\del_w\del_{\wbar})f\cdot f
-|\del_zf|^2
-|\del_{\zbar}f|^2
-|\del_wf|^2
-|\del_{\wbar}f|^2
\leq
-|\del_{\zbar}f|^2.
\]
By using the Fourier expansion in the $T$-direction,
we have the decomposition $f=f_0+f_1$,
where $f_0$ is constant in the $T$-direction,
and $\int_{T\times\{w\}} f_1=0$.
We have 
$\int_{T\times\{w\}}|f|^2=
 \int_{T\times\{w\}}|f_0|^2
+\int_{T\times\{w\}}|f_1|^2$.
Because $F(w)=\int_{T\times\{w\}}f=\int_{T\times\{w\}}f_0$
is constant,
$f_0$ and $\int_{T\times\{w\}}|f_0|^2$ are constant.
We have $C>0$ such that
\[
 \int_{T\times\{w\}}
 \bigl|\del_{\zbar}f\bigr|^2
=\int_{T\times\{w\}}
 \bigl|\del_{\zbar}f_1\bigr|^2
\geq
 C\int_{T\times\{w\}}|f_1|^2.
\]
Hence, we obtain the following inequality
on $\cnum_w$:
\[
 -\del_{w}\del_{\wbar}
 \int_{T\times\{w\}}
 |f_1|^2
\leq
 -C\int_{T\times\{w\}}
 |f_1|^2.
\]
Because $\int_{T\times\{w\}}|f_1|^2$ is bounded,
it follows from Lemma \ref{lem;17.10.7.300}
that $\int_{T\times\{w\}}|f_1|^2
=O\bigl(\exp(-\epsilon|w|)\bigr)$
for some $\epsilon>0$ as $|w|\to\infty$.
Because $\int_{T\times\{w\}}|f_1|^2$
is subharmonic,
by using the maximum principle,
we obtain that 
$\int_{T\times\{w\}}|f_1|^2\leq 0$,
and hence $f_1=0$.
It implies that $f$ is constant.
\hfill\qed

\vspace{.1in}
Take $\lambda\in\cnum$.
We have the complex structure 
$(\xi,\eta)=(z+\lambda\wbar,w-\lambda\zbar)$
on $\real^4=\cnum_z\times\cnum_w$,
which induces a complex structure on $X$.
The complex manifold is denoted by $X^{\lambda}$.
Let $(E,\delbar_E)$ be a holomorphic vector bundle
on $X^{\lambda}$ with a Hermitian metric $h_0$
such that 
(a) $|\Lambda F(h_0)|_{h_0}\leq B\varphi_X$ for $B>0$,
(b) $\Tr F(h_0)=0$,
(c) $|F(h_0)|_{h_0,g_X}$ is $L^2$.

\begin{cor}
\label{cor;17.11.29.10}
If $(E,\delbar_E,h_0)$ is analytically stable,
then there exists a Hermitian metric $h$ of $E$
such that
(i) $h$ and $h_0$ are mutually bounded,
(ii) $\det(h)=\det(h_0)$,
(iii) $\Lambda F(h)=0$,
(iv) $\delbar (h h_0^{-1})$ and $|F(h)|_{h,g_X}$ are $L^2$.
Moreover, the equality 
{\rm(\ref{eq;18.1.21.1})} holds.
If $h'$ is a Hermitian metric of $E$
satisfying the conditions (i), (ii), (iii),
then $h'=h$.
\end{cor}
\pf
We take a positive $C^{\infty}$-function $\varrho$
on $X$ such that 
$\varrho(z,w)=\log|w|^2$ $(|w|>2)$.
Clearly $\varrho$ is an exhaustion function on $X$.
Note that
$w=(1+|\lambda|^2)^{-1}(\eta+\lambda\xibar)$.
On $\{|w|>1\}$,
we have the following equality
with respect to 
the complex structure of $X^{\lambda}$:
\[
 \delbar\varrho=
 \frac{1}{1+|\lambda|^2}
 \Bigl(
 \frac{d\etabar}{\wbar}+\lambda\frac{d\xibar}{w}
 \Bigr)
=O(|w|^{-1}),
\quad
 \delbar\del\varrho
=\frac{1}{(1+|\lambda|^2)^2}
 \Bigl(
-\lambda\frac{d\xibar\,d\eta}{w^2}
+\lambdabar\frac{d\etabar\,d\xi}{\wbar^2}
 \Bigr)
=O(|w|^{-2}).
\]
Hence, 
$\delbar\varrho$ is bounded,
and $\delbar\del\varrho$ is $L^2$ and bounded
on $X^{\lambda}$.
Moreover,
$\delbar\log\varrho=O\bigl((\log|w|)^{-1}|w|^{-1}\bigr)$
is $L^2$.
Then, the claim follows from
Theorem \ref{thm;17.8.8.11},
Corollary \ref{cor;17.11.29.20},
and 
Proposition \ref{prop;17.11.30.3}.
\hfill\qed

\begin{rem}
\label{rem;17.12.1.1}
In {\rm\cite[Theorem 0.12, Proposition 5.12]{Biquard-Jardim}},
in the rank $2$ case,
Biquard and Jardim studied
the correspondence between
instantons on $X^0$ with quadratic curvature decay
and parabolic bundles on the natural compactification of $X^0$.
(See {\rm\cite[Theorem 0.12]{Biquard-Jardim}} for more details)
In the proof, they seem to mention
a generalization of {\rm\cite[Theorem 1]{Simpson88}}
in the setting where the volume of the base K\"ahler manifold
is infinite, on the basis of {\rm\cite{b}}.
Indeed, Biquard kindly replied to the author that
the volume finiteness condition is not essential in his argument
in {\rm\cite{b}}.
It is not clear to the author how their generalization
is related to Corollary {\rm\ref{cor;17.11.29.10}}.

The author studied the correspondence between 
$L^2$-instantons on $X^0$
and polystable filtered bundles with degree $0$
on the natural compactification
of $X^0$ in {\rm\cite{Mochizuki-doubly-periodic}},
as a generalization of
{\rm\cite[Theorem 0.12]{Biquard-Jardim}}.
For the construction of $L^2$-instantons
from stable filtered bundles with degree $0$,
we used the Nahm transforms
between $L^2$-instantons on $X^0$
and wild harmonic bundles
on the dual torus 
of $\cnum/\Gamma$.
Corollary {\rm\ref{cor;17.11.29.10}} should allow us 
to construct $L^2$-instantons
from stable filtered bundles with degree $0$ more directly.

We plan to study the correspondence
between $L^2$-instantons on $X^{\lambda}$
and filtered bundles
on the natural compactification of $X^{\lambda}$
for general $\lambda$
by using Corollary {\rm\ref{cor;17.11.29.10}}.
\hfill\qed
\end{rem}

\begin{rem}
Let $T^{\lor}$ denote the dual torus of
$\cnum/\Gamma$.
In Remark {\rm\ref{rem;17.12.1.1}},
$X^{\lambda}$ is naturally identified
with the moduli space of
holomorphic line bundles of degree $0$ equipped with
a $\lambda$-connection on $T^{\lor}$,
and hence
$X^{\lambda}$ is naturally an affine line bundle over 
$T^{\lor}$.
The natural compactification of $X^{\lambda}$
is obtained as the projectivization of
a holomorphic bundle
$\gbige^{\lambda}$ on $T^{\lor}$.
There exists an exact sequence
$0\lrarr \nbigo_{T^{\lor}}
\lrarr
 \gbige^0
\lrarr
 \nbigo_{T^{\lor}}
\lrarr 0$.
If $\lambda=0$,
the exact sequence has a splitting.
If $\lambda\neq 0$,
the exact sequence has no splitting.
\hfill\qed
\end{rem}

\subsubsection{Spatially periodic instantons}
\label{subsection;17.10.16.20}

Let $\Gamma$ be a lattice in $\real^3$.
We consider 
$X:=(\real^3/\Gamma)\times\real$
with the Euclidean metric $g_X$
for which $\real^3/\Gamma$ and $\real$ are orthogonal.
We take a complex structure
$\real^4\simeq\cnum^2$ for which
the multiplication of $\sqrt{-1}$ is an isometry.
It induces a complex structure on $X$.
Let $t$ be the standard coordinate on $\real$.
Let $\varphi_{\real}:\real\lrarr\openopen{0}{\infty}$ be 
a $C^{\infty}$-function on $\real$
such that $\varphi_{\real}(t)=e^{-|t|\delta}$ $(|t|>1)$
for some $\delta>0$.
Let $\varphi_X$ be the $C^{\infty}$-function
$X\lrarr\openopen{0}{\infty}$
obtained as the pull back of $\varphi_{\real}$
by the projection
$X=(\real^3/\Gamma)\times\real\lrarr\real$.

\begin{prop}
\label{prop;17.10.16.22}
$(X,g_X,\varphi_X)$ satisfies 
Assumption {\rm\ref{assumption;17.11.29.40}}.
\end{prop}
\pf
Let $f$ be a bounded function on $X$
such that $\Delta(f)\leq B\varphi$ for $B>0$.
We consider 
the function $F(t):=\int_{T^3\times t}f$,
which satisfies
$-\del_t^2F(t)\leq B\varphi_{\real}$.
By Lemma \ref{lem;17.10.16.10},
we obtain that
$\sup F(t)\leq
 C_1B+
 C_2\int_{\real} F\varphi_{\real}$.
Then, 
we obtain the estimate
$\sup f(t)\leq 
 C_3B+C_4\int_{\real}F\varphi_{\real}$
as in the case of Proposition \ref{prop;17.10.15.30}.

Suppose that $f$ satisfies the stronger condition
$\Delta(f)\leq 0$ on $X$.
We set $F:=\int_{T^3\times\{t\}}f$ as above.
We obtain $-\del_t^2F\leq 0$.
Because $F:\real\lrarr \closedopen{0}{\infty}$
is bounded and convex from below,
we obtain that $F$ is constant.
Then, we obtain that $f$ is also constant
by using the argument in the proof of 
Proposition \ref{prop;17.10.15.30}.
\hfill\qed

\vspace{.1in}

Let $(E,\delbar_E)$ be a holomorphic vector bundle on $X$
with a Hermitian metric $h_0$
such that 
(a) $|\Lambda F(h_0)|_{h_0}\leq B\varphi_X$ for $B>0$,
(b) $\Tr F(h_0)=0$,
(c) $F(h_0)$ is $L^2$.

\begin{cor}
\label{cor;17.11.29.32}
If $(E,\delbar_E,h_0)$ is analytically stable,
then there exists a Hermitian metric $h$ of $E$
such that
(i) $h$ and $h_0$ are mutually bounded,
(ii) $\det(h)=\det(h_0)$,
(iii) $\Lambda F(h)=0$,
(iv) $\delbar_E(hh_0^{-1})$ and  $F(h)$ are $L^2$.
Moreover, the equality 
{\rm(\ref{eq;18.1.21.1})} holds.
If $h'$ is a Hermitian metric of $E$
satisfying the conditions (i), (ii) and (iii),
then $h=h'$.
\end{cor}
\pf
Let $\varrho_0$ be a positive $C^{\infty}$-function
on $\real$ such that 
$\varrho_0(t)=\log|t|$ if $|t|>1$.
Let $\varrho$ be the $C^{\infty}$-function on $X$
obtained as the pull back of $\varrho_0$
by the projection $X\lrarr \real$.
We can take a complex coordinate system
$(z,w)$ such that
$\Re(z)=t$ and 
$g_X=dz\,d\zbar+dw\,d\wbar$.
On $\{t>1\}$,
we have
$\del\varrho=2^{-1}t^{-1}dz$
and 
$\delbar\del\varrho=-4^{-1}t^{-2}d\zbar\,dz$.
We also have
$\del\log\varrho=2^{-1}(\log t)^{-1}t^{-1}dz$.
Hence, $\del\varrho$ and $\delbar\del\varrho$ are bounded,
and 
$\del\log\varrho$ and $\delbar\del\varrho$ are $L^2$.
Then, the claim follows from 
Theorem \ref{thm;17.8.8.11},
Corollary \ref{cor;17.11.29.20},
and Proposition \ref{prop;17.11.30.3}.
\hfill\qed

\begin{rem}
Spatially periodic instantons were first studied 
by Charbonneau in {\rm\cite{Charbonneau-CAG}},
and more recently by 
Charbonneau-Hurtubise {\rm\cite{Charbonneau-Hurtubise2}}
and Yoshino {\rm\cite{Yoshino-Master-Thesis}}.
It is natural to expect
to obtain a correspondence between
$L^2$-instantons on $X$
and filtered bundles
on the natural compactification of $X$
depending on the holomorphic structure.
Actually, it has been established 
in {\rm\cite{Charbonneau-Hurtubise2}},
at least under a kind of genericity assumption at infinity.
Corollary {\rm\ref{cor;17.11.29.32}}
could also give another way to study
such a correspondence.
\hfill\qed
\end{rem}

\subsection{Monopoles}

\subsubsection{Periodic monopoles}

For any $T>0$,
we set $S^1_T:=\real/T\seisuu$.
If $T=1$, we denote it by $S^1$.
Take a finite subset
$Z\subset S^1_T\times\cnum_w$.
Let $\pi:
 S^1\times S^1_T\times\cnum_w\lrarr S^1\times \cnum_w$
denote the projection
$\pi(t_1,t_2,w)=(t_2,w)$.
We regard $S^1\times S^1_T$
as the quotient space $\cnum_z/(\seisuu+\sqrt{-1}T\seisuu)$
by $(t_1,t_2)\longmapsto z=t_1+\sqrt{-1}t_2$.
We set
$X:=\bigl(
 S^1\times S^1_T\times\cnum_w
 \bigr)\setminus \pi^{-1}(Z)$
with the K\"ahler metric $g_X:=dz\,d\zbar+dw\,d\wbar$.
Set $\varphi_{X}(z,w):=(1+|w|^2)^{-\delta-1}$ for $\delta>0$.

\begin{prop}
\label{prop;17.10.16.23}
$(X,g_X,\varphi_X)$ satisfies 
Assumption {\rm\ref{assumption;17.11.29.40}}.
\end{prop}
\pf
Let $f$ be an $\real_{\geq 0}$-valued function 
such that
$\Delta_X f\leq B\varphi_X$ 
on $X$ for $B\geq 0$.
As remarked in Lemma \ref{lem;17.10.16.21},
the inequality holds on $S^1\times S_T^1\times\cnum$
as distributions.
Then, we obtain the claim 
from Proposition \ref{prop;17.10.15.30}.
\hfill\qed

\vspace{.1in}
Note that $(X,g_X)$ is a hyper-K\"ahler manifold.
We consider the complex manifold $X^{\lambda}$ corresponding
to the twistor parameter $\lambda\in\cnum$.
Indeed, we regard $X$ as the quotient of 
an open subset in $\cnum^2=\{(z,w)\}\simeq\real^4$.
Local holomorphic coordinates system for $X^{\lambda}$
are given by $(\xi,\eta)=(z+\lambda\wbar,w-\lambda \zbar)$.

\vspace{.1in}
We have the natural $S^1$-action on $S^1\times S^1_T\times\cnum_w$
given by $s_1(z,w)=(z+s_1,w)$.
Let $(E,\delbar_{E})$ be an $S^1$-equivariant
holomorphic vector bundle on $X^{\lambda}$.
Let $h_0$ be an $S^1$-invariant Hermitian metric of $E$
such that
(a) $|\Lambda F(h_0)|_{h_0}\leq B\varphi_X$ for $B>0$,
(b) $\Tr F(h_0)=0$.

\begin{cor}
\label{cor;17.11.29.30}
Suppose that 
$(E,\delbar_{E},h_0)$ is analytically stable
with respect to the $S^1$-action.
\begin{itemize}
\item
There exists a unique $S^1$-invariant Hermitian metric $h$ of $E$
such that 
(i) $h$ and $h_0$ are mutually bounded,
(ii) $\det(h)=\det(h_0)$,
(iii) $\Lambda F(h)=0$,
(iv) $\delbar_E(hh_0^{-1})$ is $L^2$.
If $h'$ is a Hermitian metric satisfying 
the conditions (i), (ii) and (iii),
then $h'=h$.
\item
If $(E,\delbar_E,h_0)$ satisfies
the additional condition 
(c) $F(h_0)^{\bot}\to 0$ and $\delbar_E\Lambda F(h_0)^{\bot}\to 0$
as $|w|\to\infty$,
then $h$ satisfies the condition
(v) $F(h)\to 0$ as $|w|\to\infty$.
\end{itemize}
\end{cor}
\pf
There exists $R>0$ such that
$Z$ is contained in
$S^1_T\times\{|w|<R\}$.
For each $P\in Z$,
let $r_P:(S^1_T\times\cnum_w)\lrarr \real_{\geq 0}$
be the distance function from $P$.
We take a positive $C^{\infty}$-function $\varrho_1$
on $X$ such that
$\varrho_1(z,w)=\log|w|$ if $|w|>R$,
and 
$\varrho_1=-\log (r_P\circ\pi)$ around $\pi^{-1}(P)$.
Clearly $\varrho_1$ is an exhaustion function 
on $X\setminus\pi^{-1}(Z)$.
Let us study
$\del\log\varrho_1$ on $X^{\lambda}\setminus\pi^{-1}(Z)$.
We have already observed that
the restriction to $\{|w|>R\}$ is $L^2$.
Around each $P\in Z$,
we have
\[
 \del\log\varrho_1
=\varrho_1^{-1}\cdot(r_P\circ\pi)^{-1}\cdot
 \del(r_P\circ\pi)
=O\bigl(
 \varrho_1^{-1}\cdot (r_P\circ\pi)^{-1}
 \bigr).
\]
Hence, $\del\log\varrho_1$ is $L^2$ around $P$.
Then, the claim follows from
Theorem \ref{thm;17.8.8.11},
Proposition \ref{prop;17.8.10.40}
and Proposition \ref{prop;17.11.30.3}.
\hfill\qed

\begin{rem}
Note that 
$S^1$-equivariant instantons on $X^{\lambda}$
are equivalent to monopoles on 
$(S^1\times\cnum_w)\setminus Z$.
In {\rm\cite{Mochizuki-KH-periodic}},
by using Corollary {\rm\ref{cor;17.11.29.30}},
we study the correspondence
of monopoles on $(S^1\times \cnum_w)\setminus Z$
of GCK type
and filtered mini-holomorphic bundles on the compactification
of $(S^1\times\cnum_w)\setminus Z$
depending on $\lambda$.
The latter objects are regarded as 
difference modules with parabolic structure.
\hfill\qed
\end{rem}

\begin{rem}
In {\rm\cite{Charbonneau-Hurtubise}},
Charbonneau and Hurtubise efficiently applied
the theorem of Simpson 
for the construction of Hermitian-Einstein monopoles
with Dirac type singularity
on the product of a circle and a Riemann surface.

\hfill\qed
\end{rem}

\subsubsection{Doubly periodic monopoles}

Let $T^2$ be the quotient space
$\real^2/\Gamma$ for a lattice $\Gamma\subset\real^2$.
Take a finite subset
$Z\subset T^2\times\real$.
Let $\pi:S^1\times T^2\times\real\lrarr T^2\times\real$
be the projection.
We regard $S^1\times T^2$ as the quotient of $\real^3$.
We set $X:=(S^1\times T^2\times\real)\setminus \pi^{-1}(Z)$
with the Euclidean metric $g_X$
as in \S\ref{subsection;17.10.16.20}.
We take 
an $\real$-isomorphism
$\real\times\real^2\times\real\simeq\cnum^2$
for which the multiplication of $\sqrt{-1}$ is an isometry
with respect to $g_X$.
It induces a complex structure on $X$.
Let $\varphi_X$ be a $C^{\infty}$-function on $X$
as in \S\ref{subsection;17.10.16.20}.

\begin{prop}
\label{prop;18.1.20.1}
$(X,g_X,\varphi_X)$ satisfies 
Assumption {\rm\ref{assumption;17.11.29.40}}.
\end{prop}
\pf
It follows from Proposition \ref{prop;17.10.16.22}
as in the case of Proposition \ref{prop;17.10.16.23}.
\hfill\qed

\vspace{.1in}
Let $(E,\delbar_E)$ be an $S^1$-equivariant
holomorphic vector bundle on $X^{\lambda}$
with an $S^1$-invariant Hermitian metric $h_0$.
Suppose that
(a) $|\Lambda F(h_0)|_{h_0}\leq B\varphi_X$
for some $B>0$,
(b) $\Tr F(h_0)=0$.
By an argument similar to the proof of Corollary \ref{cor;17.11.29.30},
but by using Proposition \ref{prop;18.1.20.30}
instead of Proposition \ref{prop;17.8.10.40},
we obtain the following.

\begin{cor}
\label{cor;17.11.29.31}
Suppose that $(E,\delbar_E,h_0)$ on $X^{\lambda}$
is analytically stable with respect to the $S^1$-action.
\begin{itemize}
\item
There exists a unique $S^1$-invariant Hermitian metric $h$
of $E$
such that 
(i) $h$ and $h_0$ are mutually bounded,
(ii) $\det(h)=\det(h_0)$,
(iii) $\Lambda F(h)=0$,
(iv) $\delbar_E(h h_0^{-1})$ is $L^2$.
If $h'$ is an $S^1$-invariant Hermitian metric
satisfying the conditions (i), (ii) and (iii),
then $h'=h$.
\item
If $h_0$ satisfies the additional condition
(c) $F(h_0)^{\bot}$
 and $\delbar_E\Lambda F(h_0)^{\bot}$
 are bounded on $\{|t|>R\}$ for some $R$,
then 
$h$ satisfies the condition
(v) $F(h)$ is bounded on $\{|t|>R\}$.
\hfill\qed
\end{itemize}
\end{cor}

\begin{rem}
$S^1$-equivariant instantons on $X$
are equivalent to monopoles on $(T^2\times\real)\setminus Z$.
We plan to study the correspondence 
between monopoles on $(T^2\times\real)\setminus Z$
and filtered mini-holomorphic bundles
on the compactifications
depending on complex structures,
for which Corollary {\rm\ref{cor;17.11.29.31}} would be useful.
For appropriate complex structures,
the latter objects are regarded as 
$q$-difference modules with parabolic structure.
\hfill\qed
\end{rem}

\subsubsection{Doubly periodic monopole corresponding
 to some concrete algebraic data}
\label{subsection;18.1.20.12}

To show that Proposition \ref{prop;18.1.20.1}
(and hence Theorem \ref{thm;17.8.8.11})
is already useful for the construction of doubly periodic monopoles,
let us explain the existence of doubly periodic monopoles
corresponding to some algebraic data.
There are many similar constructions.

For simplicity,
let us consider the case
where $T^2$ is isometric to
the product $\real/\seisuu\times\real/2\pi\seisuu$.
We may identify 
$S^1\times T^2\times\real$ with the Euclidean metric
and 
$\bigl(
 \cnum_z/(\seisuu+\sqrt{-1}\seisuu)
 \bigr)\times \cnum_w^{\ast}$
with the metric $dz\,d\zbar+|w|^{-2}dw\,d\wbar$.
We regard it as the quotient of
$\cnum\times\cnum^{\ast}$
by the $(\seisuu+\sqrt{-1}\seisuu)$-action
given as 
$\kappatilde_{\chi}(z,w)=(z+\chi,w)$ 
$(\chi\in\seisuu+\sqrt{-1}\seisuu)$.
We set $T_z:=\cnum_z/(\seisuu+\sqrt{-1}\seisuu)$
and $\Gamma:=\seisuu+\sqrt{-1}\seisuu$.
We define the $\real$-action
on $\cnum_z\times\cnum_w^{\ast}$
by $\rho_s(z,w)=(z+s,w)$.
The induced  $S^1$-action on 
$T_z\times\cnum_w^{\ast}$
is also denoted by $\rho$.

\vspace{.1in}
First, we give some easy and explicit examples of rank one.
Take $\alpha\in\cnum^{\ast}$.
Let $L$ be the line bundle on $\cnum_z\times\cnum_w^{\ast}$
with a frame $e$.
We define the $\Gamma$-action on $L$
by $\kappatilde_1^{\ast}(e)=e$
and $\kappatilde_{\sqrt{-1}}^{\ast}(e)=\alpha e$.
We also define the $\real$-action on $L$ by
$\rho_s^{\ast}(e)=e$.
We obtain the $S^1$-equivariant
holomorphic bundle $\nbigl$
on $S^1\times\cnum_w^{\ast}$
induced by $L$.
Let $h_{L,\alpha}$ be the metric of $L$
given by 
\[
 h_{L,\alpha}(e,e)=|\alpha|^{2\Image(z)}.
\]
Then, we have
$\kappatilde_{\chi}^{\ast}h_{L,\alpha}=h_{L,\alpha}$
for any $\chi\in\Gamma$
and $\rho_s^{\ast}h_{L,\alpha}=h_{L,\alpha}$
for any $s\in\real$.
The curvature of
$(L,\delbar_L,h_{L,\alpha})$
is 
$\delbar\del\log|\alpha|^{2\Image(z)}=0$.
Hence, $(L,\delbar_L,h_{L,\alpha})$ gives 
an $(S^1\times\Gamma)$-equivariant
unitary flat line bundle on $\cnum_z\times \cnum_w^{\ast}$.
It induces an $S^1$-equivariant instanton
$(\nbigl,\delbar_{\nbigl},h_{\nbigl,\alpha})$
on $T_z\times\cnum_w^{\ast}$,
i.e., a monopole on $S^1\times \cnum_w^{\ast}$.

\vspace{.1in}
Let $V$ be a holomorphic vector bundle
on $\cnum_z\times\cnum_w^{\ast}$ of rank $2$
equipped with a holomorphic frame $e_1,e_2$.
Take $a\in\cnum\setminus\{0,1,-1\}$.
We define the $\Gamma$-action on $V$
as follows:
\[
 \kappatilde_1^{\ast}(e_1,e_2)=(e_1,e_2),
\quad
 \kappatilde_{\sqrt{-1}}^{\ast}(e_1,e_2)=
 (e_1,e_2)
 \left(
 \begin{array}{cc}
 w & a \\
 a & w^{-1}
 \end{array}
 \right).
\]
We also define the $\real$-action on $V$
by 
$\rho_s^{\ast}(e_1,e_2)=(e_1,e_2)$.
We obtain an $S^1$-equivariant holomorphic bundle
$\nbigv$ on $T_z\times\cnum_w^{\ast}$.

Set $L:=\det(V)$,
which is equipped with the frame
$e:=e_1\wedge e_2$.
We have  
$\kappatilde_1^{\ast}(e)=e$,
$\kappatilde_{\sqrt{-1}}^{\ast}(e)=(1-a^2)e$
and $\rho_s^{\ast}(e)=e$.
Note that 
$(L,\delbar_L,h_{L,1-a^2})$ is 
an $S^1$-equivariant unitary flat line bundle, as explained above.
We have the induced $S^1$-equivariant instanton
$(\nbigl,\delbar_{\nbigl},h_{\nbigl,1-a^2})$
on $T_z\times\cnum_w^{\ast}$.

Let us prove that
we have a family of $S^1$-invariant Hermitian metrics
$h^{(c_0,c_{\infty})}_{\nbigv}$ $((c_0,c_{\infty})\in\real^2)$
of $\nbigv$
such that
(i) $\Lambda F(h^{(c_0,c_{\infty})}_{\nbigv})=0$,
(ii) $\det(h^{(c_0,c_{\infty})}_{\nbigv})=h_{\nbigl,1-a^2}$.
It induces a family of monopoles
on $S^1\times\cnum_w^{\ast}$.

\vspace{.1in}

As a preliminary,
we construct
$S^1$-invariant metrics $h^{(c_0,c_{\infty})}_{0,\nbigv}$
of $\nbigv$ which 
satisfies the ASD-equation outside a compact subset.
Note that
the roots of the polynomial $T^2-(w+w^{-1})T+1-a^2$ are
\[
 \frac{1}{2}
 \bigl(
 w+w^{-1}\pm\sqrt{(w+w^{-1})^2-4(1-a^2)}
 \bigr).
\]

Let $U_{\infty}$ be a small neighbourhood of
$\infty$ in $\proj^1_w$.
Set $U_{\infty}^{\ast}:=U_{\infty}\setminus\{\infty\}$.
On $U_{\infty}$,
let $\delta_{\infty}(w^{-1})$ be the branch of
$\sqrt{1+(-2+4a^2)w^{-2}+w^{-4}}$
such that $\delta_{\infty}(0)=1$.
We set
\[
 \beta_1(w^{-1})=
 \frac{w}{2}\bigl(1+w^{-2}+\delta_{\infty}(w^{-1})\bigr).
\]
It is a root of $T^2-(w+w^{-1})T+1-a^2$.
The other root is
\[
 \beta_2(w^{-1})=
 \frac{w}{2}\bigl(1+w^{-2}-\delta_{\infty}(w^{-1})\bigr)
=\beta_1(w^{-1})^{-1}(1-a^2).
\]
We define a holomorphic frame 
$v_1,v_2$ of $V$
on $\cnum_z\times U_{\infty}^{\ast}$
as follows:
\[
v_i:=ae_1+(\beta_i(w^{-1})-w)e_2
\quad(i=1,2).
\]
Then, we have
$\kappatilde_{n_1+n_2\sqrt{-1}}(v_i)
=\beta_i(w^{-1})^{n_2}v_i$
and $\rho_s^{\ast}(v_i)=v_i$ $(i=1,2)$.
We have
\[
 v_1\wedge v_2=a\bigl(\beta_2(w^{-1})-\beta_1(w^{-1})\bigr)e
=-aw\delta_{\infty}(w^{-1})e.
\]

Take any $c_{\infty}\in\real$.
We define the metric $h^{(c_{\infty})}_{0,U_{\infty}^{\ast}}$
of $V_{|\cnum_z\times U_{\infty}^{\ast}}$ as follows:
\[
 h^{(c_{\infty})}_{0,U_{\infty}^{\ast}}(v_1,v_1)
=|w|^{c_{\infty}}
|\beta_1(w^{-1})|^{2\Image(z)}
 \bigl|aw\delta_{\infty}\bigr|^2,
\quad
 h^{(c_{\infty})}_{0,U_{\infty}^{\ast}}(v_2,v_2)
=|w|^{-c_{\infty}}
 |\beta_2(w^{-1})|^{2\Image(z)},
\quad
 h^{(c_{\infty})}_{0,U_{\infty}^{\ast}}(v_1,v_2)=0.
\]
Then, the metric $h^{(c_{\infty})}_{0,U_{\infty}^{\ast}}$
is invariant under the $\Gamma$-action
and the $\real$-action.
The induced metric on $L_{|\cnum_z\times U_{\infty}^{\ast}}$
is equal to $h_{L,1-a^2|\cnum_z\times U_{\infty}^{\ast}}$.
We can easily see the following for $i=1,2$:
\[
 \del_{z}\del_{\zbar}
 \log h^{(c_{\infty})}_{0,U_{\infty}^{\ast}}(v_i,v_i)
=w\del_{w}\wbar\del_{\wbar}
 \log h^{(c_{\infty})}_{0,U_{\infty}^{\ast}}(v_i,v_i)
=0.
\]
We also have
\[
 \del_z\bigl(
 \wbar\del_{\wbar}
 \log h^{(c_{\infty})}_{0,U_{\infty}^{\ast}}(v_i,v_i)
 \bigr)
=\frac{1}{2\sqrt{-1}}
 \frac{\overline{w\del_w\beta_1(w^{-1})}}
 {\overline{\beta_1(w^{-1})}},
\quad
  \del_{\zbar}\bigl(
 w\del_{w}
 \log h^{(c_{\infty})}_{0,U_{\infty}^{\ast}}(v_i,v_i)
 \bigr)
=\frac{\sqrt{-1}}{2}
 \frac{w\del_{w}\beta_1(w^{-1})} {\beta_1(w^{-1})}.
\]
Hence, we have
$\Lambda F(h^{(c_{\infty})}_{0,U_{\infty}^{\ast}})=0$,
and $F(h^{(c_{\infty})}_{0,U_{\infty}^{\ast}})$
is bounded.

\vspace{.1in}
We make a similar construction around $w=0$.
Let $U_0$ be a small neighbourhood of $0$ in $\cnum_w$.
We set $U_0^{\ast}:=U_0\setminus\{0\}$.
On $U_0$,
let $\delta_0(w)$ be the branch of
$\sqrt{1+(-2+4a^2)w^2+w^4}$
such that $\delta_0(0)=1$.
We set
\[
 \gamma_1(w)=
 \frac{w^{-1}}{2}
 \bigl(
 1+w^2+\delta_0(w)
 \bigr).
\]
It is a root of $T^2-(w+w^{-1})T+1-a^2$.
The other root is
\[
 \gamma_2(w)=
 \frac{w^{-1}}{2}
 \bigl(
 1+w^2-\delta_0(w)
 \bigr)
=\gamma_1(w)^{-1}(1-a^2).
\]
We define a holomorphic frame 
$u_1,u_2$ of $V$ on $\cnum_z\times U_0^{\ast}$
as follows:
\[
 u_i=ae_1+(\gamma_i(w)-w)e_2
\quad(i=1,2).
\]
Then, we have
$\kappatilde_{n_1+\sqrt{-1}n_2}(u_i)
=\gamma_i(w)^{n_2}u_i$
and 
$\rho_s^{\ast}(u_i)=u_i$
for $i=1,2$.
We have
$u_1\wedge u_2
=-a w \delta_{0}e$.

Take any $c_0\in\real$.
We define the metric
$h^{(c_0)}_{0,U_0^{\ast}}$
of $V_{|\cnum_z\times U_0^{\ast}}$
as follows:
\[
 h^{(c_0)}_{0,U_0^{\ast}}(u_1,u_1)
=|w|^{-c_0}|\gamma_1(w)|^{2\Image(z)}|aw^{-1}\delta_0(w)|^2,
\quad
  h^{(c_0)}_{0,U_0^{\ast}}(u_2,u_2)
=|w|^{c_0}|\gamma_2(w)|^{2\Image(z)},
\quad
  h^{(c_0)}_{0,U_0^{\ast}}(u_1,u_2)=0.
\]
Then, the metric $h^{(c_0)}_{0,U_0^{\ast}}$
is invariant under the action of
$\real$ and $\Gamma$.
The induced metric on $L_{|\cnum_z\times U_0^{\ast}}$
is equal 
to $h_{L,1-a^2|\cnum_z\times U_0^{\ast}}$.
As in the case of $h^{(c_{\infty})}_{0,U_{\infty}^{\ast}}$,
we have
$\Lambda F(h^{(c_0)}_{0,U_0^{\ast}})=0$,
and 
$F(h^{(c_0)}_{0,U_0^{\ast}})$ is bounded.

We take a Hermitian metric $h^{(c_0,c_{\infty})}_{0,V}$
of $V$ such that 
(i) it is invariant under the actions of $\real$ and $\Gamma$,
(ii) the induced metric on $L$ is equal to
 $h_{L,1-a^2}$,
(iii) the restriction to $\cnum_z\times U_{\infty}^{\ast}$
 is equal to
 $h^{(c_{\infty})}_{0,U_{\infty}^{\ast}}$,
(iv)
 the restriction to $\cnum_z\times U_{0}^{\ast}$
 is equal to
 $h^{(c_{0})}_{0,U_{0}^{\ast}}$.

We obtain the $S^1$-invariant Hermitian metric
$h^{(c_0,c_{\infty})}_{0,\nbigv}$ of $\nbigv$
induced by $h^{(c_0,c_{\infty})}_{0,V}$.
By the construction,
we have
$\det(h^{(c_0,c_{\infty})}_{0,\nbigv})=h_{\nbigl,1-a^2}$.
We also have
$\Lambda F(h^{(c_0,c_{\infty})}_{0,\nbigv})=0$
on 
$T_z\times (U_{0}^{\ast}\sqcup U_{\infty}^{\ast})$,
and 
$F(h^{(c_0,c_{\infty})}_{0,\nbigv})$
is bounded on $S^1\times\cnum_w^{\ast}$.

\begin{prop}
\label{prop;18.1.20.40}
For each $(c_0,c_{\infty})\in\real^2$,
we have a unique $S^1$-invariant
Hermitian metric $h^{(c_0,c_{\infty})}_{\nbigv}$
of $\nbigv$
such that
(i) $h^{(c_0,c_{\infty})}_{\nbigv}$
 and 
$h^{(c_0,c_{\infty})}_{0,\nbigv}$
are mutually bounded,
(ii) $\det(h^{(c_0,c_{\infty})}_{\nbigv})
 =h_{\nbigl,1-a^2}$,
(iii) $\Lambda F(h^{(c_0,c_{\infty})}_{\nbigv})=0$,
(iv)  $\delbar_E\bigl(
 h_{\nbigv}^{(c_0,c_{\infty})}
\cdot
 (h_{0,\nbigv}^{(c_0,c_{\infty})})^{-1}
 \bigr)$ is $L^2$,
(v) $F(h^{(c_0,c_{\infty})}_{\nbigv})$ is bounded.
\end{prop}
\pf
According to Corollary \ref{cor;17.11.29.31},
it is enough to prove that
$(\nbigv,\delbar_{\nbigv},h^{(c_0,c_{\infty})}_{0,\nbigv})$
is analytically stable.
Let us check a stronger condition
 that
a non-trivial $S^1$-equivariant holomorphic subbundle
of $\nbigv$
is $0$ or $\nbigv$.

Let $V'$ be the product bundle
on $\real\times\cnum_w^{\ast}$
with a frame $e_1',e_2'$.
We define the action of $\seisuu$ on $V'$
by
\[
 \kappa_1^{\ast}(e_1',e_2')
=(e_1',e_2')
 \left(
 \begin{array}{cc}
 w & a \\ a & w^{-1}
 \end{array}
 \right).
\]
We obtain the bundle $\nbigv'$
on $S^1\times\cnum_w^{\ast}$.

We define the differential operators
$\del_{V',t}$ and $\del_{V',\wbar}$ on $V'$
as follows,
for any $f\in C^{\infty}(\real\times\cnum_w^{\ast})$:
\[
 \del_{V',t}(fe_i')=\del_{t}(f)e_i',
\quad
 \del_{V',\wbar}(fe_i')=\del_{\wbar}(f)e_i'.
\]
It induces differential operators
$\del_{\nbigv',t}$
and $\del_{\nbigv',\wbar}$
on $\nbigv'$.

The restriction
$V'_{|\{0\}\times\cnum_w^{\ast}}$
and 
$\nbigv'_{|\{0\}\times\cnum_w^{\ast}}$
are holomorphic vector bundles
by $\del_{\nbigv',\wbar}$.
We have the isomorphism
$\Phi:
 V'_{|\{0\}\times\cnum_w^{\ast}}
\simeq
V'_{|\{1\}\times\cnum_w^{\ast}}$
given by the parallel transport 
along the paths
$(t,w)$ $(0\leq t\leq 1)$ for any $w\in\cnum^{\ast}$
with respect to $\del_{V',t}$.
Because $\del_{V',t}$ and $\del_{V',\wbar}$
are commutative,
$\Phi$ is holomorphic.
It induces the monodromy automorphism
$\Phi$
on $\nbigv'_{|\{0\}\times\cnum_w^{\ast}}$.

\begin{lem}
\label{lem;18.1.20.10}
Suppose that $\nbigv'_1$ is a subbundle of $\nbigv'$
such that the following holds:
\begin{itemize}
\item
 If $s$ is a local section of $\nbigv'_1$,
 then
 $\del_{\nbigv',t}s$ and $\del_{\nbigv',\wbar}s$
 are also local sections of $\nbigv'_1$.
\end{itemize}
Then, $\nbigv'_1$ is $\nbigv'$ or $0$.
\end{lem}
\pf
By the condition,
$\nbigv'_{1|\{0\}\times\cnum_w^{\ast}}$
is a holomorphic subbundle of $\nbigv'_{|\{0\}\times\cnum^{\ast}}$
such that 
$\Phi(\nbigv'_{1|\{0\}\times\cnum_w^{\ast}})
 =\nbigv'_{1|\{0\}\times\cnum_w^{\ast}}$.
Because the eigenvalues of $\Phi$ 
are not single-valued
under the assumption $a\in\cnum\setminus\{0,1,-1\}$,
we can conclude that
$\nbigv'_{1|\{0\}\times\cnum_w^{\ast}}$
is 
$\nbigv_{1|\{0\}\times\cnum_w^{\ast}}$
or $0$.
Then, we obtain the claim of Lemma \ref{lem;18.1.20.10}.
\hfill\qed

\vspace{.1in}

Let $p:\cnum_z\times \cnum_w^{\ast}
 \lrarr \real\times\cnum_w^{\ast}$
be given by
$p(z,w)=(\Image(z),w)$.
We have the $C^{\infty}$-isomorphism
$V\simeq p^{-1}(V')$
by identifying
$e_i$ and $p^{-1}(e_i')$.
By the construction,
for any section $s$ of $V'$,
we have
\begin{equation}
\label{eq;18.1.20.11}
 \del_{V,\zbar}p^{-1}(s)
=\frac{\sqrt{-1}}{2}p^{-1}(\del_{V',t}s),
\quad
 \del_{V,\wbar}p^{-1}(s)
=p^{-1}(\del_{V',\wbar}s).
\end{equation}

Let $p_1:T_z\times\cnum_w^{\ast}\lrarr S^1\times\cnum_w^{\ast}$
be the projection induced by $p$.
The isomorphism
$V\simeq p^{-1}(V')$
induces
an isomorphism
$\nbigv\simeq p_1^{-1}(V')$.

Let $\nbigv_1\subset \nbigv$ be 
an $S^1$-invariant holomorphic subbundle.
It induces a $C^{\infty}$-subbundle 
$\nbigv_1'\subset \nbigv'$.
Moreover, for any section $s$ of $\nbigv_1'$,
$\del_{\nbigv',t}s$ and $\del_{\nbigv',\wbar}s$
are also sections of $\nbigv_1'$,
which we can observe by (\ref{eq;18.1.20.11}).
Hence, by Lemma \ref{lem;18.1.20.10},
we obtain that $\nbigv'_1$ is $\nbigv'$ or $0$.
It implies that $\nbigv_1$ is $\nbigv$ or $0$.
Thus, Proposition \ref{prop;18.1.20.40} is proved.
\hfill\qed

\noindent
{\em Address\\
Research Institute for Mathematical Sciences,
Kyoto University,
Kyoto 606-8502, Japan\\
takuro@kurims.kyoto-u.ac.jp
}


\begin{thebibliography}{99}

\bibitem{Bando}
S. Bando,
{\em Einstein-Hermitian metrics on noncompact K\"{a}hler manifolds},
in
{\em Einstein metrics and Yang-Mills connections (Sanda, 1990)}, 
Lecture Notes in Pure and Appl. Math., {\bf 145}, Dekker, New York,
(1993), 27--33.

\bibitem{b}
O. Biquard,
{\it Fibr\'es de Higgs et connexions int\'egrables:
le cas logarithmique (diviseur lisse)},
Ann. Sci. \'Ecole Norm. Sup. {\bf 30} (1997), 41--96.

\bibitem{biquard-boalch}
O. Biquard and P. Boalch,
{\em Wild non-abelian Hodge theory on curves},
Compos. Math. {\bf 140} (2004),
179--204.

\bibitem{Biquard-Jardim}
O. Biquard,
M. Jardim,
{\em Asymptotic behaviour and the moduli space of 
doubly-periodic instantons}. J. Eur. Math. Soc. {\bf 3}
(2001), 335--375.

\bibitem{Charbonneau-CAG}
B. Charbonneau,
{\em From spatially periodic instantons to singular monopoles.}
Comm. Anal. Geom. {\bf 14} (2006), 183--214.

\bibitem{Charbonneau-Hurtubise}
B. Charbonneau, 
J. Hurtubise,
{\em Singular Hermitian-Einstein monopoles 
on the product of a circle and a Riemann surface},
Int. Math. Res. Not. (2011), 175--216. 

\bibitem{Charbonneau-Hurtubise2}
B. Charbonneau,
J. Hurtubise,
{\em Spatially periodic instantons:
 Nahm transform and moduli},
Comm. Math. Phys. online version.
%https://doi.org/10.1007/s00220-018-3155-3
%arXiv:1712.05075

\bibitem{Donaldson-Narasimhan-Seshadri}
S. K. Donaldson,
{\em A new proof of a theorem of Narasimhan and Seshadri},
J. Differential Geom. {\bf 18} (1983), 269--277. 

\bibitem{Donaldson-GIT}
S. K. Donaldson,
{\em Instantons and Geometric Invariant Theory},
Comm. Math. Phys. {\bf 93} (1984), 453--460.

\bibitem{Donaldson-surface}
S. K. Donaldson,
{\em Anti self-dual Yang-Mills connections over 
complex algebraic surfaces and stable vector bundles},
Proc. London Math. Soc. (3) {\bf 50} (1985), 1--26.

\bibitem{Donaldson-infinite}
S. K. Donaldson, 
{\em Infinite determinants, stable bundles and curvature},
Duke Math. J. {\bf 54} (1987), 231--247.

\bibitem{Donaldson-boundary-value}
S. K. Donaldson,
{\em Boundary value problems for Yang-Mills fields},
J. Geom. Phys. {\bf 8} (1992), 89--122. 

\bibitem{Donaldson-Kronheimer}
S. K. Donaldson,
P. B. Kronheimer,
{\em The geometry of four-manifolds},
Oxford Science Publications. 
The Clarendon Press, 
Oxford University Press, New York, 1990.
x+440 pp.

\bibitem{Gilbarg-Trudinger}
D. Gilbarg, 
N. S. Trudinger,
{\em Elliptic partial differential equations of second order},
Reprint of the 1998 edition. Classics in Mathematics. 
Springer-Verlag, Berlin, 2001. xiv+517 pp.

\bibitem{GuoI}
G.-Y. Guo,
{\em Yang-Mills fields on cylindrical manifolds 
and holomorphic bundles. I},
Comm. Math. Phys. {\bf 179} (1996), 737--775.

\bibitem{GuoII}
G.-Y. Guo,
{\em Yang-Mills fields on cylindrical manifolds 
and holomorphic bundles. II},
Comm. Math. Phys. {\bf 179} (1996), 777--788.

\bibitem{Hamilton-book}
R. S. Hamilton,
{\em Harmonic maps of manifolds with boundary},
Lecture Notes in Mathematics, {\bf 471},
Springer-Verlag, Berlin-New York, 1975. i+168 pp. 

\bibitem{Hitchin-note-vanishing}
N. Hitchin,
{\em A note on vanishing theorems};
in 
{\em Geometry and analysis on manifolds}, 373--382,
Progr. Math., {\bf 308}, Birkh\"{a}user/Springer, Cham, 2015.

\bibitem{Jarvis-construction}
S. Jarvis, 
{\em Construction of Euclidean monopoles.}
Proc. London Math. Soc. (3) {\bf 77} (1998), 193--214.

\bibitem{Kobayashi-Nagoya}
S. Kobayashi, 
{\em First Chern class and holomorphic tensor fields},
Nagoya Math. J. {\bf 77} (1980), 5--11. 

\bibitem{Kobayashi-Academy}
S. Kobayashi,
{\em Curvature and stability of vector bundles},
Proc. Japan Acad. Ser. A Math. Sci. {\bf 58} (1982), 158--162.


\bibitem{Kobayashi-vector-bundle}
S. Kobayashi,
{\em Differential geometry of complex vector bundles}.
Publications of the Mathematical Society of Japan, {\bf 15}. 
Princeton University Press, Princeton, NJ; 
Princeton University Press, Princeton, NJ, 1987. xii+305 pp.

\bibitem{Lee-manifolds}
J. M. Lee,
{\em Introduction to smooth manifolds},
Graduate Texts in Mathematics, 218. 
Springer-Verlag, New York, 2003. xviii+628 pp

\bibitem{li2}
J. Li,
{\em Hermitian-Einstein metrics and Chern number inequalities
on parabolic stable bundles over K\"{a}hler manifolds,} %\"
Comm. Anal. Geom. {\bf 8} (2000), 445--475.

\bibitem{Li-Narasimhan}
J. Li, 
M. S. Narasimhan,
{\em Hermitian-Einstein metrics on parabolic stable bundles},
Acta Math. Sin. (Engl. Ser.) {\bf 15} (1999), no. 1, 93--114.

\bibitem{Lubke1}
M. L\"{u}bke,
{\em Chernklassen von Hermite-Einstein-Vektorb\"{u}ndeln.}
Math. Ann. {\bf 260} (1982), 133--141.

\bibitem{Lubke2}
M. L\"{u}bke, 
{\em Stability of Einstein-Hermitian vector bundles.}
Manuscripta Math. {\bf 42} (1983), 245--257.


\bibitem{Lubke-Teleman}
M. L\"{u}bke,
A. Teleman,
{\em The universal Kobayashi-Hitchin correspondence 
on Hermitian manifolds},
Mem. Amer. Math. Soc. {\bf 183} (2006), no. 863, vi+97 pp.

\bibitem{mehta-seshadri}
V. B. Mehta and C. S. Seshadri,
{\em Moduli of vector bundles on curves
with parabolic structures},
Math. Ann. {\bf 248}, (1980),
205--239.

\bibitem{Mochizuki-KHI}
T. Mochizuki,
{\em Kobayashi-Hitchin correspondence for tame harmonic bundles 
and an application},
Ast\'{e}risque {\bf 309} (2006), viii+117.

\bibitem{mochi2}
T. Mochizuki,
{\em Asymptotic behaviour of tame harmonic bundles
and an application to pure twistor $D$-modules I, II},
Mem. AMS. {\bf 185} (2007),
xii+565 pp.


\bibitem{Mochizuki-KHII}
T. Mochizuki,
{\em Kobayashi-Hitchin correspondence 
for tame harmonic bundles. II},
Geom. Topol. {\bf 13} (2009), 359--455. 

\bibitem{Mochizuki-wild}
T. Mochizuki,
{\em Wild harmonic bundles and wild pure twistor $D$-modules},
Ast\'{e}risque {\bf 340},
Soci\'{e}t\'{e} Math\'{e}matique de France, Paris,
2011.

\bibitem{Mochizuki-doubly-periodic}
T. Mochizuki,
{\em Asymptotic behaviour and 
the Nahm transform of doubly periodic instantons
with square integrable curvature},
Geom. Topol. {\bf 18}, (2014), 2823--2949. 

\bibitem{Mochizuki-KH-periodic}
T. Mochizuki,
{\em Periodic monopoles and difference modules},
arXiv:1712.08981.

\bibitem{Munteanu-Sesum}
O. Munteanu, N. Sesum,
{\em The Poisson equation on complete manifolds 
with positive spectrum and applications},
Adv. Math. {\bf 223} (2010), 198--219.

\bibitem{narasimhan-seshadri}
M. S. Narasimhan and C. S.  Seshadri,
{\em Stable and unitary vector bundles 
on a compact Riemann surface}, 
Ann. of Math. (2) {\bf 82} (1965), 540--567.


\bibitem{Ni1}
L. Ni,
{\em The Poisson equation and Hermitian-Einstein metrics 
on holomorphic vector bundles over complete noncompact 
K\"ahler manifolds},
Indiana Univ. Math. J. {\bf 51} (2002), 679--704.

\bibitem{Ni-Ren}
L. Ni, H. Ren, 
{\em Hermitian-Einstein metrics for vector bundles 
on complete K\"{a}hler manifolds}, 
Trans. Amer. Math. Soc. {\bf 353} (2001), 441--456.

\bibitem{Owens}
B. Owens,
{\em Instantons on cylindrical manifolds and stable bundles},
Geom. Topol. {\bf 5} (2001), 761--797.

\bibitem{Popovici}
D. Popovici,
{\em A simple proof of a theorem by Uhlenbeck and Yau},
Math. Z. {\bf 250} (2005), 855--872.



\bibitem{sabbah3}
C. Sabbah,
{\em Harmonic metrics and connections
with irregular singularities},
Ann. Inst. Fourier (Grenoble)
{\bf 49} (1999),
1265--1291.

\bibitem{Simpson88}
C. T. Simpson,
{\em Constructing variations of Hodge structure 
using Yang-Mills theory and applications to uniformization},
J. Amer. Math. Soc. {\bf 1} (1988), 867--918.

\bibitem{Simpson90}
C. T. Simpson,
{\em Harmonic bundles on noncompact curves},
J. Amer. Math. Soc. {\bf 3} (1990), 713--770.

\bibitem{steer-wren}
B. Steer and A. Wren,
{\em The Donaldson-Hitchin-Kobayashi correspondence
 for parabolic bundles over orbifold surfaces}, 
Canad. J. Math. {\bf 53} (2001), 1309--1339.

\bibitem{Takemoto1}
F. Takemoto,
{\em Stable vector bundles on algebraic surfaces},
Nagoya Math. J. {\bf 47} (1972), 29--48.

\bibitem{Takemoto2}
F. Takemoto,
{\em Stable vector bundles on algebraic surfaces. II.}
Nagoya Math. J. {\bf 52} (1973), 173--195.

\bibitem{Uhlenbeck-Yau}
K. Uhlenbeck,
S.-T. Yau,
{\em On the existence of Hermitian-Yang-Mills connections 
in stable vector bundles},
Comm. Pure Appl. Math. {\bf 39} (1986), no. S, suppl., S257--S293.

\bibitem{Yoshino-Master-Thesis}
M. Yoshino,
{\em The Nahm transform of spatially periodic instantons,}
arXiv:1804.05565.

\end{thebibliography}
\end{document}